
\input amstex
\input amsppt.sty

\magnification1200
\hsize391.46176pt
\vsize536.00175pt

\def\ZeilBP{18}
\def\ZeilBL{17}
\def\ZeilAS{16}
\def\StemAE{15}
\def\StanAP{14}
\def\SlatAC{13}
\def\SinVAA{12}
\def\RaVeAA{11}
\def\MiRRAC{10}
\def\KratBG{9}
\def\JordAA{8}
\def\GrKPAA{7}
\def\GeViAB{6}
\def\GaRaAA{5}
\def\DoraAA{4}
\def\AnStAA{3}
\def\AnBuAA{2}
\def\AndrAW{1}

\comment
\catcode`\@=11

\def\remarkheadfont@{\smc}
\let\varindent@\indent
\def\definition{\let\savedef@\definition \let\definition\relax
  \add@missing\endproclaim \add@missing\endroster
  \add@missing\enddefinition \envir@stack\enddefinition
   \def\definition##1{\restoredef@\definition
     \penaltyandskip@{-100}\medskipamount
        {\def\usualspace{{\remarkheadfont@\enspace}}%
        \varindent@\remarkheadfont@\ignorespaces##1\unskip
        \frills@{.\remarkheadfont@\enspace}}%
        \rm \ignorespaces}%
  \nofrillscheck\definition}
\catcode`\@=13
\endcomment

\catcode`\@=11
\font\tenln    = line10
\font\tenlnw   = linew10

\newskip\Einheit \Einheit=0.5cm
\newcount\xcoord \newcount\ycoord
\newdimen\xdim \newdimen\ydim \newdimen\PfadD@cke \newdimen\Pfadd@cke

\newcount\@tempcnta
\newcount\@tempcntb

\newdimen\@tempdima
\newdimen\@tempdimb

\newdimen\@wholewidth
\newdimen\@halfwidth

\newcount\@xarg
\newcount\@yarg
\newcount\@yyarg
\newbox\@linechar
\newbox\@tempboxa
\newdimen\@linelen
\newdimen\@clnwd
\newdimen\@clnht

\newif\if@negarg

\def\@whilenoop#1{}
\def\@whiledim#1\do #2{\ifdim #1\relax#2\@iwhiledim{#1\relax#2}\fi}
\def\@iwhiledim#1{\ifdim #1\let\@nextwhile=\@iwhiledim
        \else\let\@nextwhile=\@whilenoop\fi\@nextwhile{#1}}

\def\@whileswnoop#1\fi{}
\def\@whilesw#1\fi#2{#1#2\@iwhilesw{#1#2}\fi\fi}
\def\@iwhilesw#1\fi{#1\let\@nextwhile=\@iwhilesw
         \else\let\@nextwhile=\@whileswnoop\fi\@nextwhile{#1}\fi}

\def\thinlines{\let\@linefnt\tenln \let\@circlefnt\tencirc
  \@wholewidth\fontdimen8\tenln \@halfwidth .5\@wholewidth}
\def\thicklines{\let\@linefnt\tenlnw \let\@circlefnt\tencircw
  \@wholewidth\fontdimen8\tenlnw \@halfwidth .5\@wholewidth}
\thinlines

\PfadD@cke1pt \Pfadd@cke0.5pt
\def\PfadDicke#1{\PfadD@cke#1 \divide\PfadD@cke by2 \Pfadd@cke\PfadD@cke \multiply\PfadD@cke by2}
\long\def\LOOP#1\REPEAT{\def\BODY{#1}\ITERATE}
\def\ITERATE{\BODY \let\next\ITERATE \else\let\next\relax\fi \next}
\let\REPEAT=\fi
\def\Punkt{\hbox{\raise-2pt\hbox to0pt{\hss$\ssize\bullet$\hss}}}
\def\DuennPunkt(#1,#2){\unskip
  \raise#2 \Einheit\hbox to0pt{\hskip#1 \Einheit
          \raise-2.5pt\hbox to0pt{\hss$\bullet$\hss}\hss}}
\def\NormalPunkt(#1,#2){\unskip
  \raise#2 \Einheit\hbox to0pt{\hskip#1 \Einheit
          \raise-3pt\hbox to0pt{\hss\twelvepoint$\bullet$\hss}\hss}}
\def\DickPunkt(#1,#2){\unskip
  \raise#2 \Einheit\hbox to0pt{\hskip#1 \Einheit
          \raise-4pt\hbox to0pt{\hss\fourteenpoint$\bullet$\hss}\hss}}
\def\Kreis(#1,#2){\unskip
  \raise#2 \Einheit\hbox to0pt{\hskip#1 \Einheit
          \raise-4pt\hbox to0pt{\hss\fourteenpoint$\circ$\hss}\hss}}

\def\Line@(#1,#2)#3{\@xarg #1\relax \@yarg #2\relax
\@linelen=#3\Einheit
\ifnum\@xarg =0 \@vline
  \else \ifnum\@yarg =0 \@hline \else \@sline\fi
\fi}

\def\@sline{\ifnum\@xarg< 0 \@negargtrue \@xarg -\@xarg \@yyarg -\@yarg
  \else \@negargfalse \@yyarg \@yarg \fi
\ifnum \@yyarg >0 \@tempcnta\@yyarg \else \@tempcnta -\@yyarg \fi
\ifnum\@tempcnta>6 \@badlinearg\@tempcnta0 \fi
\ifnum\@xarg>6 \@badlinearg\@xarg 1 \fi
\setbox\@linechar\hbox{\@linefnt\@getlinechar(\@xarg,\@yyarg)}%
\ifnum \@yarg >0 \let\@upordown\raise \@clnht\z@
   \else\let\@upordown\lower \@clnht \ht\@linechar\fi
\@clnwd=\wd\@linechar
\if@negarg \hskip -\wd\@linechar \def\@tempa{\hskip -2\wd\@linechar}\else
     \let\@tempa\relax \fi
\@whiledim \@clnwd <\@linelen \do
  {\@upordown\@clnht\copy\@linechar
   \@tempa
   \advance\@clnht \ht\@linechar
   \advance\@clnwd \wd\@linechar}%
\advance\@clnht -\ht\@linechar
\advance\@clnwd -\wd\@linechar
\@tempdima\@linelen\advance\@tempdima -\@clnwd
\@tempdimb\@tempdima\advance\@tempdimb -\wd\@linechar
\if@negarg \hskip -\@tempdimb \else \hskip \@tempdimb \fi
\multiply\@tempdima \@m
\@tempcnta \@tempdima \@tempdima \wd\@linechar \divide\@tempcnta \@tempdima
\@tempdima \ht\@linechar \multiply\@tempdima \@tempcnta
\divide\@tempdima \@m
\advance\@clnht \@tempdima
\ifdim \@linelen <\wd\@linechar
   \hskip \wd\@linechar
  \else\@upordown\@clnht\copy\@linechar\fi}

\def\@hline{\ifnum \@xarg <0 \hskip -\@linelen \fi
\vrule height\Pfadd@cke width \@linelen depth\Pfadd@cke
\ifnum \@xarg <0 \hskip -\@linelen \fi}

\def\@getlinechar(#1,#2){\@tempcnta#1\relax\multiply\@tempcnta 8
\advance\@tempcnta -9 \ifnum #2>0 \advance\@tempcnta #2\relax\else
\advance\@tempcnta -#2\relax\advance\@tempcnta 64 \fi
\char\@tempcnta}

\def\Vektor(#1,#2)#3(#4,#5){\unskip\leavevmode
  \xcoord#4\relax \ycoord#5\relax
      \raise\ycoord \Einheit\hbox to0pt{\hskip\xcoord \Einheit
         \Vector@(#1,#2){#3}\hss}}

\def\Vector@(#1,#2)#3{\@xarg #1\relax \@yarg #2\relax
\@tempcnta \ifnum\@xarg<0 -\@xarg\else\@xarg\fi
\ifnum\@tempcnta<5\relax
\@linelen=#3\Einheit
\ifnum\@xarg =0 \@vvector
  \else \ifnum\@yarg =0 \@hvector \else \@svector\fi
\fi
\else\@badlinearg\fi}

\def\@hvector{\@hline\hbox to 0pt{\@linefnt
\ifnum \@xarg <0 \@getlarrow(1,0)\hss\else
    \hss\@getrarrow(1,0)\fi}}

\def\@vvector{\ifnum \@yarg <0 \@downvector \else \@upvector \fi}

\def\@svector{\@sline
\@tempcnta\@yarg \ifnum\@tempcnta <0 \@tempcnta=-\@tempcnta\fi
\ifnum\@tempcnta <5
  \hskip -\wd\@linechar
  \@upordown\@clnht \hbox{\@linefnt  \if@negarg
  \@getlarrow(\@xarg,\@yyarg) \else \@getrarrow(\@xarg,\@yyarg) \fi}%
\else\@badlinearg\fi}

\def\@upline{\hbox to \z@{\hskip -.5\Pfadd@cke \vrule width \Pfadd@cke
   height \@linelen depth \z@\hss}}

\def\@downline{\hbox to \z@{\hskip -.5\Pfadd@cke \vrule width \Pfadd@cke
   height \z@ depth \@linelen \hss}}

\def\@upvector{\@upline\setbox\@tempboxa\hbox{\@linefnt\char'66}\raise
     \@linelen \hbox to\z@{\lower \ht\@tempboxa\box\@tempboxa\hss}}

\def\@downvector{\@downline\lower \@linelen
      \hbox to \z@{\@linefnt\char'77\hss}}

\def\@getlarrow(#1,#2){\ifnum #2 =\z@ \@tempcnta='33\else
\@tempcnta=#1\relax\multiply\@tempcnta \sixt@@n \advance\@tempcnta
-9 \@tempcntb=#2\relax\multiply\@tempcntb \tw@
\ifnum \@tempcntb >0 \advance\@tempcnta \@tempcntb\relax
\else\advance\@tempcnta -\@tempcntb\advance\@tempcnta 64
\fi\fi\char\@tempcnta}

\def\@getrarrow(#1,#2){\@tempcntb=#2\relax
\ifnum\@tempcntb < 0 \@tempcntb=-\@tempcntb\relax\fi
\ifcase \@tempcntb\relax \@tempcnta='55 \or
\ifnum #1<3 \@tempcnta=#1\relax\multiply\@tempcnta
24 \advance\@tempcnta -6 \else \ifnum #1=3 \@tempcnta=49
\else\@tempcnta=58 \fi\fi\or
\ifnum #1<3 \@tempcnta=#1\relax\multiply\@tempcnta
24 \advance\@tempcnta -3 \else \@tempcnta=51\fi\or
\@tempcnta=#1\relax\multiply\@tempcnta
\sixt@@n \advance\@tempcnta -\tw@ \else
\@tempcnta=#1\relax\multiply\@tempcnta
\sixt@@n \advance\@tempcnta 7 \fi\ifnum #2<0 \advance\@tempcnta 64 \fi
\char\@tempcnta}

\def\Diagonale(#1,#2)#3{\unskip\leavevmode
  \xcoord#1\relax \ycoord#2\relax
      \raise\ycoord \Einheit\hbox to0pt{\hskip\xcoord \Einheit
         \Line@(1,1){#3}\hss}}
\def\AntiDiagonale(#1,#2)#3{\unskip\leavevmode
  \xcoord#1\relax \ycoord#2\relax 
      \raise\ycoord \Einheit\hbox to0pt{\hskip\xcoord \Einheit
         \Line@(1,-1){#3}\hss}}
\def\Pfad(#1,#2),#3\endPfad{\unskip\leavevmode
  \xcoord#1 \ycoord#2 \thicklines\ZeichnePfad#3\endPfad\thinlines}
\def\ZeichnePfad#1{\ifx#1\endPfad\let\next\relax
  \else\let\next\ZeichnePfad
    \ifnum#1=1
      \raise\ycoord \Einheit\hbox to0pt{\hskip\xcoord \Einheit
         \vrule height\Pfadd@cke width1 \Einheit depth\Pfadd@cke\hss}%
      \advance\xcoord by 1
    \else\ifnum#1=2
      \raise\ycoord \Einheit\hbox to0pt{\hskip\xcoord \Einheit
        \hbox{\hskip-\PfadD@cke\vrule height1 \Einheit width\PfadD@cke depth0pt}\hss}%
      \advance\ycoord by 1
    \else\ifnum#1=3
      \raise\ycoord \Einheit\hbox to0pt{\hskip\xcoord \Einheit
         \Line@(1,1){1}\hss}
      \advance\xcoord by 1
      \advance\ycoord by 1
    \else\ifnum#1=4
      \raise\ycoord \Einheit\hbox to0pt{\hskip\xcoord \Einheit
         \Line@(1,-1){1}\hss}
      \advance\xcoord by 1
      \advance\ycoord by -1
    \fi\fi\fi\fi
  \fi\next}
\def\hSSchritt{\leavevmode\raise-.4pt\hbox to0pt{\hss.\hss}\hskip.2\Einheit
  \raise-.4pt\hbox to0pt{\hss.\hss}\hskip.2\Einheit
  \raise-.4pt\hbox to0pt{\hss.\hss}\hskip.2\Einheit
  \raise-.4pt\hbox to0pt{\hss.\hss}\hskip.2\Einheit
  \raise-.4pt\hbox to0pt{\hss.\hss}\hskip.2\Einheit}
\def\vSSchritt{\vbox{\baselineskip.2\Einheit\lineskiplimit0pt
\hbox{.}\hbox{.}\hbox{.}\hbox{.}\hbox{.}}}
\def\DSSchritt{\leavevmode\raise-.4pt\hbox to0pt{%
  \hbox to0pt{\hss.\hss}\hskip.2\Einheit
  \raise.2\Einheit\hbox to0pt{\hss.\hss}\hskip.2\Einheit
  \raise.4\Einheit\hbox to0pt{\hss.\hss}\hskip.2\Einheit
  \raise.6\Einheit\hbox to0pt{\hss.\hss}\hskip.2\Einheit
  \raise.8\Einheit\hbox to0pt{\hss.\hss}\hss}}
\def\dSSchritt{\leavevmode\raise-.4pt\hbox to0pt{%
  \hbox to0pt{\hss.\hss}\hskip.2\Einheit
  \raise-.2\Einheit\hbox to0pt{\hss.\hss}\hskip.2\Einheit
  \raise-.4\Einheit\hbox to0pt{\hss.\hss}\hskip.2\Einheit
  \raise-.6\Einheit\hbox to0pt{\hss.\hss}\hskip.2\Einheit
  \raise-.8\Einheit\hbox to0pt{\hss.\hss}\hss}}
\def\SPfad(#1,#2),#3\endSPfad{\unskip\leavevmode
  \xcoord#1 \ycoord#2 \ZeichneSPfad#3\endSPfad}
\def\ZeichneSPfad#1{\ifx#1\endSPfad\let\next\relax
  \else\let\next\ZeichneSPfad
    \ifnum#1=1
      \raise\ycoord \Einheit\hbox to0pt{\hskip\xcoord \Einheit
         \hSSchritt\hss}%
      \advance\xcoord by 1
    \else\ifnum#1=2
      \raise\ycoord \Einheit\hbox to0pt{\hskip\xcoord \Einheit
        \hbox{\hskip-2pt \vSSchritt}\hss}%
      \advance\ycoord by 1
    \else\ifnum#1=3
      \raise\ycoord \Einheit\hbox to0pt{\hskip\xcoord \Einheit
         \DSSchritt\hss}
      \advance\xcoord by 1
      \advance\ycoord by 1
    \else\ifnum#1=4
      \raise\ycoord \Einheit\hbox to0pt{\hskip\xcoord \Einheit
         \dSSchritt\hss}
      \advance\xcoord by 1
      \advance\ycoord by -1
    \fi\fi\fi\fi
  \fi\next}
\def\Koordinatenachsen(#1,#2){\unskip
 \hbox to0pt{\hskip-.5pt\vrule height#2 \Einheit width.5pt depth1 \Einheit}%
 \hbox to0pt{\hskip-1 \Einheit \xcoord#1 \advance\xcoord by1
    \vrule height0.25pt width\xcoord \Einheit depth0.25pt\hss}}
\def\Koordinatenachsen(#1,#2)(#3,#4){\unskip
 \hbox to0pt{\hskip-.5pt \ycoord-#4 \advance\ycoord by1
    \vrule height#2 \Einheit width.5pt depth\ycoord \Einheit}%
 \hbox to0pt{\hskip-1 \Einheit \hskip#3\Einheit 
    \xcoord#1 \advance\xcoord by1 \advance\xcoord by-#3 
    \vrule height0.25pt width\xcoord \Einheit depth0.25pt\hss}}
\def\Gitter(#1,#2){\unskip \xcoord0 \ycoord0 \leavevmode
  \LOOP\ifnum\ycoord<#2
    \loop\ifnum\xcoord<#1
      \raise\ycoord \Einheit\hbox to0pt{\hskip\xcoord \Einheit\Punkt\hss}%
      \advance\xcoord by1
    \repeat
    \xcoord0
    \advance\ycoord by1
  \REPEAT}
\def\Gitter(#1,#2)(#3,#4){\unskip \xcoord#3 \ycoord#4 \leavevmode
  \LOOP\ifnum\ycoord<#2
    \loop\ifnum\xcoord<#1
      \raise\ycoord \Einheit\hbox to0pt{\hskip\xcoord \Einheit\Punkt\hss}%
      \advance\xcoord by1
    \repeat
    \xcoord#3
    \advance\ycoord by1
  \REPEAT}
\def\Label#1#2(#3,#4){\unskip \xdim#3 \Einheit \ydim#4 \Einheit
  \def\lo{\advance\xdim by-.5 \Einheit \advance\ydim by.5 \Einheit}%
  \def\llo{\advance\xdim by-.25cm \advance\ydim by.5 \Einheit}%
  \def\loo{\advance\xdim by-.5 \Einheit \advance\ydim by.25cm}%
  \def\o{\advance\ydim by.25cm}%
  \def\ro{\advance\xdim by.5 \Einheit \advance\ydim by.5 \Einheit}%
  \def\rro{\advance\xdim by.25cm \advance\ydim by.5 \Einheit}%
  \def\roo{\advance\xdim by.5 \Einheit \advance\ydim by.25cm}%
  \def\l{\advance\xdim by-.30cm}%
  \def\r{\advance\xdim by.30cm}%
  \def\lu{\advance\xdim by-.5 \Einheit \advance\ydim by-.6 \Einheit}%
  \def\llu{\advance\xdim by-.25cm \advance\ydim by-.6 \Einheit}%
  \def\luu{\advance\xdim by-.5 \Einheit \advance\ydim by-.30cm}%
  \def\u{\advance\ydim by-.30cm}%
  \def\ru{\advance\xdim by.5 \Einheit \advance\ydim by-.6 \Einheit}%
  \def\rru{\advance\xdim by.25cm \advance\ydim by-.6 \Einheit}%
  \def\ruu{\advance\xdim by.5 \Einheit \advance\ydim by-.30cm}%
  #1\raise\ydim\hbox to0pt{\hskip\xdim
     \vbox to0pt{\vss\hbox to0pt{\hss$#2$\hss}\vss}\hss}%
}
\catcode`\@=13

\def\LL{\leavevmode\setbox0=\hbox{L}\hbox to\wd0{\hss\char'40L}}
\def\al{\alpha}
\def\be{\beta}
\def\ga{\gamma}

\def\ep{\varepsilon}

\def\la{\lambda}

\def\si{\sigma}
\def\ta{\tau}


\def\today{\ifcase\month\or
 January\or February\or March\or April\or May\or June\or
 July\or August\or September\or October\or November\or December\fi
 \space\number\day, \number\year}

\def\({\left(}
\def\){\right)}
\def\[{\left[}
\def\]{\right]}

\def\sgn{\operatorname{sgn}}

\def\suml{\sum\limits}

\def\prodl{\prod\limits}

\def\3{\ss}
\catcode`\@=11
\def\dddot#1{\vbox{\ialign{##\crcr
      .\hskip-.5pt.\hskip-.5pt.\crcr\noalign{\kern1.5\p@\nointerlineskip}
      $\hfil\displaystyle{#1}\hfil$\crcr}}}

\newif\iftab@\tab@false
\newif\ifvtab@\vtab@false
\def\tab{\bgroup\tab@true\vtab@false\vst@bfalse\Strich@false%
   \def\\{\global\hline@@false%
     \ifhline@\global\hline@false\global\hline@@true\fi\cr}
   \edef\l@{\the\leftskip}\ialign\bgroup\hskip\l@##\hfil&&##\hfil\cr}
\def\endtab{\cr\egroup\egroup}
\def\vtab{\vtop\bgroup\vst@bfalse\vtab@true\tab@true\Strich@false%
   \bgroup\def\\{\cr}\ialign\bgroup&##\hfil\cr}
\def\endvtab{\cr\egroup\egroup\egroup}
\def\stab{\D@cke0.5pt\null 
 \bgroup\tab@true\vtab@false\vst@bfalse\Strich@true\Let@@\vspace@
 \normalbaselines\offinterlineskip
  \openup\spreadmlines@
 \edef\l@{\the\leftskip}\ialign
 \bgroup\hskip\l@##\hfil&&##\hfil\crcr}
\def\endstab{\crcr\egroup
 \egroup}
\newif\ifvst@b\vst@bfalse
\def\vstab{\D@cke0.5pt\null
 \vtop\bgroup\tab@true\vtab@false\vst@btrue\Strich@true\bgroup\Let@@\vspace@
 \normalbaselines\offinterlineskip
  \openup\spreadmlines@\bgroup}
\def\endvstab{\crcr\egroup\egroup
 \egroup\tab@false\Strich@false}

\newdimen\htstrut@
\htstrut@8.5\p@
\newdimen\htStrut@
\htStrut@12\p@
\newdimen\dpstrut@
\dpstrut@3.5\p@
\newdimen\dpStrut@
\dpStrut@3.5\p@
\def\openup{\afterassignment\@penup\dimen@=}
\def\@penup{\advance\lineskip\dimen@
  \advance\baselineskip\dimen@
  \advance\lineskiplimit\dimen@
  \divide\dimen@ by2
  \advance\htstrut@\dimen@
  \advance\htStrut@\dimen@
  \advance\dpstrut@\dimen@
  \advance\dpStrut@\dimen@}
\def\Let@@{\relax%
    \def\\{\global\hline@@false%
     \ifhline@\global\hline@false\global\hline@@true\fi\cr}%
    \iffalse}\fi}
\def\matrix{\null\,\vcenter\bgroup
 \tab@false\vtab@false\vst@bfalse\Strich@false\Let@@\vspace@
 \normalbaselines\openup\spreadmlines@\ialign
 \bgroup\hfil$\m@th##$\hfil&&\quad\hfil$\m@th##$\hfil\crcr
 \Mathstrut@\crcr\noalign{\kern-\baselineskip}}
\def\endmatrix{\crcr\Mathstrut@\crcr\noalign{\kern-\baselineskip}\egroup
 \egroup\,}
\def\smatrix{\D@cke0.5pt\null\,
 \vcenter\bgroup\tab@false\vtab@false\vst@bfalse\Strich@true\Let@@\vspace@
 \normalbaselines\offinterlineskip
  \openup\spreadmlines@\ialign
 \bgroup\hfil$\m@th##$\hfil&&\quad\hfil$\m@th##$\hfil\crcr}
\def\endsmatrix{\crcr\egroup
 \egroup\,\Strich@false}
\newdimen\D@cke
\def\Dicke#1{\global\D@cke#1}
\newtoks\tabs@\tabs@{&}
\newif\ifStrich@\Strich@false
\newif\iff@rst

\def\Stricherr@{\iftab@\ifvtab@\errmessage{\noexpand\s not allowed
     here. Use \noexpand\vstab!}%
  \else\errmessage{\noexpand\s not allowed here. Use \noexpand\stab!}%
  \fi\else\errmessage{\noexpand\s not allowed
     here. Use \noexpand\smatrix!}\fi}
\def\format{\ifvst@b\else\crcr\fi\egroup\iffalse{\fi\ifnum`}=0 \fi\format@}
\def\format@#1\\{\def\preamble@{#1}%
 \def\Str@chfehlt##1{\ifx##1\s\Stricherr@\fi\ifx##1\\\let\Next\relax%
   \else\let\Next\Str@chfehlt\fi\Next}%
 \def\c{\hfil\noexpand\ifhline@@\hbox{\vrule height\htStrut@%
   depth\dpstrut@ width\z@}\noexpand\fi%
   \ifStrich@\hbox{\vrule height\htstrut@ depth\dpstrut@ width\z@}%
   \fi\iftab@\else$\m@th\fi\the\hashtoks@\iftab@\else$\fi\hfil}%
 \def\r{\hfil\noexpand\ifhline@@\hbox{\vrule height\htStrut@%
   depth\dpstrut@ width\z@}\noexpand\fi%
   \ifStrich@\hbox{\vrule height\htstrut@ depth\dpstrut@ width\z@}%
   \fi\iftab@\else$\m@th\fi\the\hashtoks@\iftab@\else$\fi}%
 \def\l{\noexpand\ifhline@@\hbox{\vrule height\htStrut@%
   depth\dpstrut@ width\z@}\noexpand\fi%
   \ifStrich@\hbox{\vrule height\htstrut@ depth\dpstrut@ width\z@}%
   \fi\iftab@\else$\m@th\fi\the\hashtoks@\iftab@\else$\fi\hfil}%
 \def\s{\ifStrich@\ \the\tabs@\vrule width\D@cke\the\hashtoks@%
          \fi\the\tabs@\ }%
 \def\sa{\ifStrich@\vrule width\D@cke\the\hashtoks@%
            \the\tabs@\ %
            \fi}%
 \def\se{\ifStrich@\ \the\tabs@\vrule width\D@cke\the\hashtoks@\fi}%
 \def\cd{\hfil\noexpand\ifhline@@\hbox{\vrule height\htStrut@%
   depth\dpstrut@ width\z@}\noexpand\fi%
   \ifStrich@\hbox{\vrule height\htstrut@ depth\dpstrut@ width\z@}%
   \fi$\dsize\m@th\the\hashtoks@$\hfil}%
 \def\rd{\hfil\noexpand\ifhline@@\hbox{\vrule height\htStrut@%
   depth\dpstrut@ width\z@}\noexpand\fi%
   \ifStrich@\hbox{\vrule height\htstrut@ depth\dpstrut@ width\z@}%
   \fi$\dsize\m@th\the\hashtoks@$}%
 \def\ld{\noexpand\ifhline@@\hbox{\vrule height\htStrut@%
   depth\dpstrut@ width\z@}\noexpand\fi%
   \ifStrich@\hbox{\vrule height\htstrut@ depth\dpstrut@ width\z@}%
   \fi$\dsize\m@th\the\hashtoks@$\hfil}%
 \ifStrich@\else\Str@chfehlt#1\\\fi%
 \setbox\z@\hbox{\xdef\Preamble@{\preamble@}}\ifnum`{=0 \fi\iffalse}\fi
 \ialign\bgroup\span\Preamble@\crcr}
\newif\ifhline@\hline@false
\newif\ifhline@@\hline@@false
\def\hlinefor#1{\multispan@{\strip@#1 }\leaders\hrule height\D@cke\hfill%
    \global\hline@true\ignorespaces}
\def\Item "#1"{\par\noindent\hangindent2\parindent%
  \hangafter1\setbox0\hbox{\rm#1\enspace}\ifdim\wd0>2\parindent%
  \box0\else\hbox to 2\parindent{\rm#1\hfil}\fi\ignorespaces}
\def\ITEM #1"#2"{\par\noindent\hangafter1\hangindent#1%
  \setbox0\hbox{\rm#2\enspace}\ifdim\wd0>#1%
  \box0\else\hbox to 0pt{\rm#2\hss}\hskip#1\fi\ignorespaces}
\def\item"#1"{\par\noindent\hang%
  \setbox0=\hbox{\rm#1\enspace}\ifdim\wd0>\the\parindent%
  \box0\else\hbox to \parindent{\rm#1\hfil}\enspace\fi\ignorespaces}
\let\plainitem@\item

\font@\twelverm=cmr10 scaled\magstep1
\font@\twelveit=cmti10 scaled\magstep1
\font@\twelvebf=cmbx10 scaled\magstep1
\font@\twelvei=cmmi10 scaled\magstep1
\font@\twelvesy=cmsy10 scaled\magstep1
\font@\twelveex=cmex10 scaled\magstep1

\newtoks\twelvepoint@
\def\twelvepoint{\normalbaselineskip15\p@
 \abovedisplayskip15\p@ plus3.6\p@ minus10.8\p@
 \belowdisplayskip\abovedisplayskip
 \abovedisplayshortskip\z@ plus3.6\p@
 \belowdisplayshortskip8.4\p@ plus3.6\p@ minus4.8\p@
 \textonlyfont@\rm\twelverm \textonlyfont@\it\twelveit
 \textonlyfont@\sl\twelvesl \textonlyfont@\bf\twelvebf
 \textonlyfont@\smc\twelvesmc \textonlyfont@\tt\twelvett
%
 \ifsyntax@ \def\big##1{{\hbox{$\left##1\right.$}}}%
  \let\Big\big \let\bigg\big \let\Bigg\big
 \else
  \textfont\z@=\twelverm  \scriptfont\z@=\tenrm  \scriptscriptfont\z@=\sevenrm
  \textfont\@ne=\twelvei  \scriptfont\@ne=\teni  \scriptscriptfont\@ne=\seveni
  \textfont\tw@=\twelvesy \scriptfont\tw@=\tensy \scriptscriptfont\tw@=\sevensy
  \textfont\thr@@=\twelveex \scriptfont\thr@@=\tenex
        \scriptscriptfont\thr@@=\tenex
  \textfont\itfam=\twelveit \scriptfont\itfam=\tenit
        \scriptscriptfont\itfam=\tenit
  \textfont\bffam=\twelvebf \scriptfont\bffam=\tenbf
        \scriptscriptfont\bffam=\sevenbf
  \setbox\strutbox\hbox{\vrule height10.2\p@ depth4.2\p@ width\z@}%
  \setbox\strutbox@\hbox{\lower.6\normallineskiplimit\vbox{%
        \kern-\normallineskiplimit\copy\strutbox}}%
 \setbox\z@\vbox{\hbox{$($}\kern\z@}\bigsize@=1.4\ht\z@
 \fi
 \normalbaselines\rm\ex@.2326ex\jot3.6\ex@\the\twelvepoint@}

\font@\fourteenrm=cmr10 scaled\magstep2
\font@\fourteenit=cmti10 scaled\magstep2
\font@\fourteensl=cmsl10 scaled\magstep2
\font@\fourteensmc=cmcsc10 scaled\magstep2
\font@\fourteentt=cmtt10 scaled\magstep2
\font@\fourteenbf=cmbx10 scaled\magstep2
\font@\fourteeni=cmmi10 scaled\magstep2
\font@\fourteensy=cmsy10 scaled\magstep2
\font@\fourteenex=cmex10 scaled\magstep2
\font@\fourteenmsa=msam10 scaled\magstep2
\font@\fourteeneufm=eufm10 scaled\magstep2
\font@\fourteenmsb=msbm10 scaled\magstep2
\newtoks\fourteenpoint@
\def\fourteenpoint{\normalbaselineskip15\p@
 \abovedisplayskip18\p@ plus4.3\p@ minus12.9\p@
 \belowdisplayskip\abovedisplayskip
 \abovedisplayshortskip\z@ plus4.3\p@
 \belowdisplayshortskip10.1\p@ plus4.3\p@ minus5.8\p@
 \textonlyfont@\rm\fourteenrm \textonlyfont@\it\fourteenit
 \textonlyfont@\sl\fourteensl \textonlyfont@\bf\fourteenbf
 \textonlyfont@\smc\fourteensmc \textonlyfont@\tt\fourteentt
%
 \ifsyntax@ \def\big##1{{\hbox{$\left##1\right.$}}}%
  \let\Big\big \let\bigg\big \let\Bigg\big
 \else
  \textfont\z@=\fourteenrm  \scriptfont\z@=\twelverm  \scriptscriptfont\z@=\tenrm
  \textfont\@ne=\fourteeni  \scriptfont\@ne=\twelvei  \scriptscriptfont\@ne=\teni
  \textfont\tw@=\fourteensy \scriptfont\tw@=\twelvesy \scriptscriptfont\tw@=\tensy
  \textfont\thr@@=\fourteenex \scriptfont\thr@@=\twelveex
        \scriptscriptfont\thr@@=\twelveex
  \textfont\itfam=\fourteenit \scriptfont\itfam=\twelveit
        \scriptscriptfont\itfam=\twelveit
  \textfont\bffam=\fourteenbf \scriptfont\bffam=\twelvebf
        \scriptscriptfont\bffam=\tenbf
  \setbox\strutbox\hbox{\vrule height12.2\p@ depth5\p@ width\z@}%
  \setbox\strutbox@\hbox{\lower.72\normallineskiplimit\vbox{%
        \kern-\normallineskiplimit\copy\strutbox}}%
 \setbox\z@\vbox{\hbox{$($}\kern\z@}\bigsize@=1.7\ht\z@
 \fi
 \normalbaselines\rm\ex@.2326ex\jot4.3\ex@\the\fourteenpoint@}

\catcode`\@=13

\magnification1200

\hsize13cm
\vsize19cm

\TagsOnRight

\catcode`\@=11
\def\iddots{\mathinner{\mkern1mu\raise\p@\hbox{.}\mkern2mu
    \raise4\p@\hbox{.}\mkern2mu\raise7\p@\vbox{\kern7\p@\hbox{.}}\mkern1mu}}
\catcode`\@=13
\def\CT{\operatorname{CT}}
\def\pf{\operatornamewithlimits{Pf}}
\def\SA{S}
\def\SB{T}
\def\SC{U}

\def\DetA{D}
\def\DetAa{D_{A}}
\def\DetAaa{\overline D_{A}}
\def\DetAb{D_{B}}
\def\DetAbb{\overline D_{B}}
\def\DetB{E}
\def\DetBa{E_{A}}
\def\DetBaa{\overline E_{A}}
\def\DetBb{E_{B}}
\def\DetBt{\tilde E_{B}}
\def\DetBbb{\overline E_{B}}
\def\bP{{\overline P}}
\def\fl#1{\left\lfloor#1\right\rfloor}
\def\cl#1{\left\lceil#1\right\rceil}
\def\y{{\bar y}}
\def\z{{\bar z}}
\def\q{{\tilde q}}
\def\po#1#2{(#1)_#2}

\topmatter 
\title Determinant identities and a generalization of the number of
totally symmetric self-complementary plane partitions
\endtitle 
\author C.~Krattenthaler\footnote"$^\dagger$"{Supported in part by EC's Human
Capital and Mobility Program, grant CHRX-CT93-0400 and the\linebreak 
\hbox{Austrian Science Foundation FWF, grant P10191-MAT}}
\endauthor 
\affil 
Institut f\"ur Mathematik der Universit\"at Wien,\\
Strudlhofgasse 4, A-1090 Wien, Austria.\\
e-mail: KRATT\@Pap.Univie.Ac.At\\
WWW: \tt http://radon.mat.univie.ac.at/People/kratt
\endaffil
\address Institut f\"ur Mathematik der Universit\"at Wien,
Strudlhofgasse 4, A-1090 Wien, Austria.
\endaddress
\dedicatory Submitted: September 16, 1997; Accepted: November 3, 1997\enddedicatory
\abstract We prove a constant term conjecture of Robbins and Zeilberger 
(J.~Combin\. Theory Ser.~A {\bf 66} (1994), 17--27),
by translating the problem into a determinant evaluation problem and
evaluating the determinant. This determinant generalizes the
determinant that gives the number of all totally symmetric
self-complementary plane partitions contained in a
$(2n)\times(2n)\times(2n)$ box and that was used by Andrews 
(J.~Combin\. Theory Ser.~A {\bf 66} (1994), 28--39) and Andrews and
Burge (Pacific J. Math\. {\bf 158} (1993), 1--14) to compute this
number explicitly. The evaluation of the generalized determinant is
independent of Andrews and Burge's computations, and therefore in
particular constitutes a new solution to this famous enumeration
problem. We also evaluate a related determinant, thus
generalizing another determinant identity of Andrews and Burge (loc\.
cit\.). By translating some of our determinant identities into constant term
identities, we obtain several new constant term identities.
\endabstract
\endtopmatter

\font\smcp=cmcsc8
\headline={\ifnum\pageno>1 {\smcp the electronic journal of  
combinatorics 4 (1997), \#Rxx\hfill\folio} \fi}

\document

\leftheadtext{C. Krattenthaler}
\rightheadtext{Determinant identities}

\subhead 1. Introduction\endsubhead
I started work on this paper originally hoping to find a proof
of the following conjecture of Robbins and Zeilberger
\cite{\ZeilAS, Conjecture {\bf C'=B'}} (caution: in the quotient
defining $B'$ it should read $(m+1+2j)$ instead of $(m+1+j)$), which
we state in an equivalent form.
\proclaim{Conjecture}Let $x$ and $n$ be nonnegative integers. Then
$$\multline \CT \(\frac {\prod _{0\le i<j\le n-1} ^{}(1-z_i/z_j)\prod
_{i=0} ^{n-1}(1+z_i^{-1})^{x+n-i-1}}
{\prod _{0\le i<j\le n-1} ^{}(1-z_iz_j)\prod _{i=0} ^{n-1}(1-z_i)}\)\\
=\cases \prodl _{i=0} ^{n-1}\frac {(3x+3i+1)!} {(3x+2i+1)!\,(x+2i)!}
\prodl _{i=0} ^{(n-2)/2}(2x+2i+1)!\,(2i)!&\text {if $n$ is even}
\hbox to1.5cm{\hskip1cm\hskip.56pt \rm(1.1a)\hss}\\
2^x\prodl _{i=1} ^{n-1}\frac {(3x+3i+1)!} {(3x+2i+1)!\,(x+2i)!}
\prodl _{i=1} ^{(n-1)/2}(2x+2i)!\,(2i-1)!&\text {if $n$ is odd.}
\hbox to1.5cm{\hskip1cm \rm(1.1b)\hss}
\endcases
\endmultline$$
Here, {\rm CT(Expr)} means the constant term in {\rm Expr}, i.e., the coefficient
of $z_1^0z_2^0\cdots z_n^0$ in {\rm Expr}.
\endproclaim

I thought this might be a rather boring task since 
in the case $x=0$ there existed already a proof of the
Conjecture (see \cite{\ZeilAS}). 
This proof consists of translating the constant term on
the left-hand side of (1.1) into a sum of minors of a particular
matrix (by a result \cite{\ZeilAS, Corollary {\bf D=C}} of
Zeilberger), which is known to equal the number of totally symmetric
self-complementary plane partitions contained in a
$(2n)\times(2n)\times(2n)$ box (by a result of Doran \cite{\DoraAA,
Theorem~4.1 + Proof~2 of Theorem~5.1}). 
The number of these plane partitions had been calculated by Andrews
\cite{\AndrAW} by transforming the sum of minors into a single
determinant (using a result of Stembridge \cite{\StemAE, Theorem~3.1,
Theorem~8.3}) and evaluating the determinant. Since Zeilberger shows
in \cite{\ZeilAS, Lemma {\bf D'=C'}} that the translation of the
constant term in (1.1) into a sum of minors of some matrix works for
generic $x$, and since Stembridge's result \cite{\StemAE, Theorem~3.1}
still applies to obtain a single determinant (see (2.2)), my idea was
to evaluate this determinant by
routinely extending Andrews's proof of the totally symmetric
self-complementary plane partitions conjecture, or the alternative
proofs by Andrews and Burge \cite{\AnBuAA}. However, it became clear
rather quickly that this is not possible (at least not {\it
routinely}). In fact, the aforementioned proofs take advantage of a
few remarkable coincidences, which break down if $x$ is nonzero.
Therefore I had to devise new methods and tools to solve the
determinant problem in this more general case where $x\neq 0$. 

In the course of the work on the problem, the subject became more and
more exciting as I came across an increasing number of interesting
determinants that could be evaluated, thus generalizing several
determinant identities of Andrews and Burge \cite{\AnBuAA},
which appeared in connection
with the enumeration of totally symmetric self-complementary
plane partitions. In the end, I had found a proof of the Conjecture,
but also many more interesting results. In this paper, I describe this
proof and all further results.

The proof of the Conjecture will be organized as follows. In Theorem~1,
item (3) in
Section~1 it is proved that the constant term in (1.1) equals the
positive square root of a certain determinant, actually of one
determinant, namely (2.2a), if $n$ is even, and of another
determinant, namely (2.2b), if $n$ is odd. 
In addition, Theorem~1 provides two more equivalent interpretations
of the constant term, in particular a combinatorial interpretation in terms of
shifted plane partitions, which reduces to totally symmetric
self-complementary plane partitions for $x=0$. 

The main idea that we will use to evaluate the determinants
in Theorem~1 will be to generalize them by
introducing a further parameter, $y$, see (3.1) and (4.1). 
The generalized determinants reduce to
the determinants of Theorem~1 when $y=x$. Many of our arguments do not
work without this generalization. In Section~3 we study the
two-parameter family (3.1) of determinants that contains (2.2a) as
special case. If $y=x+m$, with $m$ a fixed integer, Theorem~2 makes
it possible to evaluate the resulting determinants. This is done for
a few cases in Corollary~3, including the case $y=x$ (see (3.69)) that
we are particularly interested in. Similarly, in Section~4 we study the
two-parameter family (4.1) of determinants that contains (2.2b) as
special case. Also here,
if $y=x+m$, with $m$ a fixed integer, Proposition~5 makes
it possible to evaluate the resulting determinants. 
This is done for
two cases in Corollary~6, including the case $y=x$ (see (4.42)) that
we are particularly interested in. 
This concludes the proof of the Conjecture, which thus becomes a
theorem. It is restated as such in Theorem~11.
However, even more is possible for this second family of determinants.
In Theorem~8,
we succeed in evaluating the determinants (4.1) for {\it independent}
$x$ and $y$, taking advantage of all previous results in Section~4.

There is another
interesting determinant identity, which is related to the
aforementioned determinant identities. This is the subject of
Section~5. It generalizes a 
determinant identity of Andrews and Burge \cite{\AnBuAA}.
Finally, in Section~6 we translate
our determinant identities of Sections~4 and 5 into constant
term identities which seem to be new. Auxiliary results that are
needed in the proofs of our Theorems are collected in the Appendix.

Since a first version of this article was written, $q$-analogues of
two of the determinant evaluations in this article, Theorems~8 and
10, were found in \cite{\KratBG}. No $q$-analogues are known for the
results in Section~3. Also, it is still open whether the $q$-analogues
of \cite{\KratBG} have any combinatorial meaning. Another interesting
development is that Amdeberhan (private communication) observed that
Dodgson's determinant formula (see \cite{\ZeilBP, \ZeilBL}) can be
used to give a short inductive proof of the determinant evaluation in
Theorem~10 (and also of its $q$-analogue in \cite{\KratBG}), and
could also be used to give an inductive proof of the determinant
evaluation in Theorem~8 (and its $q$-analogue in \cite{\KratBG}) 
provided one is able to prove a certain
identity featuring three double summations.

\subhead 2. Transformation of the Conjecture into a determinant
evaluation problem\endsubhead
In Theorem~1 below we show that the constant term in (1.1) equals the
positive square root of some determinant, one if $n$ is even, another
if $n$ is odd. Also, we provide a combinatorial interpretation of the
constant term in terms of shifted plane partitions. Recall that a
shifted plane partition of shape $(\la_1,\la_2,\dots,\la_r)$ is an
array $\pi$ of integers of the form
$$\matrix 
\pi_{1,1}&\pi_{1,2}&\pi_{1,3}&\innerhdotsfor5\after\quad &\pi_{1,\la_1}\\
&\pi_{2,2}&\pi_{2,3}&\innerhdotsfor4\after\quad &\pi_{2,\la_2}\\
&&\ddots&&\vdots&&\iddots\\
&&&\pi_{r,r}&\hdots&\pi_{r,\la_r}
\endmatrix$$
such that the rows and columns are weakly decreasing. Curiously
enough, we need this combinatorial interpretation to know that we
have to choose the {\it positive\/} root once the determinant is
evaluated.
\proclaim{Theorem 1}Let $x$ and $n$ be nonnegative integers. The
constant term in {\rm(1.1)} equals
\roster
\item "\rm(1)" the sum of all $n\times n$ minors of the $n\times
(2n-1)$ matrix 
$$\(\binom {x+i}{j-i}\)_{0\le i\le n-1,\ 0\le j\le
2n+x-2},\tag2.1$$
\item "\rm(2)" the number of shifted plane partitions of shape
$(x+n-1,x+n-2,\dots,1)$, with entries between 0 and $n$, where the
entries in row $i$ are at least
$n-i$, $i=1,2,\dots,n-1$,
\item "(3)" the positive square root of
$$\hskip1cm\cases \det\limits_{0\le i,j\le n-1}\bigg(\sum\limits
 _{x+2i-j<r\le x+2j-i}
^{}\binom {2x+i+j}r\bigg)&\text {if $n$ is even,}
\hbox to1.4cm{\hskip.6cm \rm(2.2a)\hss}\\
2^{2x}\det\limits_{0\le i,j\le n-2}\(\frac
{(2x+i+j+1)!\,(3x+3i+4)(3x+3j+4)(3j-3i)}
{(x+2i-j+2)!\,(x+2j-i+2)!}\)&\text {if $n$ is odd,}
\hbox to1.4cm{\hskip.6cm\hskip2.6pt \rm(2.2b)\hss}
\endcases$$
\endroster
if the sums in {\rm(2.2a)} are interpreted by
$$\suml _{r=A+1} ^{B}\text {\rm Expr}(r)=\cases \hphantom{-}
\suml _{r=A+1} ^{B} \text {\rm Expr}(r)&A<B\\
\hphantom{-}0&A=B\\
-\suml _{r=B+1} ^{A}\text {\rm Expr}(r)&A>B.\endcases\tag2.3$$
\endproclaim

\demo{Proof} {\it ad {\rm (1)}}. This was proved by Zeilberger
\cite{\ZeilAS, Lemma {\bf D'=C'}}. (Note that we have performed a
shift of the indices $i,j$ in comparison with Zeilberger's notation.)

\smallskip
{\it ad {\rm (2)}}. Fix a minor of the matrix (2.1),
$$\det_{0\le i,j\le n-1}\(\binom {x+i}{\la_j-i}\)$$
say. By the main theorem of nonintersecting lattice paths 
\cite{\GeViAB, Cor.~2; \StemAE, Theorem~1.2} (see Proposition~A1) 
this determinant has an
interpretation in terms of nonintersecting lattice paths. By a
lattice path we mean a lattice path in the plane 
consisting of unit horizontal and vertical steps
in the positive direction. Furthermore, recall that a family of paths
is called nonintersecting if no two paths of the family have a point
in common. Now, the above determinant equals the number of all
families $(P_0,P_1,\dots,P_{n-1})$ of nonintersecting lattice paths,
where $P_i$ runs from $(-2i,i)$ to $(x-\la_i,\la_i)$,
$i=0,1,\dots,n-1$. An example with $x=2$, $n=5$, $\la_1=1$,
$\la_2=3$, $\la_3=4$, $\la_4=7$, $\la_5=9$ is displayed in
Figure~1.a. (Ignore $P_{-1}$ for the moment.)
\midinsert
\vskip10pt
\vbox{
$$
\raise0cm\hbox{$\underset\text {\raise-10pt\hbox{a. nonintersecting
lattice paths}}\to
{\Einheit0.45cm
\hbox{\hskip4.5cm}
\Gitter(5,11)(-8,-2)
\Koordinatenachsen(5,11)(-8,-2)
\thinlines
\Vektor(1,0){0}(5,0)
\Vektor(0,1){0}(0,11)
\AntiDiagonale(-9,11){14}
\thicklines
\Pfad(0,0),12\endPfad
\Pfad(-2,1),221\endPfad
\Pfad(-4,2),2121\endPfad
\Pfad(-6,3),22221\endPfad
\Pfad(-8,4),221222\endPfad
\DickPunkt(4,-2)
\DickPunkt(0,0)
\DickPunkt(-2,1)
\DickPunkt(-4,2)
\DickPunkt(-6,3)
\DickPunkt(-8,4)
\DickPunkt(1,1)
\DickPunkt(-1,3)
\DickPunkt(-2,4)
\DickPunkt(-5,7)
\DickPunkt(-7,9)
\Label\r{\eightpoint x_1}(5,0)
\Label\o{\eightpoint x_2}(0,11)
\Label\ro{P_{-1}}(4,-2)
\Label\u{P_{0}}(1,0)
\Label\r{P_{1}}(-2,2)
\Label\u{P_{2}}(-3,3)
\Label\r{P_{3}}(-6,5)
\Label\r{P_{4}}(-8,5)
\Label\o{\eightpoint x_1+x_2=x}(-9,11)
\hskip2cm}$}
\raise2cm\hbox{$\longrightarrow$}
\raise1cm\hbox{$\underset\text {\raise-10pt\hbox{b. noncrossing
lattice paths}}\to
{\hbox{\hskip3cm}
\Einheit0.5cm
\Gitter(3,7)(-4,-1)
\Koordinatenachsen(3,7)(-4,-1)
\thinlines
\Vektor(1,0){0}(3,0)
\Vektor(0,1){0}(0,7)
\AntiDiagonale(-5,7){8}
\thicklines
\Pfad(0,0),12\endPfad
\hbox to1pt{\hss}
\raise-1pt\hbox{
\Pfad(-1,0),221\endPfad}
\hbox to-2pt{\hss}
\raise1pt\hbox{
\Pfad(-2,0),2121\endPfad}
\hbox to2pt{\hss}
\raise-1pt\hbox{
\Pfad(-3,0),22221\endPfad}
\hbox to-2pt{\hss}
\raise1pt\hbox{
\Pfad(-4,0),221222\endPfad}
\DickPunkt(0,0)
\DickPunkt(-1,0)
\DickPunkt(-2,0)
\DickPunkt(-3,0)
\DickPunkt(-4,0)
\DickPunkt(1,1)
\DickPunkt(0,2)
\DickPunkt(-2,4)
\DickPunkt(-3,5)
\Label\r{\eightpoint x_1}(3,0)
\Label\o{\eightpoint x_2}(0,7)
\Label\u{P'_{0}}(1,0)
\Label\r{P'_{1}}(-1,1)
\Label\lo{P'_{2}}(-1,1)
\Label\ru{P'_{3}}(-3,4)
\Label\l{P'_{4}}(-4,2)
\Label\o{\eightpoint x_1+x_2=x}(-5,7)
\hskip2cm}$}
$$
\vskip-1cm
$$
\raise1.5cm\hbox{$\underset\text {\raise-10pt\hbox{d. shifted plane
partition}}\to
{\matrix 5&5&5&5&4&4\\
&4&3&3&3&3\\
&&3&3&3&2\\
&&&3&1&1\\
&&&&1&1\\
&&&&&0\endmatrix}$}
\raise2cm\hbox{\hskip.8cm$\longleftarrow$}
\raise0cm\hbox{$\underset\text {\raise-10pt\hbox{c. filling of the
regions}}\to
{\hbox{\hskip3.2cm}
\Einheit0.5cm
\Vektor(0,-1){2}(1,10)
\Gitter(3,7)(-4,-1)
\Koordinatenachsen(3,7)(-4,-1)
\thinlines
\Vektor(1,0){0}(3,0)
\Vektor(0,1){0}(0,7)
\AntiDiagonale(-5,7){8}
\thicklines
\Pfad(0,0),12\endPfad
\hbox to1pt{\hss}
\raise-1pt\hbox{
\Pfad(-1,0),221\endPfad}
\hbox to-2pt{\hss}
\raise1pt\hbox{
\Pfad(-2,0),2121\endPfad}
\hbox to2pt{\hss}
\raise-1pt\hbox{
\Pfad(-3,0),22221\endPfad}
\hbox to-2pt{\hss}
\raise1pt\hbox{
\Pfad(-4,0),221222\endPfad}
\PfadDicke{.5pt}
\Pfad(-4,-1),22222222\endPfad
\Label\ro{0}(1,0)
\Label\ro{1}(0,0)
\Label\ro{1}(0,1)
\Label\ro{1}(-1,0)
\Label\ro{1}(-1,1)
\Label\ro{2}(-2,0)
\Label\ro{3}(-1,2)
\Label\ro{3}(-2,1)
\Label\ro{3}(-2,2)
\Label\ro{3}(-2,3)
\Label\ro{3}(-3,0)
\Label\ro{3}(-3,1)
\Label\ro{3}(-3,2)
\Label\ro{3}(-3,3)
\Label\ro{4}(-4,0)
\Label\ro{4}(-4,1)
\Label\ro{4}(-3,4)
\Label\ro{5}(-4,2)
\Label\ro{5}(-4,3)
\Label\ro{5}(-4,4)
\Label\ro{5}(-4,5)
\hskip2cm}$}
$$
\centerline{\eightpoint Figure 1}
}
\vskip10pt
\endinsert
Hence, we see that the sum of all minors of the matrix (2.1) equals
the number of all families $(P_0,P_1,\dots,P_{n-1})$ of
nonintersecting lattice paths, where $P_i$ runs from $(-2i,i)$
to {\it some\/} point on the antidiagonal line $x_1+x_2=x$ ($x_1$
denoting the horizontal coordinate, $x_2$ denoting the vertical
coordinate), $i=0,1,\dots,n-1$. 
Next, given such a family $(P_0,P_1,\dots,P_{n-1})$ of
nonintersecting lattice paths, we shift $P_i$ by the vector $(i,-i)$,
$i=0,1,\dots,n-1$. Thus a family $(P'_0,P'_1,\dots,P'_{n-1})$ of
lattice paths is obtained, where $P'_i$ runs from $(-i,0)$ to some
point on the line $x_1+x_2=x$, see Figure~1.b. The new paths may
touch each other, but they cannot cross each other. Therefore, the paths
$P'_0,P'_1,\dots,P'_{n-1}$ cut the triangle that
is bordered by the $x_1$-axes, the line $x_1+x_2=x$,
the vertical line $x_1=-n+1$ into exactly $n+1$ regions. We fill
these regions with integers as is exemplified in Figure~1.c. To be more
precise, the region to the right of $P'_0$ is filled with 0's, the
region between $P'_0$ and $P'_1$ is filled with 1's, \dots, the
region between $P'_{n-2}$ and $P'_{n-1}$ is filled with $(n-1)$'s,
and the region to the left of $P'_{n-1}$ is filled with $n$'s.
Finally, we forget about the paths and reflect the array of integers
just obtained in an antidiagonal line, see Figure~1.d. 
Clearly, a shifted plane
partition of shape $(x+n-1,x+n-2,\dots,1)$ is obtained. Moreover, it
has the desired property that the entries in row $i$ are at least
$n-i$, $i=1,2,\dots,n-1$. It is
easy to see that each step can be reversed, which completes the proof
of (2).

\smallskip
{\it ad \rm(3)}. It was proved just before 
that the constant term in (1.1) equals the
number of all families $(P_0,P_1,\dots,P_{n-1})$ of
nonintersecting lattice paths, where $P_i$ runs from $(-2i,i)$
to some point on the antidiagonal line $x_1+x_2=x$, $i=0,1,\dots,n-1$.

Now, let first $n$ be even. By a theorem of Stembridge \cite{\StemAE,
Theorem~3.1} (see Proposition~A2, 
with $A_i=(-2i,i)$, $i=0,1,\dots,n-1$,
$I={}$(the lattice points on the line $x_1+x_2=x$)), 
the number of such families of
nonintersecting lattice paths equals the Pfaffian
$$\pf_{0\le i<j\le n-1}\big(Q(i,j)\big),\tag2.4$$
where $Q(i,j)$ is the number of all pairs $(P_i,P_j)$ of
nonintersecting lattice paths, $P_i$ running from $(-2i,i)$ to some
point on the line $x_1+x_2=x$, and $P_j$ running from $(-2j,j)$ to some
point on the line $x_1+x_2=x$. 

In order to compute the number $Q(i,j)$ for fixed $i,j$, $0\le i<j\le n-1$, 
we follow Stembridge's computation in the proof of
Theorem~8.3 in \cite{\StemAE}.
We define $b_{kl}$ to be the number of all pairs $(P_i,P_j)$ of {\it
intersecting\/} lattice paths, where $P_i$ runs from $(-2i,i)$ to
$(x-k,k)$, and where $P_j$ runs from $(-2j,j)$ to $(x-l,l)$. Since
the total number of lattice paths from $(-2i,i)$ to $x_1+x_2=x$ is
$2^{x+i}$, it follows that $2^{2x+i+j}-Q(i,j)$ is the number of pairs
of intersecting lattice paths from $(-2i,i)$ and $(-2j,j)$ to
$x_1+x_2=x$. Hence, 
$$2^{2x+i+j}-Q(i,j)=\sum _{k,l} ^{}b_{kl}=\sum _{k<l} ^{}b_{lk}+\sum
_{k\ge l} ^{}b_{lk},$$
the last equality being a consequence of the fact that
$b_{kl}=b_{lk}$, which is proved by the standard path switching
argument (find the first meeting point and interchange terminal
portions of the paths from thereon, see the proofs of 
\cite{\GeViAB, Cor.~2; \StemAE, Theorem~1.2}). When $k\le l$, every
path from $(-2i,i)$ to $(x-l,l)$ must intersect every path from
$(-2j,j)$ to $(x-k,k)$, so we have $b_{kl}=\binom{x+i}{l-i}
\binom{x+j}{k-j}$. Thus,
$$\multline
2^{2x+i+j}-Q(i,j)=\sum _{0\le k<l\le x+2i-j} ^{}\binom{x+i}{l+j-i}
\binom{x+j}{k}\\
+\sum _{0\le k\le l\le x+2i-j} ^{}\binom{x+i}{l+j-i}
\binom{x+j}{k}.
\endmultline$$
Now we replace $l$ by $x+2i-j-l$ in the first sum and $k$ by
$x+2i-j-k$ in the second sum. This leads to
$$\multline
2^{2x+i+j}-Q(i,j)=\sum _{k+l< x+2i-j} ^{}\binom{x+i}{l}
\binom{x+j}{k}\\
+\sum _{k+l\ge x+2i-j} ^{}\binom{x+i}{l+j-i}
\binom{x+j}{k+2j-2i}.
\endmultline$$
For fixed values of $r=k+l$ both sums can be simplified further by
the Vandermonde sum (see e.g\. \cite{\GrKPAA, sec.~5.1, (5.27)}), so
$$2^{2x+i+j}-Q(i,j)=\sum _{r< x+2i-j} ^{}\binom{2x+i+j}{r}
+\sum _{r\ge x+2i-j} ^{}\binom{2x+i+j}{r+3j-3i},
$$
and finally, after replacement of $r$ by $2x+i+j-r$ in the first sum,
and by $2x+4i-2j-r$ in the second sum,
$$\align
Q(i,j)&=2^{2x+i+j}-\sum _{r> x+2j-i} ^{}\binom{2x+i+j}{r}
-\sum _{r\le x+2i-j} ^{}\binom{2x+i+j}{r}\\
&=\sum _{x+2i-j<r\le x+2j-i} ^{}\binom{2x+i+j}{r}.\tag2.5
\endalign$$

As is well-known, the square of a Pfaffian equals the determinant of
the corresponding skew-symmetric matrix (see e.g\. \cite{\StemAE,
Prop.~2.2}). The quantity $Q(i,j)$, as given by (2.5), has the
property $Q(i,j)=-Q(j,i)$, due to our interpretation (2.3) of limits
of sums.
Hence, the square of the Pfaffian in (2.4) equals 
$\det_{0\le i<j\le n-1}\big(Q(i,j)\big)$, which in view of (2.5) is
exactly (2.2a).
That the Pfaffian itself is the {\it positive\/} square root of the
determinant is due to the combinatorial interpretation in item (2) of
the Theorem. Thus, item (3) is established for even $n$.

Now let $n$ be odd. Still, by the proof of (2), the constant term in
(1.1) equals the
number of all families $(P_0,P_1,\dots,P_{n-1})$ of
nonintersecting lattice paths, where $P_i$ runs from $(-2i,i)$
to some point on the antidiagonal line $x_1+x_2=x$, $i=0,1,\dots,n-1$.
However, to apply Theorem~3.1 of \cite{\StemAE} again we have to add
a ``dummy path" $P_{-1}$ of length 0, running from $(2x,-x)$ to
$(2x,-x)$, say (cf\. the Remark after Theorem~3.1 in \cite{\StemAE};
however, we order all paths {\it after\/} the dummy path). 
See Figure~1.a for the location of $P_{-1}$. We infer
that the constant term in (1.1) equals
$$\pf_{-1\le i<j\le n-1}\big(Q(i,j)\big),\tag2.6$$
where $Q(i,j)$ is the number of all pairs $(P_i,P_j)$ of
nonintersecting lattice paths, $P_i$ running from $(-2i,i)$ to 
the line $x_1+x_2=x$ if $i\ge0$, $P_{-1}$ running from
$(2x,-x)$ to $x_1+x_2=x$ (hence, to $(2x,-x)$), 
and $P_j$ running from $(-2j,j)$ to the line $x_1+x_2=x$. If $0\le
i<j\le n-1$, then 
$Q(i,j)=\sum _{x+2i-j<r\le x+2j-i} ^{}\binom{2x+i+j}{r}$ according to
the computation that led to (2.5). Moreover, we have
$Q(-1,j)=2^{x+j}$ since a pair $(P_{-1},P_j)$ is nonintersecting for
any path $P_j$ running from $(-2j,j)$ to $x_1+x_2=x$. The latter fact
is due to the location of $P_{-1}$, see Figure~1.a. Therefore, the
square of the Pfaffian in (2.6) equals
$$\det_{-1\le i,j\le n-1}\(
\SPfad(0,0),111111111111\endSPfad
\SPfad(3,-2),2222\endSPfad
\Label\r{0}(1,1)
\Label\r{-2^{x+i}}(1,-1)
\Label\l{2^{x+j}}(7,1)
\Label\l{\tsize \sum\limits _{x+2i-j<r\le x+2j-i}
^{}\binom{2x+i+j}{r}}(8,-1)
\Label\r{\hskip8pt i=-1}(13,1)
\Label\r{i\ge0}(13,-1)
\hskip5.9cm
\).
\Label\u{\hskip5pt j=-1}(-11,-2)
\Label\u{j\ge0}(-5,-2)
\hskip2cm
\tag2.7
$$
\vskip3pt

We subtract 2 times the $(j-1)$-st column from the $j$-th column,
$j=n-1,n-2,\dots,2$, in this order, and 
we subtract 2 times the $(i-1)$-st row from the $i$-th row,
$i=n-1,n-2,\dots,2$. Thus, by simple algebra, the determinant in (2.7)
is turned into
$$
\det_{-1\le i,j\le n-1}\(
\SPfad(0,2),111111111111111\endSPfad
\SPfad(0,1),111111111111111\endSPfad
\SPfad(2,-2),22222\endSPfad
\SPfad(4,-2),22222\endSPfad
\Label\o{0}(1,2)
\Label\o{0}(5,2)
\Label\o{\hbox to 3.7cm{\dotfill\hskip10pt}}(10,2)
\Label\o{0}(14,2)
\Label\o{0}(1,0)
\Label\o{\vdots}(1,-1)
\Label\o{0}(1,-2)
\Label\o{-2^{x}}(1,1)
\Label\o{2^{x}}(3,2)
\Label\o{*}(3,1)
\Label\o{*}(5,1)
\Label\o{\hbox to3.7cm{\dotfill\hskip10pt}}(10,1)
\Label\o{*}(14,1)
\Label\o{*}(3,0)
\Label\o{\vdots}(3,-1)
\Label\o{*}(3,-2)
\Label\o{\tsize \frac {(2x+i+j-1)!\,(3x+3i+1)(3x+3j+1)(3j-3i)}
{(x+2i-j+1)!\,(x+2j-i+1)!}\hskip6pt}(10,-1)
\Label\ro{\hskip8pt i\!=\!-1}(17,2)
\Label\ro{i\!=\!0}(17,1)
\Label\ro{i\!\ge\!1}(17,-1)
\hskip7.6cm
\).
\Label\r{\hskip-18pt j\!=\!-1}(-15,-3)
\Label\r{\hskip-12pt j\!=\!0}(-13,-3)
\Label\r{j\!\ge\!1}(-7,-3)
\hskip2cm
\tag2.8
$$
\vskip3pt
\noindent
By expanding this determinant along the top row, and the resulting
determinant along the left-most column, we arrive at (2.2b), upon
rescaling row and column indices.

\smallskip
Thus the proof of Theorem~1 is complete.
\quad \quad \qed
\enddemo

\remark{Remark} Mills, Robbins and Rumsey \cite{\MiRRAC, Theorem~1 +
last paragraph of p.~281}
showed that shifted plane partitions of shape $(n-1,n-2,\dots,1)$, 
where the entries in row $i$ are at least
$n-i$ and at most $n$, $i=1,2,\dots,n-1$, are in bijection with totally symmetric
self-complementary plane partitions contained in a
$(2n)\times(2n)\times(2n)$ box. Hence, by item (2) of Theorem~1, the
number (1.1) generalizes the number of these plane partitions, to
which it reduces for $x=0$.

The idea that is used in the translation of item (1) into item (2) of Theorem~1 is 
due to Doran \cite{\DoraAA, Proof of Theorem~4.1}, 
who did this translation for $x=0$. However, our presentation is
modelled after Stembridge's presentation of Doran's idea in
\cite{\StemAE, Proof of Theorem~8.3}.
 
\endremark

\subhead 3. A two-parameter family of determinants\endsubhead
The goal of this section is to evaluate the determinant in (2.2a). We
shall even consider the generalized determinant
$$\DetA(x,y;n):=\det_{0\le i,j\le n-1}\bigg(\sum
 _{x+2i-j<r\le y+2j-i}
^{}\binom {x+y+i+j}r\bigg),\tag3.1$$
for integral $x$ and $y$, which reduces to (2.2a) when $y=x$. 
In fact, many of our arguments essentially require this
generalization and would not work without it. Recall
that the sums in (3.1) have to be interpreted according to (2.3).

The main result of this section, Theorem~2 below, allows to evaluate
$\DetA(x,y;n)$ when the difference $m=y-x$ is fixed. It is done
explicitly for a number of cases in the subsequent Corollary~3,
including the case $m=0$ which gives the evaluation of (2.2a) that we
are particularly interested in. For the sake of
brevity, Theorem~2 is formulated only for $y\ge x$ (i.e., for
$m\ge0$). The corresponding result for $y\le x$ is easily obtained by
taking advantage of the fact
$$\DetA(x,y;n)=(-1)^{n}\DetA(y,x;n),\tag3.2$$
which results from transposing the matrix in (3.1) and using (2.3).
\proclaim{Theorem 2}Let $x,m,n$ be nonnegative integers with $m\le
n$. Then, with the usual notation $(a)_k:=a(a+1)\cdots(a+k-1)$,
$k\ge1$, $(a)_0:=1$,
of shifted factorials, there holds
$$\multline \DetA(x,x+m;n)=\det_{0\le i,j\le n-1}\bigg(\sum
 _{x+2i-j<r\le x+m+2j-i}
^{}\binom {2x+m+i+j}r\bigg)\\
=\prod _{i=1} ^{n-1}\(\frac {(2x+m+i)!\,(3x+m+2i+2)_i\,(3x+2m+2i+2)_i}
{(x+2i)!\,(x+m+2i)!}\)\hskip2.5cm\\
\times \frac {(2x+m)!} {(x+\fl{m/2})!\,(x+m)!}\cdot\prod _{i=0}
^{\fl{n/2}-1}(2x+2\cl{m/2}+2i+1)\cdot P_1(x;m,n),
\endmultline\tag3.3$$
where $P_1(x;m,n)$ is a polynomial in $x$ of degree $\le \fl{m/2}$.
If $n$ is odd and $m$ is even, the polynomial $P_1(x;m,n)$ is
divisible by $(2x+m+n)$. For fixed $m$, the polynomial $P_1(x;m,n)$ can
be computed explicitly by specializing $x$ to $-\fl{(m+n)/2}+t-1/2$,
$t=0,1,\dots,\fl{m/2}$, 
in the identity {\rm (3.67)}. This makes sense since
for these specializations the determinant in {\rm(3.67)} reduces to a
determinant of size at most $2t+1$, as is elaborated in Step~6 of the
proof, and since a polynomial of degree $\le \fl{m/2}$ is uniquely
determined by its values at $\fl{m/2}+1$ distinct points.
\endproclaim
\demo{Proof}The proof is divided into several steps. Our strategy is
to transform $\DetA(x,x+m;n)$ into a multiple of another determinant,
namely $\DetAb(x,x+m;n)$, by (3.6), (3.8) and (3.10), which is a polynomial in $x$, 
then identify as many 
factors of the new determinant  
as possible (as a polynomial in $x$), and finally find a bound
for the degree of the remaining polynomial factor.

For big parts of the proof we shall write $y$ for $x+m$. We feel
that this makes things more transparent.
\smallskip
{\it Step 1. Equivalent expressions for $\DetA(x,y;n)$.}
First, in the definition (3.1) of $\DetA(x,y;n)$
we subtract 2 times the $(j-1)$-st column from the $j$-th column,
$j=n-1,n-2,\dots,1$, in this order, and 
we subtract 2 times the $(i-1)$-st row from the $i$-th row,
$i=n-1,n-2,\dots,1$. By simple algebra we get
$$\multline
\DetA(x,y;n)\\
=\det_{0\le i,j\le n-1}\!\!\(
\hbox{\hskip14pt}
\raise4pt\hbox{
\SPfad(-1,0),11111111111111111\endSPfad
\SPfad(6,-2),2222\endSPfad
\Label\r{\tsize\sum\limits _{x<r\le y}
^{}\binom {x+y}r}(2,1)
\Label\r{\tsize \frac {(x+y+j)!\,(x+2y+3j+1)}
{(x-j+1)!\,(y+2j)!}}(10,1)
\Label\r{\tsize -\frac {(x+y+i)!\,(2x+y+3i+1)}
{(x+2i)!\,(y-i+1)!}\hskip5pt}(2,-1)
\Label\r{\hskip8pt\tsize \frac {(x+y+i+j-1)!\,(y-x+3j-3i)}
{(x+2i-j+1)!\,(y+2j-i+1)!}}(9,-1)
\Label\r{\hskip12pt\ssize \times(2x+y+3i+1)(x+2y+3j+1)}(11,-2)
\Label\r{\hskip-4pti\!=\!0}(18,1)
\Label\ro{i\!\ge\!1}(18,-2)
}
\hskip7.8cm
\).
\Label\r{\hskip-12pt j\!=\!0}(-15,-3)
\Label\r{j\!\ge\!1}(-7,-3)
\Label\r{\hskip12pt(3.4)}(1,-3)
\hskip1cm
\endmultline
$$
\vskip3pt

On the other hand, if in the definition (3.1) of $\DetA(x,y;n)$
we subtract $1/2$ times the $(j+1)$-st column from the $j$-th column,
$j=0,1\dots,n-2$, in this order, and if
we subtract $1/2$ times the $(i+1)$-st row from the $i$-th row,
$i=0,1,\dots,n-2$, we get
$$\multline
\DetA(x,y;n)\\
=\det_{0\le i,j\le n-1}\!\!\(
\raise0pt\hbox{
\hbox{\hskip4pt}
\SPfad(0,0),11111111111111111\endSPfad
\SPfad(8,-2),2222\endSPfad
\Label\ru{\hskip4pt\tsize \frac {1} {4}\frac {(x+y+i+j+1)!\,(y-x+3j-3i)}
{(x+2i-j+2)!\,(y+2j-i+2)!}}(3,2)
\Label\ru{\hskip0pt\ssize \times(2x+y+3i+4)(x+2y+3j+4)}(3,1)
\Label\r{\tsize \hskip0pt\frac {1} {2}\frac {(x+y+i+n)!\,(2x+y+3i+4)}
{(x+2i-n+3)!\,(y+2n-i-2)!}}(12,1)
\Label\r{\tsize -\frac {1} {2}\frac {(x+y+j+n)!\,(x+2y+3j+4)}
{(x+2n-j-2)!\,(y+2j-n+3)!}\hskip0pt}(3,-1)
\Label\r{\tsize\sum\limits _{x+n-1<r\le y+n-1}
^{}\binom {x+y+2n-2}r}(12,-1)
\Label\r{\hskip-4pt i\!\le\! n\!-\!2}(19,1)
\Label\r{\hskip-4pt i\!=\!n\!-\!1}(19,-1)
}
\hskip8.4cm
\).
\Label\r{\hskip-12pt j\!\le\! n\!-\!2}(-15,-3)
\Label\r{j\!=\!n\!-\!1}(-7,-3)
\Label\r{\hskip12pt(3.5)}(1,-3)
\hskip1cm
\endmultline
$$
\vskip3pt

\smallskip
{\it Step 2. An equivalent statement of the Theorem.}
We consider the expression (3.4). We take as many common factors out
of the $i$-th row, $i=0,1,\dots,n-1$, as possible, such that the
entries become polynomials in $x$ and $y$. To be precise, 
we take
$$\frac {(x+y+i)!\,(2x+y+3i+1)} {(x+2i)!\,(y+2n-i-1)!}$$
out of the $i$-th row, $i=1,2,\dots,n-1$, and we take 
$$\frac {(x+y)!} {\fl{(x+y)/2}!\,(y+2n-2)!}$$
out of the 0-th row. Furthermore, we take $(x+2y+3j+1)$ out of the
$j$-th column, $j=1,2,\dots,n-1$. This gives
$$\multline
\DetA(x,y;n)\\
=\frac {(x+y)!} {\fl{(x+y)/2}!\,(y+2n-2)!}
\prodl _{i=1} ^{n-1}\frac {(x+y+i)!\,(2x+y+3i+1)\,(x+2y+3i+1)} 
{(x+2i)!\,(y+2n-i-1)!}\\
\times\det_{0\le i,j\le n-1}\!\!\(
\hbox{\hskip0pt}
\raise0pt\hbox{
\SPfad(0,0),111111111111111\endSPfad
\SPfad(4,-2),2222\endSPfad
\Label\r{\hskip5pt\SA(x,y;n)}(1,1)
\Label\ru{\ssize (x+y+1)_j\,(x-j+2)_{\fl{(y-x)/2}+j-1}}(9,2)
\Label\ru{\raise6pt\hbox{$\ssize \times(y+2j+1)_{2n-2j-2}$}}(9,1)
\Label\r{\hskip5pt\ssize -(y-i+2)_{2n-2}}(1,-1)
\Label\ro{\raise-10pt\hbox{$\ssize (x+y+i+1)_{j-1}\,(x+2i-j+2)_{j-1}$}}(8,-1)
\Label\ro{\hskip0pt\ssize \times(y+2j-i+2)_{2n-2j-2}\,(y-x+3j-3i)}(10,-2)
\Label\r{\hskip-4pt i\!=\!0}(17,1)
\Label\r{\hskip-4pt i\!\ge\!1}(17,-1)
}
\hskip7.5cm
\),
\Label\r{\hskip-12pt j\!=\!0}(-15,-3)
\Label\r{j\!\ge\!1}(-7,-3)
\Label\r{\hskip12pt(3.6)}(1,-3)
\hskip1cm
\endmultline
$$
\vskip3pt
\noindent
where $\SA(x,y;n)$ is given by
$$\multline \SA(x,y;n)=\sum _{r=1}
^{\fl{(y-x)/2}}(x+r+1)_{\fl{(y-x)/2}-r}\, (y-r+1)_{2n+r-2}\\
+\sum _{r=0}
^{\cl{(y-x)/2}-1}(x+r+1)_{\fl{(y-x)/2}-r}\, (y-r+1)_{2n+r-2}.
\endmultline\tag3.7$$

For convenience, let us denote the determinant in (3.6) by
$\DetAa(x,y;n)$.
In fact, there are more factors that can be taken out of $\DetAa(x,y;n)$ 
under the
restriction that the entries of the determinant continue to be
polynomials. To this end, we multiply the $i$-th row of $\DetAa(x,y;n)$
by $(y+2n-i)_{i-1}$, $i=1,2,\dots,n-1$, divide the $j$-th column by
$(y+2j+1)_{2n-2j-2}$, $j=1,2,\dots,n-1$, 
and divide the $0$-th column by $(y+1)_{2n-2}$.
This leads to
$$\multline
\prodl _{i=1} ^{n-1}(y+2n-i)_{i-1}\prodl
_{j=1} ^{n-1}\frac {1} {(y+2j+1)_{2n-2j-2}}\,\cdot \frac {1}
{(y+1)_{2n-2}}\cdot \DetAa(x,y;n)
\\
=\det_{0\le i,j\le n-1}\!\!\(
\hbox{\hskip0pt}
\raise4pt\hbox{
\SPfad(0,0),11111111111111\endSPfad
\SPfad(4,-2),2222\endSPfad
\Label\r{\hskip5pt\SB(x,y)}(1,1)
\Label\r{\ssize (x+y+1)_j\,(x-j+2)_{\fl{(y-x)/2}+j-1}}(9,1)
\Label\r{\hskip5pt\ssize -(y-i+2)_{i-1}}(1,-1)
\Label\ro{\raise-10pt\hbox{$\ssize (x+y+i+1)_{j-1}\,(x+2i-j+2)_{j-1}$}}(8,-1)
\Label\ro{\hskip0pt\ssize \times(y+2j-i+2)_{i-1}\,(y-x+3j-3i)}(10,-2)
\Label\r{\hskip-4pt \!i\!=0}(17,1)
\Label\r{\hskip-4pt i\!\ge\!1}(17,-1)
}
\hskip7.2cm
\),
\Label\r{\hskip-12pt j\!=\!0}(-15,-3)
\Label\r{j\!\ge\!1}(-7,-3)
\Label\r{\hskip12pt(3.8)}(1,-3)
\hskip1cm
\endmultline
$$
\vskip3pt
\noindent
where $\SB(x,y)$ is given by
$$\multline \SB(x,y)=\sum _{r=1}
^{\fl{(y-x)/2}}(x+r+1)_{\fl{(y-x)/2}-r}\, (y-r+1)_{r}\\
+\sum _{r=0}
^{\cl{(y-x)/2}-1}(x+r+1)_{\fl{(y-x)/2}-r}\, (y-r+1)_{r},
\endmultline\tag3.9$$
or, if we denote the determinant in (3.8) by $\DetAb(x,y;n)$,
$$\DetAa(x,y;n)=(y+1)_{2n-2}\prod _{i=1} ^{n-1}(y+2i+1)_{n-i-1}
\cdot\DetAb(x,y;n).\tag3.10$$
A combination of (3.3), (3.8), and (3.10) then implies that Theorem~2
is equivalent to the statement: 

{\sl With $\DetAb(x,y;n)$ the determinant in {\rm (3.8)},
there holds
$$\multline
\DetAb(x,y;n)=\prod _{i=1}
^{n-1}\big((2x+y+2i+2)_{i-1}\,(x+2y+2i+2)_{i-1}\big)\\
\times\prod _{i=0}
^{\fl{n/2}-1}(2x+2\cl{(y-x)/2}+2i+1)\cdot P_1(x;y-x,n),
\endmultline\tag3.11$$
where $P_1(x;y-x,n)$
satisfies the properties that are stated in Theorem~2}.

Recall that $y=x+m$, where $m$ is a fixed nonnegative integer. In the
subsequent steps of the proof we are going to establish that (3.11)
does not hold only for integral $x$, but holds as a polynomial identity
in $x$. In order to accomplish this, we show in Step~3 that the first
product on the right-hand side of (3.11) is a factor of
$\DetAb(x,y;n)$, then we show in Step~4 that the second
product on the right-hand side of (3.11) is a factor of
$\DetAb(x,y;n)$, and finally we show in Step~5 that the degree of
$\DetAb(x,y;n)$ is at most $2\binom {n-1}2+\fl{n/2}+\fl{(y-x)/2}$, which
implies that the degree of $P_1(x;y-x,n)$ is at most
$\fl{(y-x)/2}=\fl{m/2}$. Once this is done, the proof of Theorem~2
will be complete (except for the statement about $P_1(x;y-x,n)$ for
odd $n$ and even $m$, which is proved in Step~4, and the algorithm
for computing $P_1(x;y-x,n)$ explicitly, which is described in
Step~6).

{\it Step 3. $\prod _{i=1}
^{n-1}\big((2x+y+2i+2)_{i-1}\,(x+2y+2i+2)_{i-1}\big)$
is a factor of $\DetAb(x,y;n)$.}
Here we consider the auxiliary determinant $\DetAbb(x,y,\y;n)$, which
arises from $\DetAb(x,y;n)$ (the determinant in (3.8))
by replacing each occurence of $y$ by
$\y$, except for the entries in the $0$-th row, where we only
partially replace $y$ by $\y$,
$$\multline
\DetAbb(x,y,\y;n)\\
{}:=\det_{0\le i,j\le n-1}\!\!\(
\hbox{\hskip0pt}
\raise4pt\hbox{
\SPfad(0,0),11111111111111\endSPfad
\SPfad(4,-2),2222\endSPfad
\Label\r{\hskip5pt\SB(x,y)}(1,1)
\Label\r{\ssize (x+\y+1)_j\,(x-j+2)_{\fl{(y-x)/2}+j-1}}(9,1)
\Label\r{\hskip5pt\ssize -(\y-i+2)_{i-1}}(1,-1)
\Label\ro{\raise-10pt\hbox{$\ssize (x+\y+i+1)_{j-1}\,(x+2i-j+2)_{j-1}$}}(8,-1)
\Label\ro{\hskip0pt\ssize \times(\y+2j-i+2)_{i-1}\,(\y-x+3j-3i)}(10,-2)
\Label\r{\hskip-4pt i\!=\!0}(17,1)
\Label\r{\hskip-4pt i\!\ge\!1}(17,-1)
}
\hskip7.2cm
\),
\Label\r{\hskip-12pt j\!=\!0}(-15,-3)
\Label\r{j\!\ge\!1}(-7,-3)
\Label\r{\hskip7pt(3.12)}(1,-3)
\hskip1cm
\endmultline
$$
\vskip3pt
\noindent
with $\SB(x,y)$ given by (3.9). Clearly, $\DetAbb(x,y,\y;n)$ is a
polynomial in $x$ and $\y$ (recall that $y=x+m$) which agrees with
$\DetAb(x,y;n)$ when $\y=y$. We are going to prove that
$$\DetAbb(x,y,\y;n)=\prod _{i=1}
^{n-1}\big((2x+\y+2i+2)_{i-1}\,(x+2\y+2i+2)_{i-1}\big)\cdot
P_2(x,y,\y;n),\tag3.13$$
where $P_2(x,y,\y;n)$ is a polynomial in $x$ and $\y$. Obviously,
when we set $\y=y$,
this implies that $\prod _{i=1}
^{n-1}\big((2x+y+2i+2)_{i-1}\,(x+2y+2i+2)_{i-1}\big)$ is a factor of
$\DetAb(x,y;n)$, as desired.

To prove (3.13), we first consider just one half of this product, 
$\prod _{i=1}
^{n-1}(2x+\y+2i+2)_{i-1}$. Let us concentrate on a typical factor
$(2x+\y+2i+l+1)$, $1\le i\le n-1$, $1\le l<i$. We claim that for each
such factor there is a linear combination of the rows that vanishes
if the factor vanishes. More precisely, we claim that for any $i,l$
with $1\le i\le n-1$, $1\le l<i$ there holds
$$\multline \sum _{s=l} ^{\fl{(i+l)/2}}\frac
{(2i-3s+l)} {(i-s)}\frac {(i-2s+l+1)_{s-l}} {(s-l)!}\frac
{(x+2s+1)_{2i-2s}} {(-x-2i-l+s)_{i-s}}\\
\cdot(\text {row $s$ of
$\DetAbb(x,y,-2x-2i-l-1;n)$})
=(\text {row $i$ of $\DetAbb(x,y,-2x-2i-l-1;n)$}).
\endmultline\tag3.14$$
To see this, we have to check
$$\multline \sum _{s=l} ^{\fl{(i+l)/2}}\frac
{(2i-3s+l)} {(i-s)}\frac {(i-2s+l+1)_{s-l}} {(s-l)!}\frac
{(x+2s+1)_{2i-2s}} {(-x-2i-l+s)_{i-s}}\\
\cdot(-2x-2i-l-s+1)_{s-1}=
(-2x-3i-l+1)_{i-1},
\endmultline\tag3.15$$
which is (3.14) restricted to the $0$-th column, and
$$\multline \sum _{s=l} ^{\fl{(i+l)/2}}\frac
{(2i-3s+l)} {(i-s)}\frac {(i-2s+l+1)_{s-l}} {(s-l)!}\frac
{(x+2s+1)_{2i-2s}} {(-x-2i-l+s)_{i-s}}\\
\times(-x-2i-l+s)_{j-1}\,(x+2s-j+2)_{j-1}\\
\times(-2x-2i-l+2j-s+1)_{s-1}\,(-3x-2i-l+3j-3s-1)\\
=(-x-i-l)_{j-1}\,(x+2i-j+2)_{j-1}\,(-2x-3i-l+2j+1)_{i-1}\,(-3x-l+3j-5i-1),
\endmultline\tag3.16$$
which is (3.14) restricted to the $j$-th column, $1\le j\le n-1$.
Equivalently, in terms of hypergeometric series (cf\. the Appendix
for the definition of the $F$-notation), this means to check
$$\multline 
2\,{{  ({ \textstyle 1 - 2 i - 2 l - 2 x}) _{l-1}  
      ({ \textstyle 1 + 2 l + x}) _{2 i - 2 l} }\over 
    {({ \textstyle -2 i - x}) _{i - l} }} \\
\times
{} _{5} F _{4} \!\left [ \matrix { 1 - {{2 i}\over 3} + {{2 l}\over 3},
       -{{i}\over 2} + {l\over 2}, {1\over 2} - {i\over 2} + {l\over 2}, -2 i
       - x, 2 i + 2 l + 2 x}\\ { -{{2 i}\over 3} + {{2 l}\over 3}, 1 - i + l,
       {1\over 2} + l + {x\over 2}, 1 + l + {x\over 2}}\endmatrix ;
       {\displaystyle 1}\right ]  
\\= 
  ({ \textstyle -2x - 3 i - l+1 }) _{i-1} 
\endmultline\tag3.17$$
and
$$\multline 
2\,\frac {({ \textstyle -x-2i}) _{j-1}  } 
{({ \textstyle -x-2i}) _{i - l} }
{{\left( - 3 x - 4 l + 3 j - 2 i -1  \right) 
({ \textstyle - 2 x - 2 l + 2 j - 2 i +1}) _{l-1}  
             }
    }   \\
\times ({ \textstyle x - j + 2 l + 2}) _{ 2 i + j - 2 l-1}
\cdot      {} _{6} F _{5} \!\left [ \matrix  {4\over 3} + {{2 i}\over 3} - j + {{4
       l}\over 3} + x, 1 - {{2 i}\over 3} + {{2 l}\over 3}, \\  {1\over 3} + {{2 i}\over 3} - j + {{4 l}\over 3}
       + x, -{{2 i}\over 3} + {{2 l}\over 3}, \endmatrix\right.
\hskip1.5cm\\
\hskip3cm
\left.\matrix -{{i}\over 2} +
       {l\over 2}, {1\over 2} - {i\over 2} + {l\over 2}, -1 - 2 i + j - x, 2 i
       - 2 j + 2 l + 2 x\\1 - i + l, 1 - {j\over 2} + l +
       {x\over 2}, {3\over 2} - {j\over 2} + l + {x\over 2}\endmatrix;
       {\displaystyle 1}\right ]  
\\ = 
   \left( -3x-l+ 3 j -5i-1 \right)  
     ({ \textstyle 1 - 3 i + 2 j - l - 2 x}) _{i-1} \hskip3cm \\
\times     ({ \textstyle -i - l - x}) _{j-1}  
     ({ \textstyle 2 + 2 i - j + x}) _{j-1}   
\endmultline\tag3.18$$
Now, the identity (3.17) holds since the $_5F_4$-series in (3.17) can
be summed by Corollary~A5, and the identity (3.18) holds since the
$_6F_5$-series in
(3.18) can be summed by Lemma~A6.

The product $\prod _{i=1}
^{n-1}(2x+\y+2i+2)_{i-1}$ consists of factors of the form $(2x+\y+a)$,
$4\le a\le 3n-3$. Let $a$ be
fixed. Then the factor $(2x+\y+a)$ occurs in the product $\prod _{i=1}
^{n-1}(2x+\y+2i+2)_{i-1}$ as many times as there are solutions to the
equation
$$a=2i+l+1,\quad \text {with }1\le i\le n-1,\ 1\le l<i.\tag3.19$$
For each solution $(i,l)$, we subtract the linear combination
$$\multline \sum _{s=l} ^{\fl{(i+l)/2}}\frac
{(2i-3s+l)} {(i-s)}\frac {(i-2s+l+1)_{s-l}} {(s-l)!}\frac
{(x+2s+1)_{2i-2s}} {(-x-2i-l+s)_{i-s}}\\
\cdot(\text {row $s$ of
$\DetAbb(x,y,\y;n)$})
\endmultline\tag3.20$$
of rows of $\DetAbb(x,y,\y;n)$
from row $i$ of $\DetAbb(x,y,\y;n)$. Then, by (3.14), all the entries
in row $i$ of the resulting determinant vanish for $\y=-2x-2i-l-1$.
Hence, $(2x+\y+2i+l+1)=(2x+\y+a)$ is a factor of all the entries in row
$i$, for each solution $(i,l)$ of (3.19). By taking these factors out
of the determinant we obtain
$$\DetAbb(x,y,\y;n)=(2x+\y+a)^{\#(\text {solutions $(i,l)$ of
(3.19)})}\cdot \DetAb^{(a)}(x,y,\y;n),\tag3.21$$
where $\DetAb^{(a)}(x,y,\y;n)$ is a determinant whose entries are
rational functions in $x$ and $\y$, the denominators containing
factors of the form $(x+c)$ (which come from the coefficients in the
linear combination (3.20)). Taking the limit $x\to -c$ in (3.21) then
reveals that these denominators cancel in the determinant, so that
$\DetAb^{(a)}(x,y,\y;n)$ is actually a polynomial in $x$ and $\y$.
Thus we have shown that each factor of $\prod _{i=1}
^{n-1}(2x+\y+2i+2)_{i-1}$ divides $\DetAbb(x,y,\y;n)$ with the right
multiplicity, hence the complete product divides $\DetAbb(x,y,\y;n)$.

The reasoning that $\prod _{i=1}
^{n-1}(x+2\y+2i+2)_{i-1}$ is a factor of $\DetAbb(x,y,\y;n)$ is
similar. Also here, let us concentrate on a typical factor
$(x+2\y+2j+l+1)$, $1\le j\le n-1$, $1\le l<j$. This time
we claim that for each
such factor there is a linear combination of the columns that vanishes
if the factor vanishes. More precisely, we claim that for any $j,l$
with $1\le j\le n-1$, $1\le l<j$ there holds
$$\multline \sum _{s=l} ^{\fl{(j+l)/2}}\frac
{(2j-3s+l)} {(j-s)}\frac {(j-2s+l+1)_{s-l}} {(s-l)!}
(\y+2s+1)_{2j-2s}\\
\cdot(\text {column $s$ of
$\DetAbb(-2\y-2j-l-1,y,\y;n)$})\\
=(\text {column $j$ of $\DetAbb(-2\y-2j-l-1,y,\y;n)$}).
\endmultline\tag3.22$$
This means to check
$$\multline \sum _{s=l} ^{\fl{(j+l)/2}}\frac
{(2j-3s+l)} {(j-s)}\frac {(j-2s+l+1)_{s-l}} {(s-l)!}
(\y+2s+1)_{2j-2s}\\
\times(-\y-2j-l)_s\,(-2\y-2j-l-s+1)_{\fl{(y-x)/2}+s-1}\\
=
(-\y-2j-l)_j\,(-2\y-3j-l+1)_{\fl{(y-x)/2}+j-1},
\endmultline\tag3.23$$
which is (3.22) restricted to the $0$-th row, and
$$\multline \sum _{s=l} ^{\fl{(j+l)/2}}\frac
{(2j-3s+l)} {(j-s)}\frac {(j-2s+l+1)_{s-l}} {(s-l)!}
(\y+2s+1)_{2j-2s}\\
\times(-\y-2j-l+i)_{s-1}\,(-2\y-2j-l+2i-s+1)_{s-1}\\
\times(\y+2s-i+2)_{i-1}\,(3\y+2j+l-3i+3s+1)\\
=(-\y-2j-l+i)_{j-1}\,(-2\y-3j-l+2i+1)_{j-1}\,
(\y+2j-i+2)_{i-1}\,(3\y+5j+l-3i+1),
\endmultline\tag3.24$$
which is (3.22) restricted to the $i$-th row, $1\le i\le n-1$.
If we plug
$$(-\y-2j-l)_s=\frac {(-\y-2j-l)_j} {(-\y-2j-l+s)_{j-s}}$$
into (3.23), we see that (3.23) is equivalent to (3.15) (replace $x$
by $\y$ and $i$ by $j$). Likewise, by plugging
$$(-\y-2j-l+i)_{s-1}=\frac {(-\y-2j-l+s)_{i-1}} {(-\y-2j-l+s)_{j-s}}
\frac {(-\y-2j-l+i)_{j-1}} {(-\y-j-l)_{i-1}}
$$
into (3.24), we see that (3.24) is equivalent to (3.16) (replace $x$
by $\y$ and interchange $i$ and $j$). By arguments that are similar
to the ones above, it follows that $\prod _{i=1}
^{n-1}(x+2\y+2i+2)_{i-1}$ divides $\DetAbb(x,y,\y;n)$. 

Altogether, this implies that $\prod _{i=1}
^{n-1}\big((2x+\y+2i+2)_{i-1}
(x+2\y+2i+2)_{i-1}\big)$ divides $\DetAbb(x,y,\y;n)$,
and so, as we already noted after (3.13), 
the product $\prod _{i=1}
^{n-1}\big((2x+y+2i+2)_{i-1}
(x+2y+2i+2)_{i-1}\big)$ divides $\DetAb(x,y;n)$, as desired.

\smallskip
{\it Step 4. $\prod _{i=0}
^{\fl{n/2}-1}(2x+2\cl{(y-x)/2}+2i+1)$ is a factor of $\DetAb(x,y;n)$.} 
We consider (3.5). In the determinant in (3.5) we take 
$$\frac {1} {2}\frac {(x+y+i+1)!\,(2x+y+3i+4)} {(x+2i+2)!\,(y+2n-2)!}$$
out of the $i$-th row, $i=0,1,\dots,n-2$, we take
$$\frac {(x+y+n)!} {(x+2n-2)!\,(y+2n-2)!}$$
out of the $(n-1)$-st row, and we take
$$\frac {1} {2}(y+2j+3)_{2n-2j-4}\,(x+2y+3j+4)$$
out of the $j$-th column, $j=0,1,\dots,n-2$. Then we combine with
(3.6) and (3.10) (recall that $\DetAb(x,y;n)$ is the determinant in
(3.6)), and after cancellation we obtain
$$\multline
\DetAb(x,y;n)
=\(\frac {1} {2}\)^{2n-2}\frac {(x+y+1)_n}
{(\fl{(x+y)/2}+1)_{2n-2-\fl{(y-x)/2}}\,(y+1)_{2n-2}}
\\
\times
\det_{0\le i,j\le n-1}\!\!\(
\hbox{\hskip.5cm}
\raise0pt\hbox{
\SPfad(-1,0),1111111111111111\endSPfad
\SPfad(7,-2),2222\endSPfad
\Label\ru{\ssize (x+y+i+2)_j\,(x+2i-j+3)_{j}}(2,2)
\Label\ru{\raise6pt\hbox{$\ssize
\times(y+2j-i+3)_{i}\,(y-x+3j-3i)$\hskip8pt}}(3,1)
\Label\ru{\hskip-6pt\ssize (x+y+i+2)_{n-1}\,(x+2i-n+4)_{n-1}}(11,2)
\Label\ru{\raise6pt\hbox{$\ssize \times(y+2n-i-1)_{i}$}}(10,1)
\Label\ro{\raise-10pt\hbox{$\ssize
-(x+y+n+1)_{j}\,(x+2n-j-1)_{j}$\hskip15pt}}(3,-1)
\Label\ro{\hskip0pt\ssize \times(y+2j-n+4)_{n-1}\hskip15pt}(4,-2)
\Label\r{\hskip5pt\SC(x,y;n)}(10,-1)
\Label\r{\hskip-4pt i\!\le\!n\!-\!2}(17,1)
\Label\r{\hskip-4pt i\!=\!n\!-\!1}(17,-1)
}
\hskip7.5cm
\),
\Label\r{\hskip-12pt j\!\le\!n\!-\!2}(-14,-3)
\Label\r{j\!=\!n\!-\!1}(-6,-3)
\Label\r{\hskip8pt(3.25)}(1,-3)
\hskip1cm
\endmultline
$$
\vskip3pt
\noindent
where
$$\SC(x,y;n)=\sum _{r=0} ^{y-x-1}(x+y+n+1)_{n-2}\,(x+n+r)_{n-r-1}\,
(y+n-r)_{r+n-1}.$$
The determinant on the right-hand side of (3.25) has polynomial
entries. Note that in case of the $(n-1,n-1)$-entry this is due to
$n-r-1\ge n-(y-x-1)-1=n-m\ge0$ (recall that $y=x+m$), the last
inequality being an assumption in the statement of the Theorem. The
product in the numerator of the right-hand side of (3.25) consists of
factors of the form $(x+y+a)=(2x+m+a)$ with integral $a$. Some of these factors
cancel with the denominator, but all factors of the form $(2x+2b+1)$,
with integral $b$, do not cancel, and so because of (3.25) divide
$\DetAb(x,y;n)$. These factors are
$$\prod _{i=0}
^{\cl{(m+n)/2}-\cl{m/2}-1}(2x+2\cl{m/2}+2i+1)$$
(with $m=y-x$, of course).
Since 
$$\cl{(m+n)/2}-\cl{m/2}-1=\cases \fl{n/2}&\text {$n$ odd, $m$ even}\\
\fl{n/2}-1&\text {otherwise},\endcases$$
it follows that $\prod _{i=0}
^{\fl{n/2}-1}(2x+2\cl{m/2}+2i+1)$ is a factor of $\DetAb(x,y;n)$, and
if $n$ is odd and $m$ is even $(2x+m+n)$ is an additional factor of
$\DetAb(x,y;n)$.

\smallskip
Summarizing, so far we have shown that the equation (3.11) holds,
where $P_1(x;y-x,n)=P_1(x;m,n)$ 
is {\it some\/} polynomial in $x$, that has
$(2x+m+n)$ as a factor in case that $n$ is odd and $m$ is even. It
remains to show that $P_1(x;m,n)$ is a polynomial in $x$ of degree
$\le \fl{m/2}$, and to describe how $P_1(x;m,n)$ can be
computed explicitly.

\smallskip
{\it Step 5. $P_1(x;m,n)$ is a polynomial in $x$ of degree
$\le\fl{m/2}$.}
Here we write $x+m$ for $y$ everywhere.
We shall prove that $\DetAa(x,x+m;n)$ (which is defined to be the
determinant in (3.6)) is a polynomial in $x$ of degree at most
$2\binom n2+\binom{n-1}2 +\fl{n/2}+\fl{m/2}$.
By (3.10) this would
imply that $\DetAb(x,x+m;n)$ is a polynomial in $x$ of degree at most
$2\binom {n-1}2+\fl{n/2}+\fl{m/2}$, and so, by (3.11), that
$P_1(x;m,n)$ is a polynomial in $x$ of degree at most $\fl{m/2}$, as
desired.

Establishing the claimed degree bound for $\DetAa(x,x+m;n)$ is the
most delicate part of the proof of the Theorem. We need to consider the
generalized determinant 
$$\DetAaa(x,z(1),z(2),\dots,z(n-1);n)=\DetAaa(n)$$
which arises from $\DetAa(x,x+m;n)$ by replacing each occurence of $i$
in row $i$ by an indeterminate, $z(i)$ say, $i=1,2,\dots,n-1$,
$$\multline
\DetAaa(x,z(1),z(2),\dots,z(n-1);n)\\
:=\det_{0\le i,j\le n-1}\!\!\(
\hbox{\hskip.7cm}
\raise0pt\hbox{
\SPfad(-1,0),11111111111111111\endSPfad
\SPfad(4,-2),2222\endSPfad
\Label\r{\hskip-10pt\SA(x,x+m;n)}(1,1)
\Label\ru{\ssize (2x+m+1)_j\,(x-j+2)_{\fl{m/2}+j-1}}(10,2)
\Label\ru{\raise6pt\hbox{\hskip-8pt$\ssize \times(x+m+2j+1)_{2n-2j-2}$}}(10,1)
\Label\r{\hskip-13pt\ssize -(x+m-z(i)+2)_{2n-2}}(1,-1)
\Label\ro{\raise-10pt\hbox{\hskip-10pt$\ssize (2x+m+z(i)+1)_{j-1}\,(x+2z(i)-j+2)_{j-1}$}}(9,-1)
\Label\ro{\hskip-20pt\ssize \times(x+m+2j-z(i)+2)_{2n-2j-2}\,(m+3j-3z(i))}(11,-2)
\Label\r{\hskip-8pt i\!=\!0}(18,1)
\Label\r{\hskip-8pt i\!\ge\!1}(18,-1)
}
\hskip8cm
\),
\Label\r{\hskip-12pt j\!=\!0}(-16,-3)
\Label\r{j\!\ge\!1}(-7,-3)
\Label\r{\hskip7pt(3.26)}(0,-3)
\hskip.5cm
\endmultline
$$
\vskip3pt
\noindent
where $\SA(x,x+m;n)$ is given by (3.7).

This determinant is a polynomial in $x, z(1),z(2),\dots,z(n-1)$. We
shall prove that the degree in $x$ of this determinant is at most
$2\binom n2+\binom{n-1}2 +\fl{n/2}+\fl{m/2}$, 
which clearly implies our claim upon setting $z(i)=i$,
$i=1,2,\dots,n-1$.

Let us denote the $(i,j)$-entry of $\DetAaa(n)$ by $A(i,j)$. In the
following computation we write $S_n$ for the group of all
permutations of $\{0,1,\dots,n-1\}$. By definition of the determinant
we have
$$\DetAaa(n)=\sum _{\si\in S_n} ^{}\sgn \si\prod _{j=0}
^{n-1}A(\si(j),j),$$
and after expanding the determinant along the $0$-th row,
$$\multline \DetAaa(n)=A(0,0)\sum _{\si\in S_{n-1}} ^{}\prod _{j=0}
^{n-2}A(\si(j)+1,j+1)\\
+\sum _{\ell=1} ^{n-1}(-1)^\ell A(0,\ell)\sum _{\si\in S_{n-1}}
^{}\sgn\si\prod _{j=0} ^{n-2}A(\si(j)+1,j+\chi(j\ge\ell)),
\endmultline\tag3.27$$
where $\chi(\Cal A)$=1 if $\Cal A$ is
true and $\chi(\Cal A)$=0 otherwise. Now, by Lemma~A10 we know that
for $i,j\ge1$ we have
$$A(i,j)=\sum _{p,q\ge0}
^{}2^j\al_{p,q}(j)\,x^{p}\,z(i)^{q},\tag3.28
$$
where $\al_{p,q}(j)$ is a polynomial in $j$ of degree $\le
2(2n-3-p-q)+q-1$. It should be noted that the range of the sum in
(3.28) is actually 
$$0\le p\le 2n-4,\ 0\le q\le 2n-3,\ p+q\le 2n-3.\tag3.29$$
Furthermore, by Lemma~A9 we know that for $j\ge1$ we have
$$A(0,j)=\sum _{p\ge0}
^{}2^j\be_{p}(j)\,x^{p},\tag3.30
$$
where $\be_{p}(j)$ is a polynomial in $j$ of degree $\le
2(2n+\fl{m/2}-3-p)$. Also, for $i\ge1$ let
$$A(i,0)=-(x+m-z(i)+2)_{2n-2}=\sum _{p,q\ge0}
^{}\ga_{p,q}\,x^p\,z(i)^q.\tag3.31$$
Plugging (3.28) and (3.31) into (3.27), and writing $\z(i)$ instead
of $z(i+1)$ for notational convenience, we get
$$\multline \hskip-18pt\DetAaa(n)=A(0,0)\underset q_0,\dots,q_{n-2}\ge0\to {\sum
_{p_0,\dots,p_{n-2}\ge0} ^{}}2^{\binom n2}x^{p_0+\dots+p_{n-2}}\prod _{j=0}
^{n-2}\al_{p_j,q_j}(j+1)
\sum _{\si\in S_{n-1}} ^{}\sgn\si \prod _{j=0}
^{n-2}\z(\si(j))^{q_j}\\
+\sum _{\ell=1} ^{n-1}(-1)^\ell A(0,\ell)
\underset q_0,\dots,q_{n-2}\ge0\to {\sum
_{p_0,\dots,p_{n-2}\ge0} ^{}}2^{\binom n2-\ell}
x^{p_0+\dots+p_{n-2}}\ga_{p_0,q_0}\prod
_{j=1} ^{n-2}\al_{p_j,q_j}(j+\chi(j\ge \ell))\\
\hskip2cm \times\sum _{\si\in S_{n-1}} ^{}\sgn\si \prod _{j=0}
^{n-2}\z(\si(j))^{q_j}\\
=A(0,0)\underset q_0,\dots,q_{n-2}\ge0\to {\sum
_{p_0,\dots,p_{n-2}\ge0} ^{}}2^{\binom n2}x^{p_0+\dots+p_{n-2}}\prod _{j=0}
^{n-2}\al_{p_j,q_j}(j+1)\det_{0\le i,j,\le n-2}(\z(i)^{q_j})
\hskip2cm\\
+\sum _{\ell=1} ^{n-1}(-1)^\ell A(0,\ell)
\underset q_0,\dots,q_{n-2}\ge0\to {\sum
_{p_0,\dots,p_{n-2}\ge0} ^{}}2^{\binom n2-\ell}
x^{p_0+\dots+p_{n-2}}\ga_{p_0,q_0}\prod
_{j=1} ^{n-2}\al_{p_j,q_j}(j+\chi(j\ge \ell))\\
\times\det_{0\le i,j,\le n-2}(\z(i)^{q_j}).\hskip3cm
\endmultline\tag3.32$$
The determinants in (3.32) vanish whenever $q_{j_1}=q_{j_2}$ for some
$j_1\ne j_2$. Hence, in the sequel we may assume that the summation
indices $q_0,q_1,\dots,q_{n-2}$ are pairwise distinct, in both terms
on the right-hand side of (3.32). In particular, we may assume that
in the first term
the pairs $(p_0,q_0),(p_1,q_1),\dots,(p_{n-2},q_{n-2})$ are
pairwise distinct, and that in the second term 
the pairs $(p_1,q_1),(p_2,q_2),\dots,(p_{n-2},q_{n-2})$ are
pairwise distinct. What we do next is to collect the summands in the
inner sums that are indexed by the same {\it set\/} of pairs. So,
if in addition we plug (3.30) into (3.32), we obtain
$$\multline \DetAaa(n)=A(0,0) {\sum
_{\{(p_0,q_0),\dots,(p_{n-2},q_{n-2})\}} ^{}}2^{\binom n2}
x^{p_0+\dots+p_{n-2}}\\
\times
\sum _{\ta\in S_{n-1}} ^{}\det_{0\le i,j\le n-2}(\z(i)^{q_{\ta(j)}})
\prod _{j=0} ^{n-2}\al_{p_{\ta(j)},q_{\ta(j)}}(j+1)\\
+\sum _{p,p_0,q_0\ge0} ^{}{\sum
_{\{(p_1,q_1),\dots,(p_{n-2},q_{n-2})\}} ^{}}
2^{\binom n2}x^{p+p_0+\dots+p_{n-2}}\ga_{p_0,q_0}
\sum _{\ell=1} ^{n-1}(-1)^\ell \be_p(\ell)\\
\times\sum _{\ta\in \tilde S_{n-2}} ^{}\det_{0\le i,j\le n-2}(\z(i)^{q_{\ta(j)}})
\prod _{j=1} ^{n-2}\al_{p_{\ta(j)},q_{\ta(j)}}(j+\chi(j\ge \ell))
\endmultline\tag3.33$$
where $\tilde S_{n-2}$ denotes the group of all permutations of
$\{0,1,\dots,n-1\}$ that fix $0$. Clearly, we have
$$\det_{0\le i,j\le n-2}(\z(i)^{q_{\ta(j)}})=\sgn\ta\,
\det_{0\le i,j\le n-2}(\z(i)^{q_{j}}).\tag3.34$$
Moreover, there holds
$$\align \sum _{\ell=1} ^{n-1}(-1)^\ell \be_p(\ell)&
\sum _{\ta\in \tilde S_{n-2}} ^{}\sgn\ta
\prod _{j=1} ^{n-2}\al_{p_{\ta(j)},q_{\ta(j)}}(j+\chi(j\ge \ell))\\
&=\sum _{\ell=1} ^{n-1}(-1)^\ell \be_p(\ell)
\det_{1\le i,j\le n-2}\big(\al_{p_i,q_i}(j+\chi(j\ge\ell))\big)\\
&=(-1)^{n-1}\det_{1\le i,j\le n-1}\pmatrix \al_{p_i,q_i}(j)&
\hbox{\eightpoint$i\!\le\!n\!-\!2$}\\
\be_p(j)& \hbox{\eightpoint$i\!=\!n\!-\!1$}\endpmatrix, 
\tag3.35\endalign$$
the step from the last line to the next-to-last line being just
expansion of the determinant along the bottom row. Using (3.34) and
(3.35) in (3.33) then yields
$$\multline \DetAaa(n)=A(0,0) {\sum
_{\{(p_0,q_0),\dots,(p_{n-2},q_{n-2})\}} ^{}}2^{\binom n2}
x^{p_0+\dots+p_{n-2}}\\
\hskip2cm\times
\det_{0\le i,j\le n-2}(\z(i)^{q_{j}})
\det_{0\le i,j\le n-2}\big(\al_{p_{i},q_{i}}(j+1)\big)\\
+(-1)^{n-1}\sum _{p,p_0,q_0\ge0} ^{}{\sum
_{\{(p_1,q_1),\dots,(p_{n-2},q_{n-2})\}} ^{}}
2^{\binom n2}x^{p+p_0+\dots+p_{n-2}}\ga_{p_0,q_0}\\
\times
\det_{0\le i,j\le n-2}(\z(i)^{q_{j}})
\det_{1\le i,j\le n-1}\pmatrix \al_{p_i,q_i}(j)&
\hbox{\eightpoint$i\!\le\!n\!-\!2$}\\
\be_p(j)& \hbox{\eightpoint$i\!=\!n\!-\!1$}\endpmatrix.
\endmultline\tag3.36$$
We treat the two terms on the right-hand side of (3.36) separately.
Recall that we want to prove that the degree in $x$ of $\DetAaa(n)$
is at most 
$2\binom n2+\binom{n-1}2 +\fl{n/2}+\fl{m/2}$.

What regards the first term,
$$\multline A(0,0) {\sum
_{\{(p_0,q_0),\dots,(p_{n-2},q_{n-2})\}} ^{}}2^{\binom n2}x^{p_0+\dots+p_{n-2}}\\
\hskip2cm\times
\det_{0\le i,j\le n-2}(\z(i)^{q_{j}})
\det_{0\le i,j\le n-2}\big(\al_{p_{i},q_{i}}(j+1)\big),
\endmultline\tag3.37$$
we shall prove that the degree in $x$ is actually at most 
$2\binom n2+\binom{n-1}2 +\fl{(n-1)/2}+\fl{m/2}$. Equivalently, when
disregarding $A(0,0)=\SA(x,x+m;n)$, 
whose degree in $x$ is $2n-2+\fl{m/2}$ (see (3.7)), this means to
prove that the degree in $x$ of the sum in (3.37) is at most 
$3\binom{n-1}2 +\fl{(n-1)/2}$.

So we have to examine for which indices
$p_0,\dots,p_{n-2},q_0,\dots,q_{n-2}$ the determinants in (3.37) do
not vanish. As we already noted, the first determinant does not
vanish only if the indices $q_0,q_1,\dots,q_{n-2}$ are pairwise
distinct. So, without loss of generality we may assume 
$$0\le q_0<q_1<\dots<q_{n-2}.\tag3.38$$

Turning to the second determinant in (3.37), we observe that
because of what we know about $\al_{p_i,q_i}(j+1)$ (cf\. the sentence
containing (3.28))
each row of this determinant is filled with a {\it single\/}
polynomial evaluated at $1,2,\dots,n-1$. Let $M$ be some nonnegative
integer. If we assume that among the polynomials
$\al_{p_i,q_i}(j+1)$, $i=0,1,\dots,n-2$, there are $M+1$ polynomials
of degree less or equal $M-1$, then the determinant will vanish. For,
a set of $M+1$ polynomials of maximum  degree $M-1$ is linearly
dependent. Hence, the rows in the second determinant in (3.37) will
be linearly dependent, and so the determinant will vanish. Since the
degree of $\al_{p_i,q_i}(j+1)$ as a polynomial in $j$ is at most
$2(2n-3-p_i-q_i)+q_i-1$ (again, cf\. the sentence
containing (3.28)), we have that
$$\matrix \format\l\\
\text {the number of integers $2(2n-3-p_i-q_i)+q_i-1$,
$i=0,1,\dots,n-2$,}\\
\text {that are less or equal $M-1$ is at most
$M$}.\endmatrix\tag3.39$$

Now the task is to determine the maximal value of
$p_0+p_1+\dots+p_{n-2}$ (which is the degree in $x$ of the sum in (3.37)
that we are interested in), under the conditions (3.38) and (3.39),
and the additional condition 
$$0\le p_i\le 2n-4,\ 0\le q_i\le 2n-3,\ p_i+q_i\le 2n-3,\tag3.40$$
which comes from (3.29). We want to prove that this maximal value is
$3\binom{n-1}2 +\fl{(n-1)/2}$. To simplify notation we write
$$\ep_i=2n-3-p_i-q_i.\tag3.41$$ 
Thus, since
$$\sum _{i=0} ^{n-2}p_i=(n-1)(2n-3)-\sum _{i=0} ^{n-2}(q_i+\ep_i),$$
we have to prove that the {\it minimal\/} value of
$$q_0+q_1+\dots+q_{n-2}+\ep_0+\ep_1+\dots+\ep_{n-2},\tag3.42$$
under the condition (3.38), the condition that 
$$\ep_i\ge0,\quad  i=0,1,\dots,n-2,\tag3.43$$
(which comes from the right-most inequality in (3.40) under the
substitution (3.41)), the condition
that
$$\matrix \format\l\\
\text {the number of integers $2\ep_i+q_i-1$,
$i=0,1,\dots,n-2$,}\\
\text {that are less or equal $M-1$ is at most
$M$},\endmatrix\tag3.44$$
(which is (3.39) under the substitution (3.41)), and the condition
$$\text {if $q_0=0$, then $\ep_0\ge1$},\tag3.45$$
(which comes from (3.40) and (3.41)), is $\binom
{n-1}2+\cl{(n-1)/2}$. 

As a first, simple case, we consider $q_0\ge1$. Then, from (3.38) it
follows that the sum $\sum _{i=0} ^{n-2}q_i$ alone is at least
$\binom n2=\binom{n-1}2+(n-1)$, which trivially implies our claim.
Therefore, from now on we assume that $q_0=0$. Note that this in
particular implies $\ep_0\ge1$, because of (3.45).

Next, we apply (3.44) with $M=2$. In particular, since among the first
three integers $2\ep_i+q_i-1$, $i=0,1,2$, only two can be less or
equal 1, there must be an $i_1\le 2$ with $2\ep_{i_1}+q_{i_1}-1\ge
2$. Without loss of generality we choose $i_1$ to be minimal with
this property. Now we apply (3.44) with $M=2\ep_{i_1}+q_{i_1}$.
Arguing similarly, we see that there must be an $i_2\le 2\ep_{i_1}+q_{i_1}$
with $2\ep_{i_2}+q_{i_2}-1\ge 2\ep_{i_1}+q_{i_1}$. Again, we choose
$i_2$ to be minimal with this property. This continues, until we meet an
$i_k\le 2\ep_{i_{k-1}}+q_{i_{k-1}}$ with
$2\ep_{i_{k}}+q_{i_{k}}-1\ge n-2$. That such an $i_k$ must
be found eventually is seen by applying (3.44) with $n-2$.

Let us collect the facts that we
have found so far: There exists a sequence $i_1,i_2,\dots,\mathbreak i_k$ of
integers satisfying
$$0\le i_1<i_2<\dots<i_k\le n-2\tag3.46$$
(this is because of the minimal choice for each of the $i_j$'s),
$$i_1\le 2,\ i_2\le 2\ep_{i_1}+q_{i_1},\ \dots,\ i_k\le
2\ep_{i_{k-1}}+q_{i_{k-1}},\tag3.47$$
and
$$2\ep_{i_k}+q_{i_k}-1\ge n-2.\tag3.48$$
The other inequalities are not needed later. 

Now we turn to the quantity (3.42) that we want to bound from above.
We have
$$\sum _{i=0} ^{n-2}q_i+\sum _{i=0} ^{n-2}\ep_i=
\sum _{i=0} ^{n-2}(q_i-i)+\binom{n-1}2+\sum _{i=0}
^{n-2}\ep_i.\tag3.49$$
For convenience, we write $\q_i$ for $q_i-i$ in the sequel. Because of
(3.38) we have 
$$\q_i\ge0,\quad i=0,1,\dots,n-2,\tag3.50$$
and
$$\q_i\ge \q_j,\quad \text {for }i\ge j.\tag3.51$$
For a fixed $i$ let $s$ be maximal such that $i_s\le i$. Then, because
of (3.51), there holds
$$\align q_i-i&=\q_i=(\q_i-\q_{i_s})+(\q_{i_s}-\q_{i_{s-1}})
+\dots+(\q_{i_2}-\q_{i_1})+\q_{i_1}\\
&\ge(\q_{i_s}-\q_{i_{s-1}})+\dots+(\q_{i_2}-\q_{i_1})+\q_{i_1}.
\endalign$$
Using this, (3.50) and (3.48), in (3.49), we obtain
$$\multline
\sum _{i=0} ^{n-2}q_i+\sum _{i=0} ^{n-2}\ep_i\ge\binom
{n-1}2+(n-1-i_1)\,\q_{i_1}+\sum _{s=2}
^{k}(n-1-i_s)(\q_{i_s}-\q_{i_{s-1}})\\
+\sum _{i=0} ^{n-2}\ep_i-\ep_{i_k}+\frac {n-1-q_{i_k}} {2}.
\endmultline\tag3.52$$
Now, by (3.47) we have for $q_{i_k}$ that
$$\align
q_{i_k}&=\q_{i_k}+i_k=\q_{i_1}+\sum _{s=2}
^{k}(\q_{i_s}-\q_{i_{s-1}})+i_k\tag3.53\\
&\le \q_{i_1}+\sum _{s=2}
^{k}(\q_{i_s}-\q_{i_{s-1}})+2\ep_{i_{k-1}}+q_{i_{k-1}}.
\endalign$$
A similar estimation holds for $q_{i_{k-1}}$, etc. Thus, by iteration
we arrive at
$$
q_{i_k}
\le k\q_{i_1}+\sum _{s=2}
^{k}(k-s+1)(\q_{i_s}-\q_{i_{s-1}})+2(\ep_{i_1}+\dots+\ep_{i_{k-1}})
+i_1.
$$
Using this inequality in (3.52), we get
$$\align
\sum _{i=0} ^{n-2}q_i+\sum _{i=0} ^{n-2}\ep_i&\ge\binom
{n-1}2+\frac {n-1} {2}\tag3.54a\\
&\quad \quad +(n-1-i_1-\frac {k} {2})\,\q_{i_1}+\sum _{s=2}
^{k}(n-1-i_s-\frac {k-s+1} {2})\,(\q_{i_s}-\q_{i_{s-1}})\tag3.54b\\
&\quad \quad +\sum _{i=0} ^{n-2}\ep_i-\sum _{s=1} ^{k}\ep_{i_s}
-\frac {i_1} {2}.\tag3.54c
\endalign$$
By (3.50), (3.51), and since because of (3.46) we have
$$n-1-i_s-\frac {k-s+1} {2}\ge n-1-(n-2-k+s)-\frac {k-s+1} {2}=\frac
{k-s+1} {2}\ge0,\tag3.55$$
all terms in the line (3.54b) are nonnegative. If $i_1=0$, then by
(3.43) the line (3.54c) is nonnegative. If $1\le i_1\le 2$ ($i_1$
cannot be larger because of (3.47)), then $\ep_0$ occurs in the line
(3.54c). As we already noted, we have $\ep_0\ge 1$ since we are
assuming that $q_0=0$ in which case (3.45) applies. So,
$\ep_0-i_1/2\ge 0$, which in combination with (3.43) again implies
that the line (3.54c) is nonnegative.

Hence, we conclude
$$\sum _{i=0} ^{n-2}q_i+\sum _{i=0} ^{n-2}\ep_i\ge \binom{n-1}2+\frac
{n-1} {2},$$
which is what we wanted.

The reasoning for the second term om the right-hand side of (3.36),
$$\multline (-1)^{n-1}\sum _{p,p_0,q_0\ge0} ^{}{\sum
_{\{(p_1,q_1),\dots,(p_{n-2},q_{n-2})\}} ^{}}
2^{\binom n2}x^{p+p_0+\dots+p_{n-2}}\ga_{p_0,q_0}\\
\times
\det_{0\le i,j\le n-2}(\z(i)^{q_{j}})
\det_{1\le i,j\le n-1}\pmatrix \al_{p_i,q_i}(j)&
\hbox{\eightpoint$i\!\le\!n\!-\!2$}\\
\be_p(j)& \hbox{\eightpoint$i\!=\!n\!-\!1$}\endpmatrix,
\endmultline\tag3.56$$
is similar, only slightly more complicated. We shall prove that the
degree in $x$ in (3.56) is at most
$2\binom n2+\binom{n-1}2 +\fl{n/2}+\fl{m/2}$, which by the discussion
in the first paragraph of Step~5 is what we need.

So, we have to determine the maximal value of $p+p_0+\dots+p_{n-2}$
such that the determinants in (3.56) do
not vanish. Basically, we would now more or less run through 
the same arguments as
before. Differences arise mainly in the considerations concerning the
second determinant (which is slightly different from the second
determinant in (3.37)). What has to be used here is that $\be_p(j)$
is a polynomial in $j$ of degree $\le
2(2n+\fl{m/2}-3-p)$ (see the sentence containing (3.30)). 
If we make again the substitutions
$$\ep_i=2n-3-p_i-q_i,\quad i=1,2,\dots,n-2,\tag3.57$$ 
and in addition the substitutions
$$\ep_0=2n-2-p_0-q_0,\tag3.58$$ 
and
$$\ep=2n+\fl{m/2}-3-p,\tag3.59$$ 
we obtain eventually the following conditions that are necessary to make
these two determinants not vanish: 
There must hold
$$0\le q_1<q_2<\dots<q_{n-2},\quad \text {and $q_0$ is distinct from
the other $q_i$'s,}\tag3.60$$
(this is the substitute for (3.38)),
$$\ep_i\ge0,\quad  i=0,1,\dots,n-2,\quad \text {and}\quad
\ep\ge0,\tag3.61$$
(this is the substitute for (3.43)), and finally,
$$\matrix \format\l\\
\text {the number of integers in the set $\{2\ep_i+q_i-1:
i=1,2,\dots,n-2\}\cup\{2\ep\}$}\\
\text {that are less or equal $M-1$ is at most
$M$},\endmatrix\tag3.62$$
(this is the substitute for (3.44)). 
Since by the substitutions (3.57)--(3.59) we have
$$
p+\sum _{i=0} ^{n-2}p_i=
2\binom n2+2\binom{n-1}2+(n-1)+
\fl{\frac {m} {2}}-\ep-\sum _{i=0}
^{n-2}(q_i+\ep_i),
$$
the task is to prove that the {\it minimal\/} value of
$$q_0+q_1+\dots+q_{n-2}+\ep+\ep_0+\ep_1+\dots+\ep_{n-2},\tag3.63$$
equals $\binom{n-1}2+\cl{(n-2)/2}$.

Next in the arguments for the first term on the right-hand side of
(3.36) came
the sequence of applications of (3.44). Hence, now we apply (3.62)
repeatedly. Actually, there is only one slight change, with the start.
Namely, first we apply (3.62) with $M=2\ep+1$. Since then $2\ep$ is
already less or equal $M-1$, among the first
$2\ep+1$ integers $2\ep_i+q_i-1$, $i=1,2,\dots,2\ep+1$, only $2\ep$ 
can be less or
equal $2\ep$. Hence 
there must be an $i_1\le 2\ep+1$ with $2\ep_{i_1}+q_{i_1}-1\ge
2\ep+1$. Continuing in the same manner as before, we obtain a
sequence $i_1,i_2,\dots,i_k$ of
integers satisfying
$$1\le i_1<i_2<\dots<i_k\le n-2,\tag3.64$$
$$i_1\le 2\ep+1,\ i_2\le 2\ep_{i_1}+q_{i_1},\ \dots,\ i_k\le
2\ep_{i_{k-1}}+q_{i_{k-1}},\tag3.65$$
and
$$2\ep_{i_k}+q_{i_k}-1\ge n-2.\tag3.66$$

Now we turn to the quantity (3.63) that we want to bound from above.
We want to parallel the computation (3.49)--(3.54). However, since by
(3.60) the $q_i$'s are slightly unordered (in comparison with
(3.38)), we have to modify the definition of $\q_i$. Namely, let $t$
be the uniquely determined integer such that $q_{t}<q_0<q_{t+1}$, if
existent, or $t=0$ if $q_0<q_1$, or $t=n-2$ if $q_{n-2}<q_0$. Then we
define
$$\q_i:=\cases q_0-t&\text {if }i=0,\\
q_i-i+1&\text {if }1\le i\le t\\
q_i-i&\text {if }i>t.\endcases$$
If we modify (3.49) accordingly,
$$\sum _{i=0} ^{n-2}q_i+\sum _{i=0} ^{n-2}\ep_i=
(q_0-t)+
\sum _{i=1} ^{t}(q_i-i+1)+
\sum _{i=t+1} ^{n-2}(q_i-i)+\binom{n-1}2+\sum _{i=0}
^{n-2}\ep_i,$$
all subsequent steps that lead to (3.54) can be performed without
difficulties. (A little detail is that in (3.53) the equality
$q_{i_k}=\q_{i_k}+i_k$ has to be replaced by the inequality 
$q_{i_k}\le \q_{i_k}+i_k$.) 
Also, the estimation (3.55) still holds true because of
(3.64). Hence, when we use the first inequality in
(3.65), together with (3.50), (3.51), (3.55), (3.61), in (3.54), we obtain
$$\align
\sum _{i=0} ^{n-2}q_i+\ep+\sum _{i=0} ^{n-2}\ep_i&\ge\binom
{n-1}2+\frac {n-1} {2}
+\sum _{i=0} ^{n-2}\ep_i+\ep-\sum _{s=1} ^{k}\ep_{i_s}
-\frac {2\ep+1} {2}\\
&\ge \binom
{n-1}2+\frac {n-2} {2},
\endalign$$
which is what we wanted.

The proof that the degree of the polynomial $P_1(x;m,n)$ is at most
$\fl{m/2}$ is thus complete.

\smallskip
{\it Step 6. An algorithm for the explicit computation of
$P_1(x;m,n)$.} Also here, we write $x+m$ for $y$ everywhere. A
combination of (3.11) and (3.25) yields
$$\multline
P_1(x;m,n)=\(\frac {1} {2}\)^{2n-2}\frac {\prodl _{i=1}
^{\fl{(m+n)/2}-\fl{m/2}}(2x+2\fl{m/2}+2i)}
{(x+\fl{m/2}+1)_{2n-2-\fl{m/2}}\,(x+m+1)_{2n-2}}
\\
\times\prod _{i=1}
^{n-1}\frac {1} {(3x+m+2i+2)_{i-1}\,(3x+2m+2i+2)_{i-1}}\\
\times
\det_{0\le i,j\le n-1}\!\!\(
\hbox{\hskip.5cm}
\raise0pt\hbox{
\SPfad(-1,0),1111111111111111\endSPfad
\SPfad(7,-2),2222\endSPfad
\Label\ru{\ssize (2x+m+i+2)_j\,(x+2i-j+3)_{j}}(2,2)
\Label\ru{\raise6pt\hbox{$\ssize
\times(x+m+2j-i+3)_{i}\,(m+3j-3i)$\hskip12pt}}(3,1)
\Label\ru{\hskip-6pt\ssize (2x+m+i+2)_{n-1}\,(x+2i-n+4)_{n-1}}(11,2)
\Label\ru{\raise6pt\hbox{$\ssize \times(x+m+2n-i-1)_{i}$}}(10,1)
\Label\ro{\raise-10pt\hbox{$\ssize
-(2x+m+n+1)_{j}\,(x+2n-j-1)_{j}$\hskip15pt}}(3,-1)
\Label\ro{\hskip0pt\ssize \times(x+m+2j-n+4)_{n-1}\hskip15pt}(4,-2)
\Label\r{\hskip5pt\SC(x,x+m;n)}(10,-1)
\Label\r{\hskip-4pt i\!\le\!n\!-\!2}(17,1)
\Label\r{\hskip-4pt i\!=\!n\!-\!1}(17,-1)
}
\hskip7.5cm
\).
\Label\r{\hskip-12pt j\!\le\!n\!-\!2}(-14,-3)
\Label\r{j\!=\!n\!-\!1}(-6,-3)
\Label\r{\hskip8pt(3.67)}(1,-3)
\hskip1cm
\endmultline
$$
\vskip3pt
\noindent
By Step~5, we know that the degree of $P_1(x;m,n)$ is at most
$\fl{m/2}$. Hence, if we are able to determine the value of
$P_1(x;m,n)$ at $\fl{m/2}+1$ different specializations, then we can
compute $P_1(x;m,n)$ explicitly, e.g\. by Lagrange interpolation. 

The specializations that we choose are of the form $-v-1/2$, where $v$
is some nonnegative integer. The first thing to be observed is that
if we set $x=-v-1/2$, $v$ integral, in (3.67), then the denominator
on the right-hand side of (3.67) does not vanish. So, everything is
well-defined for this type of specialization. 

Next, we observe that for $x=-v-1/2$ ``usually" (this will be
specified in a moment) a lot of entries in the determinant in (3.67)
will vanish. More precisely, since $(2x+m+i+2)_j$, which is a term in
each entry of the determinant except for the $(n-1,n-1)$-entry, vanishes
if $i\le -2x-m-2=2v-m-1$ and $i+j\ge -2x-m-1=2v-m$, the determinant
takes on the form
$$\hskip-2cm
\det_{0\le i,j\le n-1}\!\!\(
\hbox{\hskip.3cm}
\PfadDicke{.3pt}
\Pfad(0,-3),1111111\endPfad
\Pfad(0,0),1111111\endPfad
\Pfad(0,4),1111111\endPfad
\Pfad(0,-3),2222222\endPfad
\Pfad(4,-3),2222222\endPfad
\Pfad(7,-3),2222222\endPfad
\Label\ro{*}(0,0)
\Label\ro{*}(0,3)
\Label\ro{*}(3,3)
\Label\ro{$\fourteenpoint$*}(1,2)
\Label\r{$\fourteenpoint$0\hskip5pt}(5,2)
\Label\ro{0\hskip10pt}(1,0)
\Label\ro{0}(3,0)
\Label\ro{\raise5pt\hbox{$0$}}(3,2)
\Label\o{$\fourteenpoint$*}(2,-2)
\Label\ro{\Cal M}(5,-2)
\Label\ro{\raise10pt\hbox{$\iddots$}}(2,1)
\Label\ro{\raise5pt\hbox{$\iddots$}}(2,2)
\Label\ro{\raise5pt\hbox{$\iddots$\hskip5pt}}(1,1)
\Label\ro{\dots}(2,0)
\Label\ro{\dots}(1,3)
\Label\ro{\dots}(2,3)
\Label\ro{\raise5pt\hbox{$\vdots$}}(0,1)
\Label\ro{\raise5pt\hbox{$\vdots$}}(0,2)
\Label\ro{\raise10pt\hbox{$\vdots$}}(3,1)
\Label\ro{\eightpoint\leftarrow i\!=\!2v\!-\!m\!-\!1}(10,0)
\hskip3.7cm
\).
\Label\r{\hskip-9pt \eightpoint j\!=\!2v\!-\!m\!-\!1}(-5,5)
\Label\ro{\hskip-9pt \eightpoint \downarrow}(-5,4)
\tag3.68
$$
Obviously, this picture makes sense only if $-1\le 2v-m-1\le n-1$, or
equivalently, if $m/2\le v\le (m+n)/2$. It should be observed that
this constraint is met by the choices
$v=\fl{(m+n)/2},\fl{(m+n)/2}-1,\dots,\fl{(m+n)/2}-\fl{m/2}$ that are
suggested in the statement of Theorem~2. In particular, for the lower
bound this is because of the assumption $m\le n$.

Because of the 0-matrix in the upper-right block of the matrix in (3.68),
it follows that the determinant in (3.68) equals the product of the
determinant of the upper-left block times $\det(\Cal M)$. Since the
upper-left block is a triangular matrix, we obtain for the
determinant in (3.68) an expression of the type
$$\multline
\text {(product of the elements along the antidiagonal
$i+j=2v-m-1$)}\\
\times \det_{(m+n-2v)\times(m+n-2v)}(\Cal M).
\endmultline$$
In the notation of the statement of the Theorem, i.e., with
$v=\fl{(m+n)/2}-t$, the dimension of $\det(\Cal M)$ is
$(m+n)-2\fl{(m+n)/2}+2t$, which is less or equal $2t+1$.

Summarizing, we have seen that for $x=-\fl{(m+n)/2}+t-1/2$, 
$t=0,1,\dots,\mathbreak\fl{m/2}$, the determinant in (3.67) reduces to a
(well-defined) multiple of a determinant of dimension at most $2t+1$. 
Since we assume $m$ to be some fixed, {\it explicit\/} nonnegative
integer, and since $2t+1\le m+1$ ($m+1$ being a fixed bound), 
this determinant can be computed
explicitly (at least in principle), and so also the explicit value of
$P_1(x;m,n)$ at $x=-\fl{(m+n)/2}+t-1/2$, 
$t=0,1,\dots,\fl{m/2}$. So, the value of $P_1(x;m,n)$ can be
computed explicitly for $\fl{m/2}+1$ distinct specializations, which
suffices to compute $P_1(x;m,n)$ explicitly by Lagrange
interpolation.

\smallskip
This finishes the proof of Theorem~2.\quad \quad \qed
\enddemo
We have used Theorem~2 to evaluate the determinant $\DetA(x,x+m;n)$
for $m=0,1,2,3,4$. This is the contents of the next Corollary.
\proclaim{Corollary 3}Let $x$ and $n$ be nonnegative integers. Then
the determinant
$$\DetA(x,x+m;n)=\det_{0\le i,j\le n-1}\bigg(\sum
 _{x+2i-j<r\le x+m+2j-i}
^{}\binom {2x+m+i+j}r\bigg)$$
for $m=0$ equals
$$\multline
\cases \dsize\prod _{i=0} ^{n-1}\(\frac {i!\,(2x+i)!\,(3x+2i+2)_i^2}
{(x+2i)!^2}\)
\frac{\prodl _{i=0}
^{n/2-1}(2x+2i+1)} {(n-1)!!}&n\text { even}\\
0&n\text { odd}
\endcases\\
=\cases \dsize\prod _{i=0} ^{n-1}\frac {(3x+2i+2)_i^2}
{(x+2i)!^2}\prod _{i=0} ^{n/2-1}\big((2i)!^2\,(2x+2i+1)!^2\big)
&n\text { even}\\
0&n\text { odd,}
\endcases
\endmultline\tag3.69$$
for $m=1$, $n\ge1$, equals
$$\prod _{i=0} ^{n-1}\(\frac
{i!\,(2x+i+1)!\,(3x+2i+3)_i\,(3x+2i+4)_i}
{(x+2i)!\,(x+2i+1)!}\)
\frac{\prodl _{i=0}
^{\fl{n/2}-1}(2x+2i+3)} {(2\fl{n/2}-1)!!},
\tag3.70$$
for $m=2$, $n\ge2$, equals
$$\multline \prod _{i=0} ^{n-1}\(\frac
{i!\,(2x+i+2)!\,(3x+2i+4)_i\,(3x+2i+6)_i}
{(x+2i)!\,(x+2i+2)!}\)
\frac{\prodl _{i=0}
^{\fl{n/2}-1}(2x+2i+3)} {(2\fl{n/2}-1)!!}\\
\times\frac {1} {(x+1)}\cdot\cases (x+n+1)&n\text { even}\\
(2x+n+2)&n\text { odd,}\endcases
\endmultline\tag3.71$$
for $m=3$, $n\ge3$, equals
$$\multline \prod _{i=0} ^{n-1}\(\frac
{i!\,(2x+i+3)!\,(3x+2i+5)_i\,(3x+2i+8)_i}
{(x+2i)!\,(x+2i+3)!}\)
\frac{\prodl _{i=0}
^{\fl{n/2}-1}(2x+2i+5)} {(2\fl{n/2}-1)!!}\\
\times\frac {1} {(x+1)}\cdot\cases (x+2n+1)&n\text { even}\\
(3x+2n+5)&n\text { odd,}\endcases
\endmultline\tag3.72$$
and for $m=4$, $n\ge4$, equals
\vskip3pt
\vbox{{}
$$\multline \prod _{i=0} ^{n-1}\(\frac
{i!\,(2x+i+4)!\,(3x+2i+6)_i\,(3x+2i+10)_i}
{(x+2i)!\,(x+2i+4)!}\)
\frac{\prodl _{i=0}
^{\fl{n/2}-1}(2x+2i+5)} {(2\fl{n/2}-1)!!}\\
\times\frac {1} {(x+1)(x+2)}\cdot\cases (x^2+(4n+3)x+2(n^2+4n+1))
&n\text { even}\\
(2x+n+4)(2x+2n+4)&n\text { odd.}\endcases
\endmultline\tag3.73$$
\line{\hfil\hbox{\qed}\quad \quad }
}
\endproclaim
At this point we remark that (3.69) combined with Theorem~1, item (3),
(2.2a),
settles the ``$n$ even" case of the Conjecture in the introduction,
see Theorem~11.

\medskip
We have computed $P_1(x;m,n)$ for further values of $m$. Together
with the cases $m=0,1,2,3,4$ that are displayed in Corollary~3, the
results suggest that actually a stronger version of Theorem~2 is
true.
\proclaim{Conjecture}Let $x,m,n$ be nonnegative integers with $m\le
n$. Then
$$\multline \DetA(x,x+m;n)=\det_{0\le i,j\le n-1}\bigg(\sum
 _{x+2i-j<r\le x+m+2j-i}
^{}\binom {2x+m+i+j}r\bigg)\\
=\prod _{i=1} ^{n-1}\(\frac {i!\,(2x+m+i)!\,(3x+m+2i+2)_i\,(3x+2m+2i+2)_i}
{(x+2i)!\,(x+m+2i)!}\)\hskip2.4cm\\
\times \frac {(2x+m)!} {(x+\fl{m/2})!\,(x+m)!}\cdot
\frac {\prodl _{i=0}
^{\fl{n/2}-1}(2x+2\cl{m/2}+2i+1)} {(2\fl{n/2}-1)!!}\cdot P_3(x;m,n),
\endmultline\tag3.74$$
where $P_3(x;m,n)$ is a polynomial in $x$ of {\rm exact} degree $\fl{m/2}$.
In addition, if the cases $n$ even and $n$ odd are considered
separately, the coefficient of $x^e$ in $P_3(x;m,n)$ is a polynomial
in $n$ of degree $\fl{m/2}-e$ with positive integer coefficients.
\endproclaim
Note that $P_3(x;m,n)=P_1(x;m,n)\cdot(2\fl{n/2}-1)!!/\prod _{i=0} ^{n-1}i!$
(compare (3.74) and (3.3)).

Possibly, this Conjecture (at least the statement about the degree of
$P_3(x;m,n)$) can be proved by examining the considerations in Step~5
and Step~6 of the proof of Theorem~2 in more detail.

\subhead 4. Another two-parameter family of determinants\endsubhead
The goal of this section is to evaluate the determinant in (2.2b). We
shall consider the generalized determinant
$$\DetB(x,y;n):=\det\limits_{0\le i,j\le n-1}\(\frac
{(x+y+i+j-1)!\,(y-x+3j-3i)}
{(x+2i-j+1)!\,(y+2j-i+1)!}\)
,\tag4.1$$
for integral $x$ and $y$. (On the side, we remark that $\DetB(x,y;n)$
would also make sense for complex $x$ and $y$ if the factorials are
interpreted as the appropriate gamma functions. Proposition~4 below, together
with its proof, actually holds in this more general sense. This
applies also to Proposition~5, as long as $m$ is a
nonnegative integer, to Corollary~6, to Theorems~8 and 9, and their
proofs.) 
$\DetB(x,y;n)$ reduces to the determinant in
(2.2b) when $n$ is replaced by $n-1$ and $y$ is set equal to $x$, 
apart from the factor $\prod _{i=0} ^{n-1}(3x+3i+4)^2$
that can be taken out of the determinant in (2.2b).

Ultimately, in Theorem~8 at the end of this section, we shall be able
to evaluate the determinant $\DetB(x,y;n)$ completely, for
independent $x$ and $y$. This is different from the determinant
$\DetA(x,y;n)$ of the previous section. But, there is a long way to go. 
The first result of this section, Proposition~4, describes how the determinant
$\DetB(x,y;n)$ factors for independent $x$ and $y$, however, leaving
one factor undetermined. It provides the ground work for
the subsequent Proposition~5 that makes it possible to evaluate
$\DetB(x,y;n)$ when the difference $m=y-x$ is fixed. This is then done
explicitly for two cases in Corollary~6. This
includes the case $m=0$ which gives the evaluation of the
determinant in (2.2b) that we
are particularly interested in. The rest of the section is then
dedicated to the complete evaluation of the determinant
$\DetB(x,y;n)$, for independent $x$ and $y$. This is finally done in
Theorem~8. Before, in Lemma~7, we collect information
about the polynomial factor $P_4(x,y;n)$ in the factorization (4.2)
of $\DetB(x,y;n)$. The proof of Theorem~8 then combines this
information with the evaluation of $\DetB(x,x+1;n)$, which is the
second case of Corollary~6.
\proclaim{Proposition 4}Let $x,y,n$ be nonnegative integers.
Then
$$\multline \DetB(x,y;n)=\det_{0\le i,j\le n-1}\(
\frac
{(x+y+i+j-1)!\,(y-x+3j-3i)}
{(x+2i-j+1)!\,(y+2j-i+1)!}\)\\
=\prod _{i=0} ^{n-1}\(\frac {(x+y+i-1)!\,(2x+y+2i+1)_i\,(x+2y+2i+1)_i}
{(x+2i+1)!\,(y+2i+1)!}\)\cdot P_4(x,y;n),
\endmultline\tag4.2$$
where $P_4(x,y;n)$ is a polynomial in $x$ and $y$ of degree $n$.
\endproclaim

\demo{Proof} Again, the proof is divided into several steps. The strategy is
very similar to the proof of Theorem~2. First, we transform $\DetB(x,x+m;n)$ 
into a multiple of another determinant,
namely $\DetBb(x,y;n)$, by (4.3)--(4.5), which is a polynomial in $x$ and $y$, 
then identify as many 
factors of the new determinant as possible
(as a polynomial in $x$ and $y$), and finally find a bound
for the degree of the remaining polynomial factor.

\smallskip
{\it Step 1. An equivalent statement of the Theorem.}
We take as many common factors out
of the $i$-th row of $\DetB(x,y;n)$, 
$i=0,1,\dots,n-1$, as possible, such that the
entries become polynomials in $x$ and $y$. To be precise, 
we take
$$\frac {(x+y+i-1)!} {(x+2i+1)!\,(y+2n-i-1)!}$$
out of the $i$-th row, $i=0,1,\dots,n-1$. This gives
$$\multline
\DetB(x,y;n)=
\prod _{i=0} ^{n-1}\frac {(x+y+i-1)!} 
{(x+2i+1)!\,(y+2n-i-1)!}\\
\times\det_{0\le i,j\le
n-1}\big((x+y+i)_{j}\,(x+2i-j+2)_j\,(y+2j-i+2)_{2n-2j-2}\,
(y-x+3j-3i)\big).
\endmultline\tag4.3$$
For convenience, let us denote the determinant in (4.3) by
$\DetBa(x,y;n)$. 
In fact, there are more factors that can be taken out of $\DetBa(x,y;n)$ 
under the
restriction that the entries of the determinant continue to be
polynomials. To this end, we multiply the $i$-th row of $\DetBa(x,y;n)$
by $(y+2n-i)_{i}$, $i=0,1,\dots,n-1$, and divide the $j$-th column by
$(y+2j+2)_{2n-2j-2}$, $j=0,1,\dots,n-1$.
This leads to
$$\multline \prodl _{i=0} ^{n-1}(y+2n-i)_{i}\prodl
_{j=0} ^{n-1}\frac {1} {(y+2j+2)_{2n-2j-2}}\cdot\DetBa(x,y;n)
\\
=\det_{0\le i,j\le
n-1}\big((x+y+i)_{j}\,(x+2i-j+2)_j\,(y+2j-i+2)_{i}\,
(y-x+3j-3i)\big),
\endmultline\tag4.4$$
or, if we denote the determinant in (4.4) by $\DetBb(x,y;n)$,
$$\DetBa(x,y;n)=\prod _{i=0} ^{n-1}(y+2i+2)_{n-i-1}
\cdot\DetBb(x,y;n).\tag4.5$$
A combination of (4.2), (4.4), and (4.5) then implies that
Proposition~4
is equivalent to the statement: 

{\sl With $\DetBb(x,y;n)$ the determinant in {\rm(4.4)},
there holds
$$\multline
\DetBb(x,y;n)=\prod _{i=0}
^{n-1}\big((2x+y+2i+1)_{i}\,(x+2y+2i+1)_{i}\big)\cdot P_4(x,y;n),
\endmultline\tag4.6$$
where $P_4(x,y;n)$ is a polynomial in $x$ and $y$ of degree $n$}.

\smallskip
{\it Step 2. $\prod _{i=0}
^{n-1}\big((2x+y+2i+1)_{i}\,(x+2y+2i+1)_{i}\big)$
is a factor of $\DetBb(x,y;n)$.}
There are not many differences to Step~3 of the proof of Theorem~2.
So we shall be brief here.

We first consider just one half of this product, 
$\prod _{i=0}
^{n-1}(2x+y+2i+1)_{i}$. Let us concentrate on a typical factor
$(2x+y+2i+l+1)$, $0\le i\le n-1$, $0\le l<i$. We claim that 
for each
such factor there is a linear combination of the rows that vanishes
if the factor vanishes. More precisely, we claim that for any $i,l$
with $0\le i\le n-1$, $0\le l<i$ there holds
$$\multline \sum _{s=l} ^{\fl{(i+l)/2}}\frac
{(2i-3s+l)} {(i-s)}\frac {(i-2s+l+1)_{s-l}} {(s-l)!}\frac
{(x+2s+2)_{2i-2s}} {(-x-2i-l+s-1)_{i-s}}\\
\cdot(\text {row $s$ of
$\DetBb(x,-2x-2i-l-1;n)$})
=(\text {row $i$ of $\DetBb(x,-2x-2i-l-1;n)$}).
\endmultline\tag4.7$$
Restricting (4.7) to the $j$-th column, it is seen that this means to
check
$$\multline \sum _{s=l} ^{\fl{(i+l)/2}}\frac
{(2i-3s+l)} {(i-s)}\frac {(i-2s+l+1)_{s-l}} {(s-l)!}\frac
{(x+2s+2)_{2i-2s}} {(-x-2i-l+s-1)_{i-s}}\hskip1cm\\
\times(-x-2i-l+s-1)_{j}\,(x+2s-j+2)_{j}\,\\
\times(-2x-2i-l+2j-s+1)_{s}\,(-3x-2i-l+3j-3s-1)\\
=(-x-i-l-1)_{j}\,(x+2i-j+2)_{j}\,
(-2x-3i-l+2j+1)_{i}\,(-3x-l+3j-5i-1).
\endmultline$$
This is easily done by observing that it is equivalent to (3.16).
Arguments that are similar to those after (3.18) then show that 
the complete product $\prod _{i=0}
^{n-1}(2x+y+2i+1)_{i}$ divides $\DetBb(x,y;n)$.

The reasoning that $\prod _{i=0}
^{n-1}(x+2y+2i+1)_{i}$ is a factor of $\DetBb(x,y;n)$ is
similar. Also here, let us concentrate on a typical factor
$(x+2y+2j+l+1)$, $0\le j\le n-1$, $0\le l<j$. This time
we claim that for each
such factor there is a linear combination of the columns that vanishes
if the factor vanishes. More precisely, we claim that for any $j,l$
with $0\le j\le n-1$, $0\le l<j$ there holds
$$\multline \sum _{s=l} ^{\fl{(j+l)/2}}\frac
{(2j-3s+l)} {(j-s)}\frac {(j-2s+l+1)_{s-l}} {(s-l)!}
(y+2s+2)_{2j-2s}\\
\cdot(\text {column $s$ of
$\DetBb(-2y-2j-l-1,y;n)$})\\
=(\text {column $j$ of $\DetBb(-2y-2j-l-1,y;n)$}).
\endmultline\tag4.8$$
Restricting to the $i$-th row, we see that this means to check
$$\multline \sum _{s=l} ^{\fl{(j+l)/2}}\frac
{(2j-3s+l)} {(j-s)}\frac {(j-2s+l+1)_{s-l}} {(s-l)!}
(y+2s+2)_{2j-2s}\\
\times(-y-2j-l+i-1)_{s}\,(-2y-2j-l+2i-s+1)_{s}\,\\
\times(y+2s-i+2)_{i}\,(3y+2j+l-3i+3s+1)\\
=(-y-2j-l+i-1)_{j}\,(-2y-3j-l+2i+1)_{j}\,(y+2j-i+2)_{i}\,
(3y+5j+l-3i+1).
\endmultline$$
The observation that this summation is equivalent to (3.24) with
$\y=y$ establishes the claim. Similarly to as before, this eventually
shows that the complete product $\prod _{i=0}
^{n-1}(x+2y+2i+1)_{i}$ divides $\DetBb(x,y;n)$. 

Altogether, this implies that $\prod _{i=0}
^{n-1}\big((2x+y+2i+1)_{i}
(x+2y+2i+1)_{i}\big)$ divides $\DetBb(x,y;n)$,
as desired.

\smallskip
{\it Step 3. $P_4(x,y;n)$ is a polynomial in $x$ and $y$ of degree
$n$.}
We shall prove that $\DetBa(x,y;n)$ (which is defined to be the
determinant in (4.3)) is a polynomial in $x$ and $y$ of (total) degree 
$3\binom n2+n$.
By (4.5) this would
imply that $\DetBb(x,y;n)$ is a polynomial in $x$ and $y$ of degree 
$2\binom {n}2+n$, and so, by (4.6), that
$P_4(x,y;n)$ is a polynomial in $x$ and $y$ of degree $n$, as
desired.

Here we need to consider the generalized determinant
$$\DetBaa(x,y,z(0),z(1),\dots,z(n-1);n)=\DetBaa(n)$$
which arises from $\DetBa(x,y;n)$ by replacing each occurence of $i$
in row $i$ by an indeterminate, $z(i)$ say, $i=0,1,\dots,n-1$,
$$\multline
\DetBaa(x,y,z(0),z(1),\dots,z(n-1);n)=\DetBaa(n)\\
:=
\det_{0\le i,j\le
n-1}\big((x+y+z(i))_{j}\,(x+2z(i)-j+2)_j\hskip3cm\\
\cdot(y+2j-z(i)+2)_{2n-2j-2}\,
(y-x+3j-3z(i))\big).
\endmultline\tag4.9$$
This determinant is a polynomial in $x, y,z(0),z(1),\dots,z(n-1)$. We
shall prove that the degree in $x$ and $y$ of this determinant is 
$3\binom n2+n$, 
which clearly implies our claim upon setting $z(i)=i$,
$i=0,1,\dots,n-1$.

Obviously, the total degree of $\DetBaa(n)$ in
$x,y,z(0),z(1),\dots,z(n-1)$ is at most $4\binom n2+n$. However,
actually it is {\it exactly\/} equal to this upper bound, since 
the monomial 
$$y^{3\binom n2+n}z(0)^0z(1)^1\cdots z(n-1)^{n-1}$$
occurs only in the product of the main diagonal of the determinant
with nonzero coefficient, and therefore cannot cancel. On the other
hand, when $z(i_1)=z(i_2)$ for some $i_1\ne i_2$, the rows $i_1$ and
$i_2$ in $\DetBaa(n)$ are identical. Hence $\DetBaa(n)$ vanishes in
this case. This shows that
the product $\prod _{0\le i<j\le n-1} ^{}(z(j)-z(i))$ divides
$\DetBaa(n)$. Therefore, the degree in $x$ and $y$ of $\DetBaa(n)$
equals $\(4\binom n2+n\)-\binom n2=3\binom n2+n$, which 
is what we need.

\smallskip
This finishes the proof of Proposition~4.\quad \quad \qed
\enddemo

As we already mentioned, it is possible to compute the polynomial
$P_4(x,y;n)$ explicitly, see Theorem~8. However, we are not yet in
the position to do so. First we restrict ourselves to the situation
that was studied in Section~3 for the determinant (3.1), namely when
the difference $m=y-x$ is a fixed integer. In this situation, 
it turns out that the
degree of $P_4(x,y;n)=P_4(x,x+m;n)$, now as a polynomial in $x$ of course,
shrinks significantly. To be precise, the degree can be at most $\fl{n/2}$. 
We prove this fact in Proposition~5 below. Moreover, we identify
several factors of $P_4(x,x+m;n)$, which will be of great help in the
proof of Corollary~6.
\proclaim{Proposition 5}Let $x,m,n$ be nonnegative integers with $m\le
n$. Then, 
$$\multline \DetB(x,x+m;n)=\det_{0\le i,j\le n-1}\(
\frac
{(2x+m+i+j-1)!\,(m+3j-3i)}
{(x+2i-j+1)!\,(x+m+2j-i+1)!}\)\\
=\prod _{i=0} ^{n-1}\(\frac {(2x+m+i-1)!\,(3x+m+2i+1)_i\,(3x+2m+2i+1)_i}
{(x+2i+1)!\,(x+m+2i+1)!}\)\\
\times(x+m)_{\fl{n/2}-\fl{m/2}}\cdot P_5(x;m,n),
\endmultline\tag4.10$$
where $P_5(x;m,n)$ is a polynomial in $x$ of degree $\le \fl{m/2}$.
\endproclaim

\demo{Proof}Much of the required work has already been done. In
particular, if we set $y=x+m$ in (4.2) and compare with (4.10), we
obtain
$$P_4(x,x+m;n)=(x+m)_{\fl{n/2}-\fl{m/2}}\cdot P_5(x;m,n).\tag4.11$$
So what remains to prove is that $(x+m)_{\fl{n/2}-\fl{m/2}}$ is a factor
of $P_4(x,x+m;n)$, 
and that the degree of $P_4(x,x+m;n)$ is at most
$\fl{n/2}$.
The first fact is established in Step~1, the second in
Step~2.

\smallskip
{\it Step 1. $(x+m)_{\fl{n/2}-\fl{m/2}}$ is a factor of $P_4(x,x+m;n)$.} 
Fix an integer $v$, 
$$m\le v\le \frac {m+n-1} {2}.\tag4.12$$
We shall show that $x+v$ is a factor of $P_4(x,x+m;n)$, or
equivalently, that $x=-v$ is a zero of $P_4(x,x+m;n)$. The polynomial
$P_4(x,x+m;n)$ is defined by (4.6), with $y=x+m$, i.e., by
$$\multline
\DetBb(x,x+m;n)=\prod _{i=0}
^{n-1}\big((3x+m+2i+1)_{i}\,(3x+2m+2i+1)_{i}\big)\cdot P_4(x,x+m;n),
\endmultline\tag4.13$$
where $\DetBb(x,x+m;n)$ is the determinant in (4.4) with $y=x+m$. So,
we would like to set $x=-v$ in (4.13), prove that
$\DetBb(-v,-v+m;n)$ equals 0, that the product on the right-hand side of
(4.13) is nonzero, and conclude that therefore $P_4(-v,-v+m;n)$ must
be 0. However, the product on the right-hand side of (4.13)
unfortunately (usually) {\it is\/} 0 for $x=-v$. This makes it necessary
to cancel first all factors $(x+v)$ that occur in the product, and
only then set $x=-v$.

To accomplish this, we have to ``generate" these factors on the
left-hand side. This is done by reading through Step~2 of the proof
of Proposition~4 with $y=x+m$. To make this more precise, observe that
$x+v$ divides a typical factor $3x+m+2i+l+1$, $0\le i\le n-1$, $0\le
l<i$, of the first half of the product in (4.13) if and only if
$3v=m+2i+l+1$. So, for each solution $(i,l)$ of
$$3v=m+2i+l+1,\quad  \text {with }0\le i\le n-1,\ 0\le l<i,\tag4.14$$
we subtract the linear combination
$$\multline \sum _{s=l} ^{\fl{(i+l)/2}}\frac
{(2i-3s+l)} {(i-s)}\frac {(i-2s+l+1)_{s-l}} {(s-l)!}\frac
{(x+2s+2)_{2i-2s}} {(-x-2i-l+s-1)_{i-s}}\\
\cdot(\text {row $s$ of
$\DetBb(x,x+m;n)$})
\endmultline\tag4.15$$
of rows of $\DetBb(x,x+m;n)$ from row $i$ of $\DetBb(x,x+m;n)$. Let
us denote the resulting determinant by $\DetBt(x;m,n)$. By (4.7), the
effect is that $(3x+m+2i+l+1)=(3x+3v)$ (the equality being due to
(4.14)), and hence $(x+v)$, is a factor of each
entry of the $i$-th row of $\DetBt(x;m,n)$, for each solution $(i,l)$
of (4.14). For later use we record that the $(i,j)$-entry of
$\DetBt(x;m,n)$, $(i,l)$ a solution of (4.14), reads
$$\multline (2x+m+i)_{j}\,(x+2i-j+2)_j\,(x+m+2j-i+2)_{i}\,
(m+3j-3i)\\
-\sum _{s=l} ^{\fl{(i+l)/2}}\frac
{(2i-3s+l)} {(i-s)}\frac {(i-2s+l+1)_{s-l}} {(s-l)!}\frac
{(x+2s+2)_{2i-2s}} {(-x-2i-l+s-1)_{i-s}}\hskip2cm\\
\times(2x+m+s)_{j}\,(x+2s-j+2)_j\,(x+m+2j-s+2)_{s}\,
(m+3j-3s).\endmultline\tag4.16$$
Similar considerations concern the second half of the product in
(4.13). Omitting the details, for each solution $(j,l)$ of
$$3v=2m+2j+l+1,\quad  \text {with }0\le j\le n-1,\ 0\le l<j,\tag4.17$$
we subtract the linear combination
$$\multline \sum _{s=l} ^{\fl{(j+l)/2}}\frac
{(2j-3s+l)} {(j-s)}\frac {(j-2s+l+1)_{s-l}} {(s-l)!}
(x+m+2s+2)_{2j-2s}\\
\cdot(\text {column $s$ of
$\DetBt(x;m,n)$})
\endmultline$$
of columns of $\DetBt(x;m,n)$ (we definitely mean $\DetBt(x;m,n)$,
and not $\DetBb(x;m,n)$) from column $j$ of $\DetBt(x,x+m;n)$. By
(4.8), each entry of the $j$-th column of the new determinant will
have $(x+v)$ as a factor. We remark that entries that were changed by
a row {\it and\/} column operations will now have $(x+v)^2$ as a
factor. Now we take $(x+v)$ out of the $i$-th row, for each solution
$(i,l)$ of (4.14), and we take $(x+v)$ out of the $j$-th column, 
for each solution $(j,l)$ of (4.17). We denote the resulting
determinant by $\DetBbb(x;m,n)$. Thus, we have
$$\DetBb(x,x+m;n)=(x+v)^{\#\text {(solutions $(i,l)$ of (4.14))}+
\#\text {(solutions $(j,l)$ of (4.17))}}\DetBbb(x;m,n).$$
Plugging this into (4.13), we see that now all factors $(x+v)$ can be
cancelled on both sides, so that we obtain
$$\DetBbb(x;m,n)=C(x;m,n)\,P_4(x,x+m;n),$$
for some $C(x;m,n)$ that does not vanish for $x=-v$. Hence, if we are
able to prove that $\DetBbb(-v;m,n)=0$, it would follow that
$P_4(-v,-v+m;n)=0$, which is what we want to establish.

So we are left with showing that $\DetBbb(-v;m,n)=0$. This will be
implied by the following two claims: The matrix of which $\DetBbb(-v;m,n)$ is the
determinant has a block form (see (4.18)), where 

Claim 1. the
upper-right block, consisting of the entries that are in one of the
rows $0,1,\dots,2v-m$ and one of the columns
$2v-m+1,2v-m+2,\dots,n-1$, is a zero matrix, and where 

Claim 2. the
determinant of the upper-left block, $\Cal N$, 
consisting of the entries that are in one of the
rows $0,1,\dots,2v-m$ and one of the columns
$0,1,\dots,2v-m$, equals 0.
\vskip3pt
$$
\PfadDicke{.3pt}
\Pfad(0,-3),11111\endPfad
\Pfad(0,0),11111\endPfad
\Pfad(0,2),11111\endPfad
\Pfad(0,-3),22222\endPfad
\Pfad(2,-3),22222\endPfad
\Pfad(5,-3),22222\endPfad
\Label\r{\Cal N\hskip15pt}(1,1)
\Label\r{$\fourteenpoint$0\hskip5pt}(3,1)
\Label\o{$\fourteenpoint$*}(1,-2)
\Label\ro{$\fourteenpoint$*}(3,-2)
\Label\r{\eightpoint\hskip4pt\leftarrow i\!=\!2v\!-\!m}(6,0)
\hbox{\hskip3.7cm}
\Label\r{\hskip-18pt \eightpoint j\!=\!2v\!-\!m}(-5,3)
\Label\ru{\hskip-25pt \eightpoint \downarrow}(-5,3)
\tag4.18
$$
\vskip3pt
For, the determinant of a block matrix of the form (4.18) equals the
product of the determinants of the upper-left block and the
lower-right block, the first determinant being equal to $0$ by
Claim~2. 

Claim~1 is most obvious for all the entries that did not change in
the transition from $\DetBb(x,x+m;n)$ to $\DetBbb(x;m,n)$. For, the
$(i,j)$-entry of $\DetBb(x,x+m;n)$, by its definition in (4.4), is
$$(-2v+m+i)_{j}\,(-v+2i-j+2)_j\,(-v+m+2j-i+2)_{i}\,
(m+3j-3i).\tag4.19$$
Clearly, if $0\le i\le 2v-m$ and $2v-m+1\le j\le n-1$, we have
$(-2v+m+i)_j=0$, and so the complete expression in (4.19) is 0. 

On the
other hand, let us consider an $(i,j)$-entry of $\DetBbb(x;m,n)$
that changed in
the transition from $\DetBb(x,x+m;n)$ to $\DetBbb(x;m,n)$. First we
want to know, where such an entry could be located. If it changed
under a row operation, then $(i,l)$ is a solution of (4.14), for some
$l$. By (4.14) we have
$$m+2i+1\le m+2i+l+1=3v\quad \text {and}\quad 3v=m+2i+l+1\le m+3i,$$
and so,
$$v-\frac {m} {3}\le i\le \frac {3v-m-1} {2}.\tag4.20$$
If the $(i,j)$-entry 
changed
under a column operation, then $(j,l)$ is a solution of (4.17), for some
$l$. Similar arguments then give, using (4.17), that
$$v-\frac {2m} {3}\le j\le \frac {3v-2m-1} {2}.\tag4.21$$
In particular we have $j<2v-m$, so an $(i,j)$-entry that is located
in the upper-right block, which we are currently interested in, did
not change under a column operation.  

But it could have changed under a row operation. Such an
$(i,j)$-entry is given by (4.16) divided by $(x+v)$. (Recall that
(4.16) was the expression for an $(i,j)$-entry that changed under a
row operation {\it before\/} we factored $(x+v)$ out of the $i$-th row.)
Thus, it can be written as
$$\multline \frac {(2x+m+i)_{2v-m-i+1}} {(x+v)}\bigg(
(2x+2v+1)_{i+j-2v+m-1}\,(x+2i-j+2)_j\\
\cdot (x+m+2j-i+2)_{i}\,(m+3j-3i)
-\sum _{s=l} ^{\fl{(i+l)/2}}\frac
{(2i-3s+l)} {(i-s)}\frac {(i-2s+l+1)_{s-l}} {(s-l)!}\\
\cdot\frac
{(x+2s+2)_{2i-2s}} {(-x-2i-l+s-1)_{i-s}}\,
(2x+m+s)_{i-s}\,(2x+2v+1)_{j+s-2v+m-1}\,\\
(x+2s-j+2)_j\,(x+m+2j-s+2)_{s}\,
(m+3j-3s)\bigg).
\endmultline\tag4.22$$
We have to show that (4.22) vanishes for $x\to-v$. Because of the
denominators, it is not even evident that (4.22) is well-defined when
$x\to-v$. However, by (4.20) we have $2v-m-i\ge(v-m+1)/2\ge0$, the
last inequality being due to our assumption $v\ge m$. Hence,
$$\frac {(2x+m+i)_{2v-m-i+1}} {(x+v)}=(2x+m+i)_{2v-m-i}\cdot 2,$$
and so the first term in (4.22) is well-defined when $x\to-v$.
Furthermore, the denominator in the sum in (4.22) (neglecting the
terms that do not depend on $x$) when $x\to-v$
becomes
$$(v-2i-l+s-1)_{i-s}=(v-2i-l+s-1)\cdots(v-i-l-2).\tag4.23$$
By (4.14) and (4.20) we have
$v-i-l-2=-2v+i+m-1\le \frac {1} {2}(-v+m-3)<0,$
again the last inequality being due to our assumption $v\ge m$.
Therefore, all the terms in (4.23) are nonzero, which means that the
denominator in the sum in (4.22) is nonzero when $x\to-v$. Hence,
(4.22) is well-defined for $x\to-v$. To demonstrate that it actually
vanishes for $x\to-v$, we show that the second term in (4.22) (the
term in big parentheses) equals 0 for $x=-v$. 

To see this, set $x=-v$, and by (4.14) replace $l$ by $3v-2i-m-1$ in
the sum (4.22), and then convert it into hypergeometric notation, to
obtain
$$\multline 2 \left( 3 + 6 i + 3 j + 4 m - 9 v \right)  
    ({ \textstyle 1}) _{v+j - 2 i-2}  \\
 \times ({ \textstyle -4v+ 2 i + 2 j + 2 m+3}) _{ 3v- 2 i - m -1}  
  ({ \textstyle 5v-4 i - j - 2 m }) _{ -6v+6 i + j + 2 m +2} \\
\times {} _{6} F _{5} \!\left [ \matrix { -2 i - j
   - {{4 m}\over 3} + 3 v, {1\over 3} - 2 i - {{2 m}\over 3} + 2 v,
    -{1\over 2} -
   {{3 i}\over 2} - {m\over 2} + {{3 v}\over 2},}\\ { -1 - 2 i - j -
   {{4 m}\over 3} + 3 v, -{2\over 3} -
   2 i - {{2 m}\over 3} + 2 v, -3 i - m + 3 v, 
   }
   \endmatrix\right.\\
   \left.\matrix  -{{3 i}\over 2} - {m\over
   2} + {{3 v}\over 2},  
   -1 - 2 i + j + v,  -2 - 2 i - 2 j - 2 m + 4 v\\
   -2 i - {j\over 2} - m + {{5 v}\over 2},
   {1\over 2} - 2 i - {j\over 2} - m + {{5 v}\over 2}\endmatrix ; 
   {\displaystyle   1}\right ].
\endmultline\tag4.24$$
The $_6F_5$-series can be summed by means of Lemma~A6. Then, after
simplification, (4.24) becomes
$$(i+j-2v+m-1)!\,(-v+2i-j+2)_j\,
 (-v+m+2j-i+2)_{i}\,(m+3j-3i),
$$
which is exactly the first term in big parentheses in (4.22) for
$x=-v$.
Therefore, the term in big parentheses in (4.22) vanishes for $x=-v$.
This settles Claim~1.

Next we turn to Claim~2. We have to prove that the determinant of the
matrix $\Cal N$, consisting of the entries of $\DetBbb(-v;m,n)$ that
are in one of the rows $0,1,\dots,2v-m$ and one of the columns
$0,1,\dots,2v-m$ (recall (4.18)), equals 0. We do this by locating
enough zeros in the matrix $\Cal N$.

We concentrate on the entries that did not change in the transition
from $\DetBb(x,x+m;n)$ to $\DetBbb(x;m,n)$. 
For the location of the various regions in the matrix $\Cal N$
that we are going to describe, always consult Figure~2 which 
gives a rough sketch.
\vskip10pt
\vbox{
$$\Cal N\ =\ 
\raise2.25cm\hbox{$
\Einheit.25cm
\thinlines
\PfadDicke{.35pt}
\Pfad(0,0),111111111111111111\endPfad
\Pfad(0,-6),111111111111111111\endPfad
\Pfad(0,-10),111111111111111111\endPfad
\Pfad(0,-15),111111111111111111\endPfad
\Pfad(0,-18),111111111111111111\endPfad
\Pfad(0,-18),222222222222222222\endPfad
\Pfad(3,-18),222222222222222222\endPfad
\Pfad(8,-18),222222222222222222\endPfad
\Pfad(12,-18),222222222222222222\endPfad
\Pfad(18,-18),222222222222222222\endPfad
\catcode`\@=11
\hskip3\Einheit\hbox to0pt{\Line@(1,-2)9\hss}\hskip-3\Einheit
\raise-6\Einheit\hbox to0pt{\Line@(2,-1){18}\hss}
\catcode`\@=13
\Diagonale(0,-18){18}
\Vektor(-1,0)8(23,-3)
\Vektor(-1,0)9(23,-13)
\Vektor(1,0)6(-2,-2)
\Label\r{\eightpoint i=0}(19,0)
\Label\r{\eightpoint i=\cl{\frac {v-1} {2}}}(20,-6)
\Label\r{\eightpoint\hskip8pt i=\cl{v-\frac {m} {3}}-2+\chi(m\equiv 0(3))}(25,-10)
\Label\r{\eightpoint\hskip17pt i=\fl{\frac {3v-m-1} {2}}+1}(21,-15)
\Label\r{\eightpoint\hskip3pt i=2v-m}(20,-18)
\Label\o{\eightpoint 0}(0,0)
\Label\o{\eightpoint \raise3pt\hbox{$j=$}}(0,1)
\Label\u{\eightpoint j=}(3,-18)
\Label\u{\eightpoint \raise-10pt\hbox{$\cl{\frac {v-m-1} {2}}$}}(3,-19)
\Label\o{\eightpoint \cl{v-\frac {2m} {3}}-1}(8,0)
\Label\o{\eightpoint \raise10pt\hbox{$j=$}}(8,1)
\Label\u{\eightpoint j=}(12,-18)
\Label\u{\eightpoint \raise-10pt\hbox{$\fl{\frac {3v-2m-1} {2}}+1$}}(12,-19)
\Label\o{\eightpoint 2v-m}(18,0)
\Label\o{\eightpoint \raise5pt\hbox{$j=$}}(18,1)
\Label\r{\eightpoint\hskip-12pt i=\frac {v+j-2} {2}}(26,-13)
\Label\r{\eightpoint j=\frac {v-m+i-2} {2}}(-7,-2)
\Label\r{\eightpoint\hskip5pt i+j=2v-m}(26,-3)
\Label\r{$\hskip8pt I$}(4,-8)
\Label\r{$II\hskip2pt $}(14,-8)
\Label\ro{$\hskip10pt IV$}(14,-17)
\Label\ro{$III$}(5,-17)
\hskip7cm
$}
$$
\vskip6pt
\centerline{\eightpoint Figure 2}
}
\vskip10pt
By earlier
considerations, an $(i,j)$-entry did not change if $i$ is outside
the range (4.20), i.e.,
$$0\le i\le \cl{v-\frac {m} {3}}-1\quad \text {or}\quad \fl{\frac
{3v-m-1} {2}}+1\le i\le n-1,\tag4.25$$
and if $j$ is outside
the range (4.21), i.e.,
$$0\le j\le \cl{v-\frac {2m} {3}}-1\quad \text {or}\quad \fl{\frac
{3v-2m-1} {2}}+1\le j\le n-1.\tag4.26$$
As we already noted, such an $(i,j)$-entry is given by (4.19).
The first term in (4.19) vanishes if and only if
$$i\le 2v-m\quad \text {and}\quad i+j>2v-m.\tag4.27$$
The second term in (4.19) vanishes if and only if
$$\cl{\frac {v-1} {2}}\le i\le \frac {v+j-2} {2}.\tag4.28$$
The third term in (4.19) vanishes if and only if
$$\cl{\frac {v-m-1} {2}}\le j\le \frac {v-m+i-2} {2}.\tag4.29$$
Finally, the fourth term in (4.19) vanishes if and only if
$$m\equiv 0\pmod3\quad \text {and}\quad i=j+\frac {m} {3}.\tag4.30$$

Now we claim that in the following four regions of $\Cal N$ all the
entries are 0, except for the case $v=m=1$, which we treat
separately. Again, to get an idea of the location of these regions,
consult Figure~2.

\smallskip
Region I: All $(i,j)$-entries with
$$\multline 
\cl{\frac {v-1} {2}}\le i\le \cl{v-\frac {m} {3}}-2+\chi(m\equiv
0\pmod 3)\\
\quad \text {and}\quad \cl{\frac {v-m-1} {2}}\le j\le 
\cl{v-\frac {2m} {3}}-1,
\endmultline\tag4.31$$
where again $\chi(\Cal A)$=1 if $\Cal A$ is
true and $\chi(\Cal A)$=0 otherwise.

Region II: All $(i,j)$-entries with
$$\multline
\cl{\frac {v-1} {2}}\le i\le \cl{v-\frac {m} {3}}-2+\chi(m\equiv
0\pmod 3)\\
\quad \text {and}\quad \fl{\frac {3v-2m-1} {2}}+1\le j\le
2v-m.
\endmultline\tag4.32$$

Region III: All $(i,j)$-entries with
$$\fl{\frac {3v-m-1} {2}}+1\le i\le 2v-m\quad \text {and}\quad 
\cl{\frac {v-m-1} {2}}\le j\le 
\cl{v-\frac {2m} {3}}-1.\tag4.33$$

Region IV: All $(i,j)$-entries with
$$\fl{\frac {3v-m-1} {2}}+1\le i\le 2v-m\quad \text {and}\quad 
\fl{\frac {3v-2m-1} {2}}+1\le j\le
2v-m.\tag4.34$$

\smallskip
Instantly we observe that all four regions satisfy (4.25) and (4.26).
So, all the entries in these regions are given by (4.19). Hence, to
verify that all these entries are 0 we have to show that for each
entry one of (4.27)--(4.30) is true. Of course, we treat the four
regions separately.

\smallskip
ad Region I. First let $i\le \cl{v-m/3}-2$. In case that $i\le
j+m/3$, we have
$$i\le \frac {i+j+\frac {m} {3}} {2}\le \frac {\cl{v-\frac {m}
{3}}-2+j+\frac {m} {3}} {2}\le \frac {v-\frac {m} {3}+\frac {2}
{3}-2+j+\frac {m} {3}} {2}=\frac {v+j-2} {2}+\frac {1} {3}.$$
Combined with (4.31), this implies that (4.28) is satisfied. On the
other hand, in case that $i>j+m/3$, or equivalently, 
$$i\ge j+\frac {m} {3}+\frac {1} {3},\tag4.35$$
we have, using the last inequality in (4.31),
$$\multline
j\le \frac {i+j-\frac {m} {3}-\frac {1} {3}} {2}\le 
\frac {i+\cl{v-\frac {2m} {3}}-1-\frac {m} {3}-\frac {1} {3}} {2}\\
\le\frac {i+v-\frac {2m} {3}+\frac {2} {3}-1-\frac {m} {3}-\frac {1} {3}}
{2}=\frac {v-m+i-2} {2}+\frac {2} {3}.
\endmultline$$
Combined with (4.31), this implies that (4.29) is satisfied, unless
$j=(v-m+i-1)/2$. But if we plug this into (4.35), we obtain $i\ge
v-m/3-1/3$, a contradiction to our assumption $i\le \cl{v-m/3}-2$.

Collecting our results so far, we have seen that if
$m\equiv1,2\pmod3$, then each $(i,j)$-entry in region~I satisfies
(4.28) or (4.29). If $m\equiv0\pmod3$, region~I also contains entries
from row $i=v-m/3-1$. First let $j\le v-(2m)/3-2$. Then it is
immediate that (4.29) is satisfied. If $j=v-(2m)/3-1$, then (4.30) is
satisfied. 
This shows that if $m\equiv0\pmod3$ then an $(i,j)$-entry in region~I
satisfies (4.28), (4.29), or (4.30).

ad Region II. Here, by (4.32), we have
$$i+j\ge \cl{\frac {v-1} {2}}+\fl{\frac {3v-2m-1} {2}}+1=2v-m.$$
Hence, (4.27) is satisfied, except when $i=\cl{(v-1)/2}$ and
$j=\fl{(3v-2m-1)/2}+1$. But in that case there holds (4.28), apart
from a few exceptional cases. For, if $v\ne0,2$ (and $v\ge0$ of
course) then
$$\cl{\frac {v-1} {2}}\le \frac {v+\fl{\frac {v-1} {2}}-1} {2}.$$
By the assumption $v\ge m$ it follows that
$$\cl{\frac {v-1} {2}}\le \frac {v+\(\fl{\frac {3v-2m-1} {2}}+1\)-2}
{2},\tag4.36$$
which is nothing but (4.28) with the current choices of $i$ and $j$.
Thus, (4.28) is satisfied except when $v=m=0$, or $v=2$ and
$m=0,1,2$. (There are no more cases because $v\ge m$.) Starting 
from the back, the case $v=m=2$ does not bother us, since in that
case region~II is empty (there is no $i$ satisfying (4.32)).
By inspection, it is seen that (4.36), and hence (4.28), also holds
if $v=2$ and $m=0$ or $1$. Finally, in case $v=m=0$ we have
$i=\cl{(v-1)/2}=0$ and $j=\fl{(3v-2m-1)/2}+1=0$. Hence, (4.30) is
satisfied.

ad Region III. By (4.33) we have
$$i+j\ge \fl{\frac {3v-m-1} {2}}+1+\cl{\frac {v-m-1} {2}}=2v-m.$$
Hence again, (4.27) is satisfied, except when $i=\fl{(3v-m-1)/2}+1$
and $j=\cl{(v-m-1)/2}$. In that case there holds (4.29), apart from a
few special cases. For, if $w$ is a positive integer, then
$$\cl{\frac {w-1} {2}}\le \frac {w+\fl{\frac {w-1} {2}}} {2}.$$
Setting $w=v-m$ in this inequality we obtain for $v>m$ and $v\ge1$
the inequality
$$\cl{\frac {v-m-1} {2}}\le \frac {v-m+\(\fl{\frac {3v-m-1}
{2}}+1\)-2} {2}.\tag4.37$$
This is exactly (4.29) with the current choices for $i$ and $j$.
Thus, (4.29) is satisfied except when $v=m$ or $v=0$. (Recall that
there are no more cases because of $v\ge m$.) But (4.37), and hence
(4.29), holds in more cases. Namely, by inspection, if $v=m$, then
(4.37) holds for $v\ge2$.
So, the only cases in which (4.37) is not true are $v=m=0$ and
$v=m=1$. The case $v=m=0$ does not bother us, since in that case
region~III is empty (there is no $j$ satisfying (4.33)). The case
$v=m=1$ is the exceptional case that is treated separately.

ad Region~IV. By (4.34) we have
$$i+j\ge \fl{\frac {3v-m-1} {2}}+1+\fl{\frac {3v-2m-1} {2}}+1\ge 3v-\frac
{3m} {2}+1>2v-m,$$
the last inequality being again due to the assumption $v\ge m$.
Hence, (4.27) is satisfied.

\smallskip
Consequently, if we are not in the case $v=m=1$, then the rows
$\cl{(v-1)/2},\dots,\mathbreak 
\cl{v-m/3}-2+\chi(m\equiv0\pmod3),\fl{(3v-m-1)/2}+
1,\dots,2v-m$ are rows with zeros in columns
$\cl{(v-m-1)/2},\dots,\cl{v-(2m)/3}-1,\fl{(3v-2m-1)/2}+1,\dots,2v-m$.
These are
$$\cl{v-\frac {m} {3}}-1+\chi(m\equiv0\kern-5pt\pmod3)-\cl{\frac {v-1} {2}}
+2v-m-\fl{\frac {3v-m-1} {2}}
\tag4.38$$
rows, containing possibly nontrivial entries in only
$$\cl{\frac {v-m-1} {2}}+\fl{\frac {3v-2m-1} {2}}-\cl{v-\frac {2m}
{3}}+1\tag4.39$$
columns. By simple algebra, the difference between (4.38) and (4.39)
equals
$$m+\cl{-\frac {m} {3}}+\cl{-\frac {2m}
{3}}+\chi(m\equiv0\kern-5pt\pmod3).\tag4.40$$
As is easily verified, the expression (4.40) equals 1 always. So we
have found $N+1$ rows (with $N$ the expression in (4.39)) that
actually live in $\Bbb R^N$ ($\Bbb R$ denoting the set of real numbers). 
Hence, they must be linearly dependent.
This implies that the determinant of $\Cal N$ must be 0.

Finally we settle the case $v=m=1$. The matrix $\Cal N$ then is a
$2\times 2$ matrix (cf\. Figure~2) in which column 1 vanishes. For,
$i=0$ and $j=1$ satisfy (4.25), (4.26), and (4.28), while $i=1$ and
$j=1$ satisfy (4.25), (4.26), and (4.27). Hence, $\det(\Cal N)=0$.

\smallskip
Altogether, this establishes that $P_4(x,x+m;n)$ vanishes for $x=-v$,
$m\le v\le (n+m-1)/2$, so that $(x+m)_{\fl{(n-m+1)/2}}$ is a factor
of $P_4(x,x+m;n)$. Since
$$\fl{\frac {n-m+1} {2}}=\cases \fl{\frac {n} {2}}-\fl{\frac {m}
{2}}+1&\text {$n$ odd, $m$ even}\\
\fl{\frac {n} {2}}-\fl{\frac {m} {2}}&\text {otherwise,}\endcases$$
it follows that $(x+m)_{\fl{n/2}-\fl{m/2}}$ is a factor of
$P_4(x,x+m;n)$, as desired.

\smallskip
{\it Step 2. The degree of $P_4(x,x+m;n)$ is at most
$\fl{n/2}$.}
Fortunately, this was, implicitly, already proved in Step~5 of the
proof of Theorem~2. To see this, we examine parts of Step~5 in more
detail and relate them to $P_4(x,x+m;n)$. 

First we go back to (3.37). In what follows after (3.37), it is shown
that the degree in $x$ of the polynomial (3.37) is at most $2\binom
n2+\binom{n-1}2+ \fl{(n-1)/2}+\fl{m/2}$. Now, a closer look at the
manipulations before (3.37) unfolds that (3.37) equals the first term
on the right-hand side of (3.27). However, this first term on the
right-hand side of (3.27) is nothing else but the $(0,0)$-entry of
the determinant in (3.26) times the minor of the same determinant
consisting of rows $1,2,\dots,n-1$, and columns $1,2,\dots,n-1$. This
minor is
$$\multline \det_{1\le i,j\le
n-1}\big((2x+m+z(i)+1)_{j-1}\,(x+2z(i)-j+2)_{j-1}\\
\cdot
(x+m+2j-z(i)+2)_{2n-2j-2}\,(m+3j-3z(i))\big)\\
=\det_{0\le i,j\le
n-2}\big((2x+m+z(i+1)+1)_{j}\,(x+2z(i+1)-j+1)_{j}\hskip2cm\\
\cdot
(x+m+2j-z(i+1)+4)_{2n-2j-4}\,(m+3j+3-3z(i+1))\big).
\endmultline\tag4.41$$
The $(0,0)$-entry of the determinant in (3.27) is $\SA(x,x+m;n)$, which
by its definition (3.7) is a polynomial in $x$ of degree
$2n-2+\fl{m/2}$ (as is also remarked after (3.37)). So, the degree in
$x$ of the determinant in (4.41) is at most
$3\binom{n-1}2+\fl{(n-1)/2}$. Now it should be noted that upon
replacing $n$ by $n+1$, $x$ by $x-1$, and upon setting $z(i+1)=i+1$,
$i=0,1,\dots,n-1$, the right-hand side of (4.41) turns into the
determinant in (4.3) with $y=x+m$, which by definition is
$\DetBa(x,x+m;n)$. Hence, the degree in $x$ of $\DetBa(x,x+m;n)$ is at most
$3\binom n2+\fl{n/2}$, and so, by (4.5) and (4.6), the degree in $x$
of $P_4(x,x+m;n)$ is at most $\fl{n/2}$.

\smallskip
This finishes the proof of Proposition~5.\quad \quad \qed
\enddemo
In the next step, we use Proposition~5 to evaluate the determinant
$\DetB(x,x+m;n)$ for $m=0$ and $m=1$. The case $m=0$ is the one that
we need for the evaluation of the determinant in (2.2b), the case
$m=1$ is needed for the evaluation of $\DetB(x,y;n)$, for independent
$x$ and $y$, in the proof of Theorem~8.
\proclaim{Corollary 6}Let $x$ and $n$ be nonnegative integers. Then
the determinant
$$\DetB(x,x+m;n)=\det_{0\le i,j\le n-1}\(
\frac
{(2x+m+i+j-1)!\,(m+3j-3i)}
{(x+2i-j+1)!\,(x+m+2j-i+1)!}\)$$
for $m=0$ equals
$$\multline \cases \dsize\frac {n!} {(n/2)!}\prod _{i=0} ^{n-1}\(\frac
{i!\,(2x+i-1)!\,(3x+2i+1)_i^2}
{(x+2i+1)!^2}\)\cdot
(x)_{n/2}&n\text { even}\\
0&n\text { odd}\endcases\\
=\cases \dsize\prod _{i=0} ^{n-1}\frac
{(3x+2i+1)_i^2}
{(x+2i+1)!^2}\prod _{i=0} ^{n/2-1}\big((2x+2i)!^2\,(2i+1)!^2\big)&n\text {
even}\\0&n\text { odd,}\endcases
\endmultline\tag4.42$$
and for $m=1$, $n\ge1$, equals
$$\frac {n!} {\fl{n/2}!}\prod _{i=0} ^{n-1}\(\frac
{i!\,(2x+i)!\,(3x+2i+2)_i\,(3x+2i+3)_i}
{(x+2i+1)!\,(x+2i+2)!}\)\cdot
(x+1)_{\fl{n/2}}.\tag4.43$$
\endproclaim

\demo{Proof}By Proposition~5 we know exactly how $\DetB(x,x+m;n)$
factors, except for the polynomial $P_5(x;m,n)$. However, also by
Proposition~5, for $m=0,1$ the degree of $P_5(x;m,n)$ is at most $0$.
Hence, $P_5(x;m,n)$ is a constant for $m=0$ and for $m=1$.

A combination of (4.6), with $y=x+m$, and (4.11) yields
$$\multline P_5(x;m,n)=\frac {1} {(x+m)_{\fl{n/2}-\fl{m/2}}}
\prod _{i=0}
^{n-1}\frac {1} {(3x+m+2i+1)_{i}\,(3x+2m+2i+1)_{i}}\\
\times \det_{0\le i,j\le
n-1}\big((2x+m+i)_{j}\,(x+2i-j+2)_j\,(x+m+2j-i+2)_{i}\,
(m+3j-3i)\big).
\endmultline\tag4.44$$
Thus, the (constant) value of $P_5(x;m,n)$ can be determined, 
by finding an appropriate special value for $x$, which allows to
evaluate the determinant in (4.44). 

We choose $x=-\fl{(n+m)/2}+1/2$, $m=0,1$. 
With this choice for $x$, the denominator
on the right-hand side of (4.44) does not vanish. So, everything is
well-defined for this specialization. In addition, 
since $(2x+m+i)_j$, which is a term in
each entry of the determinant, vanishes
if $i\le -2x-m=2\fl{(n+m)/2}-m-1$ and 
$i+j\ge -2x-m+1=2\fl{(n+m)/2}-m$, for $m=0$ and $m=1$ the determinant
takes on the form (3.68), where the submatrix $\Cal M$ is empty or a
$1\times 1$ matrix. Therefore the determinant can be easily evaluated
for this specialization. This gives $P_5(x;m,n)$ for $m=0$ and $m=1$.
Substitution of these values for $P_5(x;m,n)$ into (4.10)
yields the expressions (4.42) and (4.43).
\quad \quad \qed
\enddemo
The evaluation for the special case $m=x=0$ is implicitly in
\cite{\AnBuAA}. (It is equivalent to the determinant evaluation for
$\det(v(n))$ in Section~4 of \cite{\AnBuAA}.)

At this point we remark that (4.42) combined with Theorem~1, item (3),
(2.2b),
settles the ``$n$ odd" case of the Conjecture in the Introduction,
see Theorem~11.

\medskip
To be able to evaluate the determinant $\DetB(x,y;n)$ of
(4.1) completely, for independent $x$ and $y$, we need one more
auxiliary result. It locates several zeros of the polynomial factor
$P_4(x,y;n)$ of $\DetB(x,y;n)$ (recall (4.2)).
\proclaim{Lemma 7}If $u,v$ are nonnegative integers with $u+v\le
n-1$, then $P_4(-u,-v;n)=0$, with $P_4(x,y;n)$ the polynomial in
{\rm(4.2)}.
\endproclaim
\demo{Proof}Let $u,v$ be nonnegative integers with $u+v\le n-1$. 
The polynomial
$P_4(x,y;n)$ is defined by (4.6), 
$$\multline
\DetBb(x,y;n)=\prod _{i=0}
^{n-1}\big((2x+y+2i+1)_{i}\,(x+2y+2i+1)_{i}\big)\cdot P_4(x,y;n),
\endmultline\tag4.45$$
where $\DetBb(x,y;n)$ is the determinant in (4.4). What 
we would like to do is to set $x=-u$ and $y=-v$ in (4.45), prove that
$\DetBb(-u,-v;n)$ equals 0, that the product on the right-hand side of
(4.45) is nonzero, and conclude that therefore $P_4(-u,-v;n)$ must
be 0. However, the product on the right-hand side of (4.45)
unfortunately (usually) {\it is\/} 0 for $x=-u$ and $y=-v$. 
So we are in exactly the same situation as in Step~1 of Proposition~5.
The specialization of $P_4(x,y;n)$ that we are considering here is
very different, though. Curiously enough, the arguments of Step~1 of
Proposition~5 can still be used here, word by word,
with suitable replacements of
parameters. To get convinced that this is indeed the case, it will
suffice to do the very beginning. Soon it will become clear that
everything runs in parallel with Step~1 of Proposition~5.

To begin with, we set $y=-v$ in (4.45). Before setting $x=-u$, we
have to cancel all factors of the form $x+u$ that occur in the
product on the right-hand side of (4.45).
To accomplish this, we have to ``generate" these factors on the
left-hand side. Here, this is done by reading through Step~2 of the proof
of Proposition~4 with $y=-v$. To make this more precise, observe that
$x+u$ divides a typical factor $2x-v+2i+l+1$, $0\le i\le n-1$, $0\le
l<i$, of the first half of the product in (4.45) if and only if
$2u=-v+2i+l+1$. Therefore, if we recall (4.7), for each solution $(i,l)$ of
$$2u=-v+2i+l+1,\quad  \text {with }0\le i\le n-1,\ 0\le l<i,\tag4.46$$
we subtract the linear combination
$$\multline \sum _{s=l} ^{\fl{(i+l)/2}}\frac
{(2i-3s+l)} {(i-s)}\frac {(i-2s+l+1)_{s-l}} {(s-l)!}\frac
{(x+2s+2)_{2i-2s}} {(-x-2i-l+s-1)_{i-s}}\\
\cdot(\text {row $s$ of
$\DetBb(x,-v;n)$})
\endmultline\tag4.47$$
of rows of $\DetBb(x,-v;n)$ from row $i$ of $\DetBb(x,-v;n)$. 
By (4.7), the
effect is that $(2x-v+2i+l+1)=2(x+u)$ (the equality being due to
(4.46)), is a factor of each
entry of the $i$-th row of the new determinant, for each solution $(i,l)$
of (4.46).

Now it should be observed that (4.46) is exactly equivalent to
(4.14) with the replacements $v\to u$ and $m\to u-v$, while (4.47) is
exactly (4.15) with the replacements $m\to y-x$ and $y\to -v$, and
the determinant $\DetB(x,-v;n)$ that we are considering here is
exactly the determinant $\DetB(x,x+m;n)$ that is considered in Step~1
of Proposition~5, with the same replacements. 
This observation makes it apparent that similar replacements in the
rest of Step~1 of Proposition~5 will produce a valid proof of Lemma~7. In
particular, in (4.17), in the statements of Claim~1 and Claim~2, in
(4.20), (4.21), Figure~2, (4.25)--(4.40), the replacements $v\to u$
and $m\to u-v$ yield what we need here. We leave the details to the
reader.\quad \quad \qed
\enddemo
Now we are in the position to prove the promised full evaluation of
the determinant $\DetB(x,y;n)$.
\proclaim{Theorem 8}Let $x,y,n$ be nonnegative integers.
Then
$$\align \DetB(x,y;n)&=\det_{0\le i,j\le n-1}\(
\frac
{(x+y+i+j-1)!\,(y-x+3j-3i)}
{(x+2i-j+1)!\,(y+2j-i+1)!}\)\\
&=\prod _{i=0} ^{n-1}\(\frac {i!\,(x+y+i-1)!\,(2x+y+2i+1)_i\,(x+2y+2i+1)_i}
{(x+2i+1)!\,(y+2i+1)!}\)\\
&\hskip2cm\cdot \sum _{k=0} ^{n}(-1)^k\binom nk (x)_k\,(y)_{n-k}.
\tag4.48\endalign$$
\endproclaim
\demo{Proof}Obviously, the Theorem is equivalent to the assertion
that with
$$P_6(x,y;n)=\bigg(\prod _{i=0} ^{n-1}i!\bigg)
\cdot \sum _{k=0} ^{n}(-1)^k\binom nk (x)_k\,(y)_{n-k}
\tag4.49$$ 
there holds $P_4(x,y;n)=P_6(x,y;n)$, where $P_4(x,y;n)$ is the polynomial
in (4.2).

For the proof of this assertion we check the following properties for
$P_6(x,y;n)$:
\roster
\item $P_6(x,y;n)$ is a polynomial in $x$ and $y$ of (total) degree
$n$.
\item $P_6(-u,-v;n)=0$ for all nonnegative integers $u$ and $v$ with
$u+v\le n-1$.
\item $P_6(y,x;n)=(-1)^nP_6(x,y;n)$.
\item $\dsize P_6(x,x+1;n)=\frac {n!} {\fl{n/2}!}
\bigg(\prod _{i=0} ^{n-1}i!\bigg)\cdot (x+1)_{\fl{n/2}}$.
\endroster
It should be noted that all these properties are also satisfied by
$P_4(x,y;n)$. This is because of Proposition~4 for (1), because of
Lemma~7 for (2), because of 
$$\DetB(x,y;n)=(-1)^{n}\DetB(y,x;n)\tag4.50$$
(if combined with (4.2)) for (3) (identity (4.50) results from
transposing the matrix in (4.1)), 
and because of (4.43) (if combined with (4.2)) for (4). Since we also
show that
\roster
\item "(5)"The conditions (1)--(4) determine a polynomial in $x$ and
$y$ uniquely,
\endroster
the assertion follows.

\smallskip
{\it ad\/} (1). This is obvious from the definition (4.49).

\smallskip
{\it ad\/} (2). We have $(-u)_k=0$ for $k>u$. Hence, if $k>u$ the
corresponding summand in the sum in (4.49) vanishes for $x=-u$ and
$y=-v$. Now let $k\le u$. Because of $u+v\le n-1$ it follows that
$k<n-v$, or equivalently, $n-k>v$. But this implies $(-v)_{n-k}=0$.
Therefore also any summand with $k\le u$ vanishes for $x=-u$ and
$y=-v$. Thus, $P_6(-u,-v;n)=0$, as desired.

\smallskip
{\it ad\/} (3). This is obvious from the definition (4.49).

\smallskip
{\it ad\/} (4). Setting $y=x+1$ in (4.49), we get
$$\frac {P_6(x,x+1;n)} {\prod _{i=0} ^{n-1}i!}=\sum _{k=0} ^{n}(-1)^k\binom
nk (x)_k\, (x+1)_{n-k},$$
or in hypergeometric notation (cf\. the Appendix for the definition
of the $F$-notation),
$$\frac {P_6(x,x+1;n)} {\prod _{i=0} ^{n-1}i!}=(x+1)_n\cdot
{}_2F_1\!\[\matrix -n,x\\-n-x\endmatrix; -1\].$$
Next we use the contiguous relation
$${} _{2} F _{1} \!\left [ \matrix { a, b}\\ { c}\endmatrix ; {\displaystyle
   z}\right ]  = {} _{2} F _{1} \!\left [ \matrix { a , b+1}\\ {
    c}\endmatrix ; {\displaystyle z}\right ]  - 
   {{az }\over
    {c}}
   {} _{2} F _{1} \!\left [ \matrix { a + 1, b+1}\\ { c+1}\endmatrix ;
        {\displaystyle z}\right ]
\tag4.51$$
to obtain
$$
\frac {P_6(x,x+1;n)} {\prod _{i=0} ^{n-1}i!}
=(x+1)_n\,\bigg(
{}_2F_1\!\[\matrix -n,1+x\\-n-x\endmatrix; -1\]+
\frac {n} {n+x}
{}_2F_1\!\[\matrix 1-n,1+x\\1-n-x\endmatrix; -1\]\bigg).
$$
To 
the $_2F_1$-series we apply the quadratic transformation (see \cite{\RaVeAA,
(3.2)})
$$
{} _{2} F _{1} \!\left [ \matrix { a, b}\\ { 1 + a - b}\endmatrix ;
   {\displaystyle z}\right ]  
=  {{{{\left( 1 +z \right) }^{-a}}}}
  {{  {} _{2} F _{1} \!\left [ \matrix { {a\over 2}, {1\over 2} + {a\over
       2}}\\ { 1 + a - b}\endmatrix ; {\displaystyle {{4 z}\over {{{\left( 1
       + z \right) }^2}}}}\right ] }}.
\tag4.52$$
This gives
$$\multline \frac {P_6(x,x+1;n)} {\prod _{i=0} ^{n-1}i!}=(x+1)_n
\lim _{z\to-1}\,\bigg((1+z)^n{}_2F_1\!\[\matrix -\frac {n} {2},\frac {1} {2}-\frac {n} {2}\\
-n-x\endmatrix; \frac {4z} {(1+z)^2}\]\\
+\frac {n} {n+x}(1+z)^{n-1}
{}_2F_1\!\[\matrix \frac {1} {2}-\frac {n} {2},1-\frac {n} {2}\\1-n-x
\endmatrix; \frac {4z} {(1+z)^2}\]\bigg).
\endmultline$$
Now, when performing the limit, only one term survives 
on the right-hand side, either in the first $_2F_1$-series or in the
second, depending on whether $n$ is odd or even. After
simplification, it is seen that both cases result in
$$\frac {P_6(x,x+1;n)} {\prod _{i=0} ^{n-1}i!}=
\frac {n!} {\fl{n/2}!}(x+1)_{\fl{n/2}},
$$
which is what we want.

\smallskip
{\it ad\/} (5). Let $Q(x,y)$ be a polynomial in $x$ and $y$
satisfying conditions (1)--(4). Because of (1), $Q(x,y)$ can be
written in the form
$$Q(x,y)=\underset i+j\le n\to{\sum _{i,j\ge0}
^{}}a_{ij}\,(x)_i\,(y)_j,\tag4.53$$
with uniquely determined coefficients $a_{ij}$. Now, in (4.53) 
we set $x=0$ and $y=-v$, $0\le v\le n-1$. Because of (2), we obtain
$0=\sum _{j=0} ^{v}a_{0j}\,(-v)_j.$
{}From this system of equations we get $a_{0j}=0$ for $0\le j\le
n-1$. Similarly, by using (2) with $x=-1,-2,\dots,-(n-1)$, we get
$a_{ij}=0$ whenever $i+j\le n-1$. 

Thus, $Q(x,y)$ can be written in the form
$$Q(x,y)=\sum _{k=0} ^{n}b_k\,(x)_k\,(y)_{n-k},\tag4.54$$
where we set $b_k:=a_{k,n-k}$.

Now we apply (3). Since the coefficients $b_k$ in the expansion (4.54)
are uniquely determined, we get $b_{k}=(-1)^nb_{n-k}$, and so
$$\align Q(x,y)&=\sum _{k=0} ^{\fl{n/2}}b_k\,(x)_k\,(y)_{n-k}+
\sum _{k=0} ^{\cl{n/2}-1}b_{n-k}\,(x)_{n-k}\,(y)_{k}\\
&=\sum _{k=0} ^{\fl{n/2}}b_k\,(x)_k\,(y)_{n-k}+
(-1)^n\sum _{k=0} ^{\cl{n/2}-1}b_{k}\,(x)_{n-k}\,(y)_{k}.
\endalign$$
Finally we set $y=x+1$ in this equation and use condition (4). This
leads to
$$\multline \frac {n!} {\fl{n/2}!}
\bigg(\prod _{i=0} ^{n-1}i!\bigg)(x+1)_{\fl{n/2}}=(x+1)_{\fl{n/2}}\sum _{k=0}
^{\fl{n/2}}b_k\,(x)_k\,(x+\fl{n/2}+1)_{n-k-\fl{n/2}}\\
+(-1)^n(x+1)_{\fl{n/2}}\sum _{k=0}
^{\cl{n/2}-1}b_k\,x\cdot(x+\fl{n/2}+1)_{n-k-\fl{n/2}-1}\,(x+1)_k,
\endmultline$$
and after cancellation,
$$\multline \frac {n!} {\fl{n/2}!}
\bigg(\prod _{i=0} ^{n-1}i!\bigg)=\sum _{k=0}
^{\fl{n/2}}b_k\,(x)_k\,(x+\fl{n/2}+1)_{\cl{n/2}-k}
\\+
(-1)^n\sum _{k=1}
^{\cl{n/2}}b_{k-1}\,(x)_k\cdot(x+\fl{n/2}+1)_{\cl{n/2}-k}.
\endmultline\tag4.55$$

We distinguish between $n$ being even or odd. First let $n$ be even.
It is straight-forward to see that the polynomials
$$(x)_k\,(x+n/2+1)_{n/2-k},\quad  k=0,1,\dots,n/2,$$ 
are linearly
independent. Hence, by comparison of coefficients, equation (4.55) is
equivalent to a system of equations of the form
$$b_0=c_0,\ b_1+b_0=c_1,\ b_2+b_1=c_2,\ \dots,\
b_{n/2}+b_{n/2-1}=c_{n/2},$$
where $c_0,c_1,\dots,c_{n/2}$ are certain uniquely determined
constants. This system of equations has a unique solution, which
implies that $Q(x,y)$ is uniquely determined.

The case of odd $n$ is handled similarly. Here, the polynomials
$$(x)_k\,(x+(n+1)/2)_{(n+1)/2-k},\quad  k=0,1,\dots,(n+1)/2,$$ 
are linearly independent. Hence, equation (4.55) is
equivalent to a system of equations of the form
$$\multline
b_0=c'_0,\ b_1+b_0=c'_1,\ b_2+b_1=c'_2,\ \dots,\\
b_{(n-1)/2}+b_{(n-3)/2}=c'_{(n-1)/2},\ b_{(n-1)/2}=c'_{(n+1)/2}
\endmultline$$
where $c'_0,c'_1,\dots,c'_{(n+1)/2}$ are certain uniquely determined
constants. Again, this system of equations has a unique solution, which
implies that $Q(x,y)$ is uniquely determined also in this case.

This completes the proof of the Theorem.\quad \quad \qed

\enddemo
Now we can also say explicitly what the polynomial factor $P_5(x;m,n)$
of $\DetB(x,x+m;n)$ in Proposition~5 is.
\proclaim{Theorem 9}Let $x,m,n$ be nonnegative integers with $m\le
n$. Then
$$\multline \DetB(x,x+m;n)=\det_{0\le i,j\le n-1}\(
\frac
{(2x+m+i+j-1)!\,(m+3j-3i)}
{(x+2i-j+1)!\,(x+m+2j-i+1)!}\)\\
=\frac {n!} {\fl{n/2}!}\prod _{i=0} ^{n-1}
\(\frac {i!\,(2x+m+i-1)!\,(3x+m+2i+1)_i\,(3x+2m+2i+1)_i}
{(x+2i+1)!\,(x+m+2i+1)!}\)\\
\times(x+m)_{\fl{n/2}-\fl{m/2}}\cdot \sum _{k\ge0} ^{}\binom
m{2k+\chi(n\text { is odd})}\\
\times(\fl{n/2}-k+1)_k\,
(x+\cl{m/2}+\fl{n/2})_{\fl{m/2}-k},
\endmultline\tag4.56$$
again with $\chi(\Cal A)$=1 if $\Cal A$ is
true and $\chi(\Cal A)$=0 otherwise.

\endproclaim
\demo{Proof}We put $y=x+m$ in Theorem~8. Comparison of (4.48) and
(4.56) then reveals that we have to show
$$\multline
\sum _{k=0} ^{n}(-1)^k\binom nk (x)_k\,(x+m)_{n-k}\\
=\frac {n!} {\fl{n/2}!}\sum _{k\ge0} ^{}\binom
m{2k+\chi(n\text { is odd})}(\fl{n/2}-k+1)_k\,(x+\cl{m/2}+\fl{n/2})_{\fl{m/2}-k}.
\endmultline\tag4.57$$
Actually, this was already done for the special case $m=1$ when we
checked condition (4) in the proof of Theorem~8. Therefore we have to 
generalize what we did there.

First we write the left-hand side of (4.57) in hypergeometric
notation (cf\. the Appendix for the definition of the $F$-notation),
$$
(x+m)_n\cdot {}_2F_1\!\[\matrix -n,x\\1-n-x-m\endmatrix; -1\].$$
Iteration of the contiguous relation (4.51) then turns this
expression into
$$
(x+m)_n\sum _{k\ge0} ^{}\binom mk \frac {(-n)_k} {(1-n-x-m)_k}
{}_2F_1\!\[\matrix -n+k,x+m\\1+k-n-x-m\endmatrix; -1\].$$
Now we can again apply the quadratic transformation (4.52) to obtain
$$\multline
\sum _{k\ge0} ^{}\binom mk \frac {(x+m)_n\,(-n)_k} {(1-n-x-m)_k}\\
\times\lim _{z\to-1}\((1+z)^{n-k}{}_2F_1\!\[\matrix -\frac {n} {2}+\frac {k}
{2},-\frac {n} {2}+\frac {k} {2}+\frac {1} {2}\\1+k-n-x-m\endmatrix; \frac
{4z} {(1+z)^2}\]\),
\endmultline$$
and when expanding the $_2F_1$-series according to its definition and
simplifying a little,
$$
\sum _{k\ge0} ^{}\binom mk \sum _{\ell\ge0} ^{}\frac {(x+m)_n\,(-n)_{k+2\ell}} 
{(1-n-x-m)_{k+\ell}}\frac {(-1)^\ell} {\ell!}
\lim _{z\to-1}(1+z)^{n-k-2\ell}.
$$
Now, the limit is nonzero only if $n$ and $k$ have the same parity
and if $\ell=(n-k)/2$, and in that case it is 1. By substituting
$2k+\chi(n\text { is odd})$ for $k$ and by little manipulation we arrive
finally at the right-hand side of (4.57). Thus, (4.57) is established, and
therefore the Theorem.\quad \quad \qed
\enddemo

In Section~3 we formulated a Conjecture about the ``extra" polynomial factor
$P_1(x;m,n)$ that occurs in the evaluation of the determinant
$\DetA(x,x+m;n)$ in Theorem~2. 
An analogous result seems to hold for the ``extra" polynomial factor $P_5(x;m,n)$
of $\DetB(x,x+m;n)$ as given in Proposition~5, which was identified
as
$$\frac {n!} {\fl{n/2}!}\bigg(\prod _{i=0} ^{n-1}i!\bigg)\cdot
\sum _{k\ge0} ^{}\binom
m{2k+\chi(n\text { is odd})}(\fl{n/2}-k+1)_k\,(x+\cl{m/2}+\fl{n/2})_{\fl{m/2}-k}
$$
by Theorem~9.
\proclaim{Conjecture}Let $x,m,n$ be nonnegative integers with $m\le
n$. Then the polynomial
$$
\sum _{k\ge0} ^{}\binom
m{2k+\chi(n\text { is odd})}(\fl{n/2}-k+1)_k\,(x+\cl{m/2}+\fl{n/2})_{\fl{m/2}-k},
\tag4.58$$
a polynomial in $x$ of exact degree $\fl{m/2}$, satisfies:
If the cases $n$ even and $n$ odd are considered
separately, its coefficient of $x^e$ is a polynomial
in $n$ of degree $\fl{m/2}-e$ with positive integer coefficients.
\endproclaim

\subhead 5. A related determinant identity\endsubhead
In this section we derive a determinant identity that is somewhat
related to the determinant identities of the previous sections (see
the paragraph after the proof of Theorem~10 for an account of this
relationship).
Special cases of this identity appeared previously in the paper
\cite{\AnBuAA} of Andrews and Burge, also in connection with the 
enumeration of totally symmetric self-complementary plane partitions.

\medskip
In \cite{\AnBuAA, sec.~4}, Andrews and Burge show that the
determinants in (2.2) for $x=0$ (which give the enumeration of
totally symmetric self-complementary plane partitions) can be
transformed by elementary row and column operations into the
determinant
$$\det_{0\le i,j\le n-1}\(\binom {i+j+1} {2j-i} +\binom {i+j}
{2j-i-1}\),\tag5.1$$
and in Theorem~2 of their paper (see also \cite{\AnStAA, Theorem~3}) 
provide an evaluation even for
$$\det_{0\le i,j\le n-1}\(\binom {x+i+j+1} {2j-i+1} +\binom {x+i+j}
{2j-i}\).\tag5.2$$
(We changed the notation of \cite{\AnBuAA} slightly.
In particular, we replaced $x$ by $x-2$.) Then they observe that the
determinant (5.1) reduces to the determinant (5.2) with $x=2$ and
with $n$ replaced by $n-1$, and thus provide another proof of the
totally symmetric self-complementary plane partitions conjecture.
However, there is even a two-parameter generalization of (5.2),
(namely the determinant in (5.4) below), that
can be evaluated. This two-parameter generalization 
is the subject of our next theorem. We
formulate it only for integral $x$ and $y$. But in fact, with a
generalized definition of factorials and binomials (cf\.
\cite{\GrKPAA, sec.~5.5, (5.96), (5.100)}; 
see also the remarks after (4.1)), Theorem~10,
together with its proof, would also hold for complex $x$ and $y$. 

\proclaim{Theorem 10}Let $x,y,n$ be nonnegative integers. Then there
holds
$$\multline \det_{0\le i,j\le n-1}\(\frac {(x+y+i+j-1)!}
{(x+2i-j)!\,(y+2j-i)!}\)\\
=\prod _{i=0} ^{n-1}\frac {i!\,(x+y+i-1)!\,(2x+y+2i)_i\,(x+2y+2i)_i}
{(x+2i)!\,(y+2i)!},
\endmultline\tag5.3$$
or equivalently,
$$\multline \det_{0\le i,j\le n-1}\(
\binom {x+y+i+j} {y+2j-i} +\binom {x+y+i+j-1}
{y+2j-i-1}\)\\
=\prod _{i=0} ^{n-1}\frac
{i!\,(x+y+i-1)!\,(2x+y+2i)_i\,(x+2y+2i)_{i+1}}
{(x+2i)!\,(y+2i)!}.
\endmultline\tag5.4$$
\endproclaim
\demo{Proof}The equivalence of (5.3)  and (5.4) is obvious from the
simple fact
$$\binom {x+y+i+j} {y+2j-i} +\binom {x+y+i+j-1}
{y+2j-i-1}=\frac {(x+y+i+j-1)!\,(x+2y+3j)}
{(x+2i-j)!\,(y+2j-i)!}.$$

We are going to prove (5.3). Our procedure is very similar to the
preceding proofs of Theorems~2 and 4, only that things are much
simpler here. Actually, in the research process it was the other way
round. This proof was found first and provided (some of) the
inspiration for the later proofs of Theorems~2 and 4. 

\smallskip
{\it Step 1. An equivalent statement of the Theorem.}
We take as many common factors out
of the $i$-th row of the determinant in (5.3), $i=0,1,\dots,n-1$, 
as possible, such that the
entries become polynomials in $x$ and $y$. Thus we obtain
$$\multline
\prodl _{i=0} ^{n-1}\frac {(x+y+i-1)!} {(x+2i)!\,(y+2n-i-2)!}\\
\times\det_{0\le i,j\le
n-1}\big((x+y+i)_{j}\,(x+2i-j+1)_j\,(y+2j-i+1)_{2n-2j-2}\big).
\endmultline$$
Comparing with (5.3), we see that (5.3) is equivalent to
$$\multline \det_{0\le i,j\le
n-1}\big((x+y+i)_{j}\,(x+2i-j+1)_j\,(y+2j-i+1)_{2n-2j-2}\big)\\
=\prod _{i=0}
^{n-1}\big(i!\,(y+2i+1)_{n-i-1}\,(2x+y+2i)_i\,(x+2y+2i)_i\big).
\endmultline\tag5.5$$
Let us denote the determinant in (5.5) by $P_8(x,y;n)$.

We are going to establish (5.5), and thus (5.3), by showing in
Steps~2 and 3 that the right-hand side of (5.5) divides $P_8(x,y;n)$
as a polynomial in $x$ and $y$,
by showing in Step~4 that the (total) degree in $x$ and $y$ of $P_8(x,y;n)$
is $3\binom n2$, so that $P_8(x,y;n)$ is a constant multiple of the
right-hand side of (5.5), and by showing that this constant equals 
1, also in Step~4.

\smallskip
{\it Step 2. $\prod _{i=0} ^{n-1}(y+2i+1)_{n-i-1}$
is a factor of $P_8(x,y;n)$.}
We multiply the $i$-th row of $P_8(x,y;n)$, which is the determinant
in (5.5),
by $(y+2n-i-1)_{i}$, $i=0,1,\dots,n-1$, and divide the $j$-th column by
$(y+2j+1)_{2n-2j-2}$, $j=0,1,\dots,n-1$.
This leads to
$$\multline P_8(x,y;n)=\prodl _{i=0} ^{n-1}\frac {1} {(y+2n-i-1)_{i}}\prodl
_{j=0} ^{n-1}(y+2j+1)_{2n-2j-2}\\
\hskip2cm\times\det_{0\le i,j\le
n-1}\big((x+y+i)_{j}\,(x+2i-j+1)_j\,(y+2j-i+1)_{i}\big)\\
=\prod _{i=0} ^{n-1}(y+2i+1)_{n-i-1}\cdot\det_{0\le i,j\le
n-1}\big((x+y+i)_{j}\,(x+2i-j+1)_j\,(y+2j-i+1)_{i}\big).
\endmultline$$
Since the determinant in the last line is a polynomial in $x$ and
$y$, we infer that $\prod _{i=0} ^{n-1}(y+2i+1)_{n-i-1}$ divides
$P_8(x,y;n)$.

\smallskip
{\it Step 3. $\prod _{i=0}
^{n-1}\big((2x+y+2i)_{i}\,(x+2y+2i)_{i}\big)$
is a factor of $P_8(x,y;n)$.}
We proceed in the spirit of Step~3 of the proof of Theorem~2. So it
is not necessary to provide all the details. The basic ideas will
suffice.

First, let us concentrate on a typical factor
$(2x+y+2i+l)$, $0\le i\le n-1$, $0\le l<i$, of the first half of the
product, $\prod _{i=0}
^{n-1}(2x+y+2i)_{i}$. We claim that 
for each
such factor there is a linear combination of the rows that vanishes
if the factor vanishes. More precisely, we claim that for any $i,l$
with $0\le i\le n-1$, $0\le l<i$ there holds
$$\multline \hskip-10pt\sum _{s=l} ^{\fl{(i+l)/2}}\frac
{(2i-3s+l)} {(i-s)}\frac {(i-2s+l+1)_{s-l}} {(s-l)!}\frac
{(x+2s+1)_{2i-2s}} {(-x-2i-l+s)_{i-s}}(-2x-2i-l+2n-s-1)_s\\
\cdot(\text {row $s$ of
$P_8(x,-2x-2i-l;n)$})
=(\text {row $i$ of $P_8(x,-2x-2i-l;n)$}).
\endmultline$$
Restricting to the $j$-th column, it is seen that this means to
check
$$\multline \hskip-15pt\sum _{s=l} ^{\fl{(i+l)/2}}\frac
{(2i-3s+l)} {(i-s)}\frac {(i-2s+l+1)_{s-l}} {(s-l)!}\frac
{(x+2s+1)_{2i-2s}} {(-x-2i-l+s)_{i-s}}(-2x-2i-l+2n-s-1)_{s-i}\\
\times(-x-2i-l+s)_{j}\,(x+2s-j+1)_{j}\,
(-2x-2i-l+2j-s+1)_{2n-2j-2}\\
=(-x-i-l)_{j}\,(x+2i-j+1)_{j}\,
(-2x-3i-l+2j+1)_{2n-2j-2}.
\endmultline$$
This is easily done by observing that it is equivalent to (3.15) with
$x$ replaced by $x-j$.
Arguments that are similar to those after (3.18) then show that 
the complete product $\prod _{i=0}
^{n-1}(2x+y+2i)_{i}$ divides $P_8(x,y;n)$.

The reasoning that $\prod _{i=0}
^{n-1}(x+2y+2i)_{i}$ is a factor of $P_8(x,y;n)$ is
similar. Also here, let us concentrate on a typical factor
$(x+2y+2j+l)$, $0\le j\le n-1$, $0\le l<j$. This time
we claim that for each
such factor there is a linear combination of the columns that vanishes
if the factor vanishes. More precisely, we claim that for any $j,l$
with $0\le j\le n-1$, $0\le l<j$ there holds
$$\multline \sum _{s=l} ^{\fl{(j+l)/2}}\frac
{(2j-3s+l)} {(j-s)}\frac {(j-2s+l+1)_{s-l}} {(s-l)!}
\cdot(\text {column $s$ of
$P_8(-2y-2j-l,y;n)$})\\
=(\text {column $j$ of $P_8(-2y-2j-l,y;n)$}).
\endmultline$$
Restricting to the $i$-th row, we see that this means to check
$$\multline \sum _{s=l} ^{\fl{(j+l)/2}}\frac
{(2j-3s+l)} {(j-s)}\frac {(j-2s+l+1)_{s-l}} {(s-l)!}\\
\times(-y-2j-l+i)_{s}\,(-2y-2j-l+2i-s+1)_{s}\,
(y+2s-i+1)_{2n-2s-2}\\
=(-y-2j-l+i)_{j}\,(-2y-3j-l+2i+1)_{j}\,(y+2j-i+1)_{2n-2j-2}.
\endmultline\tag5.6$$
The observation that this summation is equivalent to (3.23) with
$y$ replaced by $x+2$ and $\y$ replaced by $y-i$ 
establishes the claim. Similarly to as before, this eventually
shows that the complete product $\prod _{i=0}
^{n-1}(x+2y+2i)_{i}$ divides $P_8(x,y;n)$. 

Altogether, this implies that $\prod _{i=0}
^{n-1}\big((2x+y+2i)_{i}
(x+2y+2i)_{i}\big)$ divides $P_8(x,y;n)$,
as desired.

\smallskip
{\it Step 4. $P_8(x,y;n)$ is a polynomial in $x$ and $y$ of degree
$3\binom n2$, and the evaluation of the multiplicative constant.}
We consider the generalized determinant
$$\bP_8(x,y,z(0),z(1),\dots,z(n-1);n)=\bP_8(x,y;n)$$
which arises from $P_8(x,y;n)$ by replacing each occurence of $i$
in row $i$ by an indeterminate, $z(i)$ say, $i=0,1,\dots,n-1$,
$$\multline 
\bP_8(x,y,z(0),z(1),\dots,z(n-1);n)=\bP_8(n)\\
=\det_{0\le i,j\le
n-1}\big((x+y+z(i))_{j}\,(x+2z(i)-j+1)_j\,(y+2j-z(i)+1)_{2n-2j-2}\big).
\endmultline$$
This determinant is a polynomial in $x, y,z(0),z(1),\dots,z(n-1)$ of
(total) degree at most $4\binom n2$. 

When $z(i_1)=z(i_2)$ for some $i_1\ne i_2$, the rows $i_1$ and
$i_2$ in $\bP_8(n)$ are identical. Hence $\bP_8(n)$ vanishes in
this case. This shows that
the product $\prod _{0\le i<j\le n-1} ^{}(z(j)-z(i))$ divides
$\bP_8(n)$. Moreover, the argument in the second half of Step~3 
shows that also $\prod _{i=0} ^{n-1}(x+2y+2i)_i$ divides
$\bP_8(n)$, just replace $i$ by $z(i)$ in (5.6). Thus we obtain
that
$$\bP_8(n)=\prod _{0\le i<j\le n-1} ^{}(z(j)-z(i))
\prod _{i=0} ^{n-1}(x+2y+2i)_i\cdot
Q(x,y,z(0),z(1),\dots,z(n-1);n),\tag5.7$$
where $Q(x,y,z(0),z(1),\dots,z(n-1);n)$ is a polynomial in
$x,y,z(0),z(1),\dots,z(n-1)$ of total degree at most $4\binom n2-2\binom
n2=2\binom n2$. By comparing coefficients of
$$y^{3\binom n2}z(0)^0z(1)^1\cdots z(n-1)^{n-1}$$
on both sides of (5.7), it is seen that the coefficient of
$y^{2\binom n2}$ in $Q(x,y,z(0),z(1),\dots,\mathbreak z(n-1);n)$ equals 1. 

Now we set $z(i)=i$, $i=0,1,\dots,n-1$, in (5.7). Then $\bP_8(n)$ on
the left-hand side reduces to $P_8(x,y;n)$. By Steps~2 and 3 
we know that 
$$\prod _{i=0} ^{n-1}\big((y+2i+1)_{n-i-1}\,(2x+y+2i)_i\big)\tag5.8$$ 
divides $P_8(x,y;n)$. Therefore, by (5.7),
it also divides $Q(x,y,0,1,\dots,n-1;n)$. Since the degree in $x$ and
$y$ of this
factor is $2\binom n2$, which at the same time is an upper bound for
the degree in $x$ and $y$ of
$Q(x,y,0,1,\dots,n-1;n)$, as we saw before, $Q(x,y,0,1,\dots,n-1;n)$
is a constant multiple of (5.8). Moreover, the coefficient of
$y^{2\binom n2}$ in (5.8) equals 1, which we already know to be the
coefficient of $y^{2\binom n2}$ in $Q(x,y,0,1,\dots,n-1;n)$. Hence, 
$Q(x,y,0,1,\dots,n-1;n)$ agrees with (5.8), which by (5.7) completes
the proof of (5.5), and hence of the Theorem.\quad \quad \qed
\enddemo
It should be noted that the two-parameter determinant in (5.4) carries a
strong relationship to the determinants $\DetA(x,y;n)$ in (3.1) and 
$\DetB(x,y;n)$ in (4.1). Namely, the $(i,j)$-entry of the determinant
in (5.4), with $x$ replaced by $x+1$, equals the $(i,j)$-entry of 
$\DetA(x,y;n)$ minus 2 times the $(i,j-1)$-entry of $\DetA(x,y;n)$,
while $(2x+y+3i+1)(x+2y+3j+1)$ times 
the $(i,j)$-entry of $\DetB(x,y;n)$ equals the $(i,j)$-entry of
the determinant in (5.4) with $x$ replaced by $x+1$ 
minus 2 times the $(i-1,j)$-entry of the
same determinant. So the determinant in (5.4) is somehow ``in between"
the determinants $\DetA(x,y;n)$ and $\DetB(x,y;n)$.

\subhead 6. Constant term identities\endsubhead
In this section we translate some of our determinant identities 
into constant term identities.

Of course, we start by stating the Conjecture of the Introduction,
now as a theorem.
\proclaim{Theorem 11}Let $x$ and $n$ be nonnegative integers. Then
there holds
$$\multline \CT \(\frac {\prod _{0\le i<j\le n-1} ^{}(1-z_i/z_j)\prod
_{i=0}
^{n-1}(1+z_i^{-1})^{x+n-i-1}}
{\prod _{0\le i<j\le n-1} ^{}(1-z_iz_j)\prod _{i=0} ^{n-1}(1-z_i)}\)\\
=\cases \prodl _{i=0} ^{n-1}\frac {(3x+3i+1)!} {(3x+2i+1)!\,(x+2i)!}
\prodl _{i=0} ^{(n-2)/2}(2x+2i+1)!\,(2i)!&\text {if $n$ is even}\\
2^x\prodl _{i=1} ^{n-1}\frac {(3x+3i+1)!} {(3x+2i+1)!\,(x+2i)!}
\prodl _{i=1} ^{(n-1)/2}(2x+2i)!\,(2i-1)!&\text {if $n$ is odd.}
\endcases
\endmultline\tag6.1$$
Also, both the sum of all $n\times n$ minors of the $n\times
(2n-1)$ matrix\linebreak 
$\(\binom {x+i}{j-i}\)_{0\le i\le n-1,\ 0\le j\le
2n+x-2}$, and the number of shifted plane partitions of shape
$(x+n-1,x+n-2,\dots,1)$, with entries between 0 and $n$, 
where the entries in row $i$ are at least
$n-i$, $i=1,2,\dots,n-1$, equal the right-hand side in
(6.1).
\endproclaim

\demo{Proof}For even $n$, equation (6.1) follows from a combination
of Theorem~1, item (3), (2.2a), and (3.69). 
For odd $n$, equation (6.1) follows from a combination
of Theorem~1, item (3), (2.2b), and (4.42) or (4.56) with $m=0$. 
The other claims are due
to Theorem~1, items (1) and (2), respectively.\quad \quad \qed
\enddemo
Next we translate the determinant identities of Section~4 into
constant term identities.
\proclaim{Theorem 12}Let $x,y,m,n$ be nonnegative integers with $m\le
n$. Then there holds
$$\multline
\CT\Big(\prod _{i=0}
^{n-1}\big((1+z_i)^{x+y+i-1}(1+2z_i)(z_i+2)(z_i-1)
z_i^{-y-2n+i+1}\big)\\
\times\prod _{0\le i<j\le n-1}
^{}\big((z_i-z_j)(z_i+z_j+z_iz_j)\big)\Big)\\
=\prod _{i=0} ^{n-1}\(\frac
{i!\,(x+y+i-1)!\,(2x+y+2i+1)_{i+1}\,(x+2y+2i+1)_{i+1}}
{(x+2i+1)!\,(y+2i+1)!}\)\\
\cdot \sum _{k=0} ^{n}(-1)^k\binom nk (x)_k\,(y)_{n-k}.
\endmultline\tag6.2$$
If $y=x+m$, with $m$ a fixed nonnegative integer, then the
constant term in {\rm(6.2)} equals
$$\multline 
\frac {n!} {\fl{n/2}!}
\prod _{i=0} ^{n-1}\(\frac
{i!\,(2x+m+i-1)!\,(3x+m+2i+1)_{i+1}\,(3x+2m+2i+1)_{i+1}}
{(x+2i+1)!\,(x+m+2i+1)!}\)\\
\times(x+m)_{\fl{n/2}-\fl{m/2}}\cdot \sum _{k\ge0} ^{}\binom
m{2k+\chi(n\text { is odd})}\\
\times(\fl{n/2}-k+1)_k\,(x+\cl{m/2}+\fl{n/2})_{\fl{m/2}-k}.
\endmultline\tag6.3$$
\endproclaim

\demo{Proof}It is routine to verify that
$$\multline \frac {(x+y+i+j-1)!\,(2x+y+3i+1)\,(x+2y+3j+1)\,(y-x+3j-3i)}
{(x+2i-j)!\,(y+2j-i)!}\\
=\CT\big((1+z)^{x+y+i+j-1}(1+2z)(z+2)(z-1)z^{-y-2j+i-1}\big).
\endmultline$$ 
Consequently, taking determinants we obtain
$$\align \det_{0\le i,j\le n-1}&\!\!\(
\frac {(x+y+i+j-1)!\,(2x+y+3i+1)\,(x+2y+3j+1)\,(y-x+3j-3i)}
{(x+2i-j)!\,(y+2j-i)!}\)\\\tag6.4\\
&=\det_{0\le i,j\le n-1}\(
\CT\big((1+z_i)^{x+y+i+j-1}(1+2z_i)(z_i+2)(z_i-1)z_i^{-y-2j+i-1}\big)\)\\
&=\CT\Big(\prod _{i=0} ^{n-1}\big(
(1+z_i)^{x+y+i-1}(1+2z_i)(z_i+2)(z_i-1)z_i^{-y+i-1}\big)\\
&\hskip3cm\times\det_{0\le i,j\le n-1}\Big(\Big(\frac {1+z_i}
{z_i^2}\Big)^j\Big)\Big)\\
&=\CT\Big(\prod _{i=0} ^{n-1}\big(
(1+z_i)^{x+y+i-1}(1+2z_i)(z_i+2)(z_i-1)z_i^{-y-2n+i+1}\big)\\
&\hskip3cm\times\prod _{0\le i,j\le n-1} ^{}\big((z_i-z_j)(z_i+z_j+z_iz_j)\big) \Big),
\endalign$$
where we used the Vandermonde determinant identity in the last step.
Obviously, the last line agrees exactly with the left-hand side of
(6.2). Thus, by taking factors that depend only on $i$, respectively
only on $j$, out of the determinant in (6.4) and applying
Theorems~8 and 9 to the resulting determinant,
all the assertions of the Theorem follow immediately.\quad \quad \qed
\enddemo
Finally, we translate the determinant identity of Section~5 into a
constant term identity. 
\proclaim{Theorem 13}Let $x,y,n$ be nonnegative integers. Then there
holds
$$\multline \CT\Big(\prod _{i=0} ^{n-1}\big(
(1+z_i)^{x+y+i-1}(1+2z_i)z_i^{-y-2n+i+2}\big)
\prod _{0\le i,j\le n-1} ^{}\big((z_i-z_j)(z_i+z_j+z_iz_j)\big) \Big)\\
=\prod _{i=0} ^{n-1}\frac
{i!\,(x+y+i-1)!\,(2x+y+2i)_i\,(x+2y+2i)_{i+1}}
{(x+2i)!\,(y+2i)!}.
\endmultline\tag6.5$$
\endproclaim
\demo{Proof}We observe
$$\binom {x+y+i+j} {y+2j-i} +\binom {x+y+i+j-1}
{y+2j-i-1}=\CT\big((1+z)^{x+y+i+j-1}(1+2z)z^{-y-2j+i}\big),$$ 
and then proceed in the same way as in the proof of Theorem~12. The
reader will have no difficulties to fill in the details.\quad \quad
\qed
\enddemo

\head \tenpoint\bf Appendix\endhead
Here we provide auxiliary results that are needed in the proofs of
our Theorems. 

We start by recalling the theorems about nonintersecting lattice
paths that we need in the proof of Theorem~1. The main theorem
of nonintersecting lattice paths 
\cite{\GeViAB, Cor.~2; \StemAE, Theorem~1.2} is the following.
\proclaim{Proposition~A1}Let $A_0,A_1,\dots,A_{n-1}$ and
$E_0,E_1,\dots,E_{n-1}$ be lattice points with the ``compatibility"
property that, given $i<j$ and $k<l$, any lattice path from $A_i$ to $E_l$
meets any lattice path from $A_j$ to $E_k$. Then the number of all
families $(P_0,P_1,\dots,P_{n-1})$ of nonintersecting lattice paths,
where $P_i$ runs from $A_i$ to $E_i$, $i=0,1,\dots,n-1$, is given by
the determinant
$$\det_{0\le i,j\le n-1}\big(\vert \Cal P(A_i\to E_j)\vert\big),$$
where $\vert\Cal P(A\to E)\vert$ denotes the number of all lattice paths from
$A$ to $E$.
\endproclaim
The second result about nonintersecting lattice paths is
Stembridge's enumeration \cite{\StemAE,
Theorem~3.1} of nonintersecting lattice paths
when the end points of the lattice paths are allowed to vary.
\proclaim{Proposition~A2}Let $A_0,A_1,\dots,A_{n-1}$ be lattice
points, and let $I=\{\dots,E_k,\mathbreak E_{k+1},\dots\}$ 
be a {\rm totally ordered} set of lattice points, again
with the ``compatibility"
property that, given $i<j$ and $k<l$, any lattice path from $A_i$ to $E_l$
meets any lattice path from $A_j$ to $E_k$. Then the number of all
families $(P_0,P_1,\dots,P_{n-1})$ of nonintersecting lattice paths,
where $P_i$ runs from $A_i$ to some point of $I$, $i=0,1,\dots,n-1$, is given by
the Pfaffian
$$\pf_{0\le i<j\le n-1}\big(Q(i,j)\big),$$
where $Q(i,j)$ is the number of all pairs $(P_i,P_j)$ of
nonintersecting lattice paths, $P_i$ running from $A_i$ to some
point of $I$, and $P_j$ running from $A_j$ to some
point of $I$. 
\endproclaim

Next we prove some identities for hypergeometric series. 
We use the usual hypergeometric notation
$${}_r F_s\!\left[\matrix a_1,\dots,a_r\\ b_1,\dots,b_s\endmatrix; 
z\right]=\sum _{k=0} ^{\infty}\frac {\po{a_1}{k}\cdots\po{a_r}{k}}
{k!\,\po{b_1}{k}\cdots\po{b_s}{k}} z^k\ ,$$
where the shifted factorial
$(a)_k$ is given by $(a)_k:=a(a+1)\cdots(a+k-1)$,
$k\ge1$, $(a)_0:=1$, as before. 

To begin with, in Lemma~A3 
we quote a result of Andrews and Burge \cite{\AnBuAA,
Lemma~1}. This $_4F_3$-summation was derived in \cite{\AnBuAA} from a
similar $_4F_3$-summation due to Bailey. We provide an alternative
proof here, showing that, in fact, Andrews and Burge's summation
follows easily from a
transformation formula due to Singh. (The above mentioned
$_4F_3$-summation of Bailey's, as well as Lemma~2 of \cite{\AnBuAA}
do also follow from Singh's transformation formula.)
\proclaim{Lemma A3}Let $n$ be a positive integer. Then

\vskip3pt
\vbox{{}
$${}_4F_3\!\[\matrix -\frac {n} {2},\frac {1} {2}-\frac {n} {2},-A,A+B\\
1-n,\frac {B} {2},\frac {1} {2}+\frac {B} {2}\endmatrix; 1\]=\frac {(A+B)_n} {(B)_n}+
\frac {(-A)_n} {(B)_n}.\tag A.1$$
}
\endproclaim
\demo{Proof} In \cite{\SinVAA, main theorem}, Singh proves the following
transformation formula (actually a $q$-analogue thereof, see also 
\cite{\GaRaAA, (3.10.13); Appendix (III.21)}):
$$
{} _{4} F _{3} \!\left [ \matrix { a, b, c, d}\\ { {1\over 2} + a + b, {{c +
   d}\over 2}, {{1 + c + d}\over 2}}\endmatrix ; {\displaystyle 1}\right ]  
=  {} _{3} F _{2} \!\left [ \matrix { 2 a, 2 b, c}\\ { {1\over 2} + a + b, c +
   d}\endmatrix ; {\displaystyle 1}\right ] ,
$$
provided both series terminate.

Now, let first $n$ be odd. In Singh's transformation we
choose $a=-n/2+\ep$, $b=1/2-n/2$, $c=-A$, and $d=A+B$. Thus, we
obtain
$$
{} _{4} F _{3} \!\left [ \matrix {  - {n\over 2}+{\ep}, {1\over 2} - {n\over
   2}, -A, A + B}\\ { 1 + {\ep} - n, {B\over 2}, {1\over 2} + {B\over
   2}}\endmatrix ; {\displaystyle 1}\right ]  = 
  {} _{3} F _{2} \!\left [ \matrix { 2 {\ep} - n, 1 - n, -A}\\ { 1 + 
{\ep} - n, B}\endmatrix ; {\displaystyle 1}\right ] .
$$
Note that because $n$ is odd both series do indeed terminate. We may
express the sum on the right-hand side explicitly,
$$
{} _{4} F _{3} \!\left [ \matrix { - {n\over 2}+{\ep}, {1\over 2} - {n\over
   2}, -A, A + B}\\ { 1 + {\ep} - n, {B\over 2}, {1\over 2} + {B\over
   2}}\endmatrix ; {\displaystyle 1}\right ]  = 
  \sum_{k = 0}^{n-1}{{({ \textstyle -A}) _{k} \,
        ({ \textstyle 1 - n}) _{k} \,({ \textstyle 2 {\ep} - n}) _{k} }
       \over {({ \textstyle 1}) _{k} \,({ \textstyle B}) _{k} \,
        ({ \textstyle 1 + {\ep} - n}) _{k} }}
$$
and then let $\ep$ tend to zero. This gives
$$\align
{} _{4} F _{3} \!\left [ \matrix { -{{n}\over 2}, {1\over 2} - {n\over 2}, -A,
   A + B}\\ { 1 - n, {B\over 2}, {1\over 2} + {B\over 2}}\endmatrix ;
   {\displaystyle 1}\right ]  &= 
  \sum_{k = 0}^{n-1}{{({ \textstyle -A}) _{k} \,
        ({ \textstyle -n}) _{k} }\over 
      {({ \textstyle 1}) _{k} \,({ \textstyle B}) _{k} }}\\
&=    \sum_{k = 0}^{n}{{({ \textstyle -A}) _{k} \,
          ({ \textstyle -n}) _{k} }\over 
        {({ \textstyle 1}) _{k} \,({ \textstyle B}) _{k} }} 
  -{{({ \textstyle -A}) _{n} \,({ \textstyle -n}) _{n} }\over 
      {({ \textstyle 1}) _{n} \,({ \textstyle B}) _{n} }} 
\\
&={}_2F_1\!\[\matrix -A,-n\\B\endmatrix; 1\]
  -{{({ \textstyle -A}) _{n} \,({ \textstyle -n}) _{n} }\over 
      {({ \textstyle 1}) _{n} \,({ \textstyle B}) _{n} }}.
\endalign$$
The $_2F_1$-series can be evaluated by the hypergeometric form of
Vandermonde summation (see \cite{\SlatAC, (1.7.7); Appendix
(III.4)}),
$$
{} _{2} F _{1} \!\left [ \matrix { a, -n}\\ { c}\endmatrix ; {\displaystyle
   1}\right ]  = {{({ \textstyle c-a}) _{n} }\over 
    {({ \textstyle c}) _{n} }},
$$
where $n$ is a nonnegative integer. Applying this, we get
$$
{} _{4} F _{3} \!\left [ \matrix { -{{n}\over 2}, {1\over 2} - {n\over 2}, -A,
   A + B}\\ { 1 - n, {B\over 2}, {1\over 2} + {B\over 2}}\endmatrix ;
   {\displaystyle 1}\right ]  
=  {{({ \textstyle A + B}) _{n} }\over {({ \textstyle B}) _{n} }} - 
   {{({ \textstyle -A}) _{n} \,({ \textstyle -n}) _{n} }\over 
     {{n!} \,({ \textstyle B}) _{n} }},
$$
which reduces to (A.1) since $n$ is odd and thus $(-n)_n=-n!$.

The case that $n$ is even is treated similarly. One would start by
choosing $a=-n/2$, $b=1/2-n/2+\ep$, $c=-A$, and $d=A+B$ in Singh's 
transformation formula. We leave the details to the reader.
\quad \quad \qed
\enddemo
{}From the $_4F_3$-summation in Lemma~A3 we derive a summation for a $_5F_4$-series.
\proclaim{Lemma A4}Let $n$ be a positive integer. Then
$$\multline
  {} _{5} F _{4} \!\left [ \matrix { 1 - {{2 n}\over 3}, -{{n}\over 2},
   {1\over 2} - {n\over 2}, -A, A + B}\\ { -{{2 n}\over 3}, 1 - n, {1\over 2}
   + {B\over 2}, 1 + {B\over 2}}\endmatrix ; {\displaystyle 1}\right ]
\\
=\frac {1} {2}{{\left( A - B - 2 n \right) }\over
 {\left( 2 A + B \right) } }{{ ({ \textstyle -A}) _{n} }\over 
     {  ({ \textstyle 1 + B}) _{n} }} + 
\frac {1} {2}  {{\left( A + 2 B + 2 n \right) }\over
  {\left( 2 A + B \right) }} {{ ({ \textstyle A + B}) _{n} }\over 
     {  ({ \textstyle 1 + B}) _{n} }} .
\endmultline\tag A.2
$$
\endproclaim
\demo{Proof}We transform the $_5F_4$-series by the contiguous
relation
$$\multline 
{} _{5} F _{4} \!\left [ \matrix { 1 - {{2 n}\over 3}, -{{n}\over 2}, {1\over
  2} - {n\over 2}, -A, A + B}\\ { -{{2 n}\over 3}, 1 - n, {1\over 2} + {B\over
  2}, 1 + {B\over 2}}\endmatrix ; {\displaystyle 1}\right ] \\
={{3 A \left( B + n \right) }\over{{4 n\left( 2 A + B \right)  }}}
 {} _{4} F _{3} \!\left [ \matrix { -{{n}\over 2},
       {1\over 2} - {n\over 2}, 1 - A, A + B}\\ { 1 - n, {1\over 2} + {B\over
       2}, 1 + {B\over 2}}\endmatrix ; {\displaystyle 1}\right ]
\hskip3cm\\
     - 
  {{3 B \left( 2 A + B + {n\over 3} \right)  }\over 
    {4 n\left( 2 A + B \right)  }}
      {} _{4} F _{3} \!\left [ \matrix { -{{n}\over 2}, {1\over 2} - {n\over
       2}, -A, A + B}\\ { 1 - n, {B\over 2}, {1\over 2} + {B\over
       2}}\endmatrix ; {\displaystyle 1}\right ] 
\\ +
  {{3  \left( B + {{5 n}\over 3} \right)  \left( A + B \right) }\over 
    {4 n\left( 2 A + B \right)  }}
      {} _{4} F _{3} \!\left [ \matrix { -{{n}\over 2}, {1\over 2} - {n\over
       2}, -A, 1 + A + B}\\ { 1 - n, {1\over 2} + {B\over 2}, 1 + {B\over
       2}}\endmatrix ; {\displaystyle 1}\right ] .
\endmultline$$
Now each of the $_4F_3$-series can be summed by means of Lemma~1.
Some manipulation then leads to (A.2).\quad \quad \qed
\enddemo
The special case that is of particular importance in Step~3 of the
proofs of
Theorems~2 and 10, and in Step~2 of the proof of Proposition~4, is $A=2n+B$.
\proclaim{Corollary A5}Let $n$ be a positive integer. Then

\vskip3pt
\vbox{{}
$$
{} _{5} F _{4} \!\left [ \matrix { 1 - {{2 n}\over 3}, -{{n}\over 2}, {1\over
  2} - {n\over 2}, -2n - B, 2n + 2 B}\\ { -{{2 n}\over 3}, 1 - n, {1\over 2} +
  {B\over 2}, 1 + {B\over 2}}\endmatrix ; {\displaystyle 1}\right ] =
\frac {1} {2}   {{  ({ \textstyle 2n + 2 B}) _{n} }\over 
     {   ({ \textstyle 1 + B}) _{n} }} .
\tag A.3
$$
\line{\hfil\hbox{\qed}\quad \quad }
}
\endproclaim
Finally we move one step further to a $_6F_5$-summation.
\proclaim{Lemma~A6}Let $n$ be a positive integer. Then
$$\multline {} _{6} F _{5} \!\left [ \matrix { {4\over 3} + {{2 n}\over 3} + B, 1 -
  {{2 n}\over 3}, -{{n}\over 2}, {1\over 2} - {n\over 2}, -1 - 2 n - B, 2 n
  + 2 B}\\ { {1\over 3} + {{2 n}\over 3} + B, -{{2 n}\over 3}, 1 - n, 1 +
  {B\over 2}, {3\over 2} + {B\over 2}}\endmatrix ; {\displaystyle 1}\right ] 
\\
=\frac {1} {2}{{\left( 1 + 5 n + 3 B \right)  }\over
 {\left( 1 + 2 n + 3 B \right)}}
{{   ({ \textstyle  2 n + 2 B}) _{ n} }\over 
   {  ({ \textstyle 2 + B}) _{n} }}
\endmultline\tag A.4$$

\endproclaim
\demo{Proof}We use the contiguous relation
$$\multline {} _{6} F _{5} \!\left [ \matrix { {4\over 3} + {{2 n}\over 3} + B, 1 -
  {{2 n}\over 3}, -{{n}\over 2}, {1\over 2} - {n\over 2}, -1 - 2 n - B, 2 n
  + 2 B}\\ { {1\over 3} + {{2 n}\over 3} + B, -{{2 n}\over 3}, 1 - n, 1 +
  {B\over 2}, {3\over 2} + {B\over 2}}\endmatrix ; {\displaystyle 1}\right ] 
\\
={{4 \left( n + B \right)  \left( 1 + 2 n + 2 B \right)  }
\over     {  \left( 1 + 4 n + 3 B \right) \left( 1 + 2 n + 3 B
\right)}}\hskip7cm\\
\hskip4cm\times
      {} _{5} F _{4} \!\left [ \matrix { 1 - {{2 n}\over 3}, -{{n}\over 2},
       {1\over 2} - {n\over 2}, -1 - 2 n - B, 2 + 2 n + 2 B}\\ {
       -{{2 n}\over 3}, 1 - n, 1 + {B\over 2}, {3\over 2} + {B\over
       2}}\endmatrix ; {\displaystyle 1}\right ]  
\\+ 
  {{\left( 1 + B \right)  \left( 1 + 2 n + B \right)  }
\over     {  \left( 1 + 4 n + 3 B \right) \left( 1 + 2 n + 3 B \right)}}
      {} _{5} F _{4} \!\left [ \matrix { 1 - {{2 n}\over 3}, -{{n}\over 2},
       {1\over 2} - {n\over 2}, -2 n - B, 2 n + 2 B}\\ { -{{2 n}\over 3},
       1 - n, {1\over 2} + {B\over 2}, 1 + {B\over 2}}\endmatrix ;
       {\displaystyle 1}\right ] .\hskip-8pt
\endmultline$$
Each of the $_5F_4$-series can be summed by means of Corollary~A5.
After little manipulation we arrive at (A.4).\quad \quad \qed
\enddemo

The next Lemmas provide the means for finding degree bounds in Step~5
of the proof of Theorem~2 and Step~2 of the proof of Proposition~5.

As usual, given nonnegative integers $n$ and $k$, 
we write 
$$e_k(x_1,\dots,x_n):=\sum _{1\le i_1<\dots< i_k\le n}
^{}x_{i_1}\cdots x_{i_k}$$ 
for the {\it elementary symmetric
function\/} of order $k$ in $x_1,\dots,x_n$. In particular, this
definition implies $e_k(x_1,x_2,\dots,x_n)\equiv 0$ if $n<k$, since then
the defining sum is empty. The following Lemma (together with its
proof) holds with this understanding of the definition of elementary
symmetric functions.
\proclaim{Lemma A7}Let $a$ and $n$ be fixed integers, $n\ge0$. Then,
as $k$ varies through the nonnegative integers, 
$e_k(a,a+1,\dots,a+n-1)$ is a polynomial in $n$ of degree $2k$.
\endproclaim
\demo{Proof}By induction on $k$. The assertion is trivially true for
$k=0$. If we assume that the assertion is true for $k$, we have
$$\align e_{k+1}(a,a+1,\dots,a+n-1)&=\sum _{a\le i_1<\dots<i_{k+1}\le
a+n-1}
^{}i_1i_2\cdots i_{k+1}\\
&=\sum _{0\le i_{k+1}\le n-1} ^{}(a+i_{k+1})e_k(a,a+1,\dots,a+i_{k+1}-1).\tag
A.5
\endalign$$
By induction hypothesis, $e_k(a,a+1,\dots,a+i_{k+1}-1)$ is some polynomial
$p(i_{k+1})$ in $i_{k+1}$ of degree $2k$. Therefore
$e_{k+1}(a,a+1,\dots,a+n-1)$, which by (A.5) is the indefinite sum of a
polynomial of degree $2k+1$, is a polynomial of degree $2k+2$ (cf\.
e.g\. \cite{\JordAA, sec.~32, Example on p.~103}).\quad \quad \qed
\enddemo
\proclaim{Lemma A8}Let $a$ and $p$ be fixed integers, $p\ge0$. Then,
as $j$ varies through the nonnegative integers, 
the coefficient of $x^{j-p}$ in $(x+a)_j$ is a polynomial in $j$ of
degree $2p$.
\endproclaim
\demo{Proof}By definition of shifted factorials we have
$$\align (x+a)_j&=\prod _{i=0} ^{j-1}(x+a+i)\\
&=\sum _{p\ge0} ^{}x^{j-p}\,e_p(a,a+1,\dots,a+j-1).
\endalign$$ 
(Note that we need not give an upper bound for the sum, since
$e_p(a,a+1,\dots,a+j-1)\equiv0$ for $p>j$, see the
paragraph before Lemma~A7 that contains the definition of elementary
symmetric functions).
Therefore, by Lemma~A7, the coefficient of $x^{j-p}$ in $(x+a)_j$ is
a polynomial in $j$ of degree $2p$.\quad \quad \qed
\enddemo
\proclaim{Lemma A9}Let $n$, $m$, and $p$ be fixed integers $0\le p\le
2n+\fl{m/2}-3$. Then, as $j$ varies through the integers, $1\le
j\le n-1$, the coefficient of $x^p$ in
$$(2x+m+1)_j\,(x-j+2)_{\fl{m/2}+j-1}\,(x+m+2j+1)_{2n-2j-2}$$
is $2^j$ times a polynomial in $j$ of degree $\le
2(2n+\fl{m/2}-3-p)$.
\endproclaim
\demo{Proof} We have
$$(2x+m+1)_j=\sum _{r\ge0} ^{}x^{j-r}2^{j-r}b_r(j),$$
where, by Lemma~A8, $b_r(j)$ is a polynomial in $j$ of degree $2r$. Similarly, we
have
$$\align
(x-j+2)_{\fl{m/2}+j-1}&=(-1)^{\fl{m/2}+j-1} 
(-x-\fl{m/2})_{\fl{m/2}+j-1}\\
&=(-1)^{\fl{m/2}+j-1} \sum _{s\ge0}
^{}x^{\fl{m/2}+j-1-s}(-1)^{\fl{m/2}+j-1-s}c_s(j)\\
&= \sum _{s\ge0}
^{}x^{\fl{m/2}+j-1-s}(-1)^{s}c_s(j),
\endalign$$
where, by Lemma~A8, $c_s(j)$ is a polynomial in $(\fl{m/2}+j-1)$ of degree
$2s$, and as such is a polynomial in $j$ of degree $2s$.
Finally, we
have
$$\align (x+m+2j+1)_{2n-2j-2}&=(-x-m-2n+2)_{2n-2j-2}\\
&=\sum _{t\ge0} ^{}x^{2n-2j-2-t}(-1)^td_t(j),
\endalign$$
where, by Lemma~A8, $d_t(j)$ is a polynomial in $2n-2j-2$ of degree
$2t$, and as such is a polynomial in $j$ of degree $2t$.

Putting things together we get
$$\multline
(2x+m+1)_j\,(x-j+2)_{\fl{m/2}+j-1}\,(x+m+2j+1)_{2n-2j-2}\\
=\sum _{k\ge0} ^{}x^{2n+\fl{m/2}-3-k}\,2^j\underset r+s+t=k\to{\sum _{r,s,t\ge0}
^{}}(-1)^{s+t}2^{-r}b_r(j)\,c_s(j)\,d_t(j).
\endmultline$$
Now, the (finite) range of the inner sum does not depend on $j$.
Hence, by what we know about $b_r(j)$, $c_s(j)$, and $d_t(j)$, the inner
sum is a polynomial in $j$ of degree at most
$2r+2s+2t=2k$. By replacing $k$ by $2n+\fl{m/2}-3-p$ we get the
assertion of the Lemma.\quad \quad \qed
\enddemo

\proclaim{Lemma A10}Let $n$, $m$, $p$, and $q$ be fixed integers $0\le p\le
2n-4$, $0\le q\le 2n-3$, $p+q\le 2n-3$. 
Then, as $j$ varies through the integers, $1\le
j\le n-1$, the coefficient of $x^pz^q$ in
$$(2x+m+z+1)_{j-1}\,(x+2z-j+2)_{j-1}\,
(x+m+2j-z+2)_{2n-2j-2}\,(m+3j-3z)
\tag A.6$$
is $2^j$ times a polynomial in $j$ of degree $\le
2(2n-3-p-q)+q-1$.
\endproclaim
\demo{Proof} We have
$$\align (2x+m+z+1)_{j-1}&=\sum _{r\ge0} ^{}(2x+z)^{j-1-r}b_r(j)\\
&=\sum _{r\ge0} ^{}\sum _{h\ge0} ^{}\binom{j-1-r}h
2^{j-1-r-h}x^{j-1-r-h}z^{h}b_r(j),
\endalign$$
where, by Lemma~A8, $b_r(j)$ is a polynomial in $j-1$ of degree $2r$,
and as such a polynomial in $j$ of degree $2r$. 
Similarly, we have
$$\align (x+2z-j+2)_{j-1}&=\sum _{s\ge0} ^{}(x+2z)^{j-1-s}c_s(j)\\
&=\sum _{s\ge0} ^{}\sum _{k\ge0} ^{}\binom{j-1-s}k
x^{j-1-s-k}2^{k}z^{k}c_s(j),
\endalign$$
where, by Lemma~A8, $c_s(j)$ is a polynomial in $j$ of degree $2s$.
Finally, we have
$$\align (x+m+2j-z+&2)_{2n-2j-2}=(-x+z-m-2n+1)_{2n-2j-2}\\
&=\sum _{t\ge0} ^{}(-x+z)^{2n-2j-2-t}d_t(j)\\
&=\sum _{t\ge0} ^{}\sum _{l\ge0}
^{}\binom{2n-2j-2-t}l(-1)^{t+l}
x^{2n-2j-2-t-l}z^{l}d_t(j),
\endalign$$
where, by Lemma~A8, $d_t(j)$ is a polynomial in $2n-2j-2$ of degree
$2t$, and as such is a polynomial in $j$ of degree $2t$.

Putting things together we get
$$\multline (2x+m+z+1)_{j-1}\,(x+2z-j+2)_{j-1}\,
(x+m+2j-z+2)_{2n-2j-2}\,(m+3j-3z)\\
=(m+3j-3z)\sum _{p,q\ge0} ^{}x^pz^q\,2^j
\underset h+k+l=q\to{\sum _{r+s+t=2n-4-p-q} ^{}}
(-1)^{t+l}\,2^{-1-r-h+k}\\
\times\binom{j-1-r}h \binom{j-1-s}k \binom{2n-2j-2-t}l 
b_r(j)c_s(j)d_t(j).
\endmultline$$
Again, the (finite) range of the inner sum does not depend on $j$.
Hence, by what we know about $b_r(j)$, $c_s(j)$, and $d_t(j)$, and
since a binomial $\binom{j+a}\nu$ is a polynomial in $j$ of degree
$\nu$, the inner
sum is a polynomial in $j$ of degree at most
$2r+2s+2t+h+k+l=2(2n-4-p-q)+q$. Combining with $(m+3j-3z)$, we see that the
coefficient of $x^pz^q$ in (A.6) is a polynomial in $j$ of degree at
most $2(2n-3-p-q)+q-1$, as desired.\quad \quad \qed

\enddemo

\enddocument